\documentclass{amsart}
\usepackage{fullpage}
\usepackage{amsmath}
\usepackage{color}
\usepackage{graphicx}
\usepackage{multirow}
\usepackage[english]{babel}
\usepackage{tabulary}
\usepackage{amssymb}
\usepackage{epsfig}
\newtheorem{Proposition}{Proposition}[section]
 
 \newtheorem{Corollary}[Proposition]{Corollary}
 \newtheorem{Lemma}[Proposition]{Lemma}
 
 \newtheorem{Theorem}{Theorem}[section]

\def\blackslug{\hbox{\hskip 1pt \vrule width 4pt height 8pt depth 1.5pt
\hskip 1pt}}
\def\qed{\quad\blackslug\lower 8.5pt\null\par}
 
\def\CC{\mathbb{C}}
 \def\RR{\mathbb{R}}
 \def\NN{\mathbb{N}}
\def\ZZ{\mathbb{Z}}

\def\Re{\mathrm{Re}}
\def\Im{\mathrm{Im}}

\newcommand{\mb}{\mathbf}

\title{A proof for the mode stability of a self-similar wave map}

\author{O. Costin}
\address{Mathematics Department, The Ohio State University, Columbus OH 43210, USA}
\email{costin.9@osu.edu}
\author{R. Donninger}
\address{Rheinische Friedrich-Wilhelms-Universit\"at Bonn,
Mathematisches Institut, Endenicher Allee 60, D-53115 Bonn, Germany}
\email{donninge@math.uni-bonn.de}
\thanks{Roland Donninger is supported by the Alexander von 
Humboldt Foundation via a Sofja Kovalevskaja Award endowed by the German
Federal Ministry of Education and Research.} 
\author{X. Xia}
\address{Mathematics Department, The Ohio State University, Columbus OH 43210, USA}
\email{xia.48@osu.edu}

\begin{document}

\begin{abstract}
We study the fundamental self-similar solution to the SU(2) sigma model, found
by Shatah and Turok-Spergel. We give a rigorous proof for its mode stability.
Based on earlier results by the second author,
the present paper constitutes the last building block for a completely rigorous proof
of the nonlinear stability of the Shatah-Turok-Spergel wave map.
\end{abstract}
\maketitle

\section{Introduction}
\noindent A map $u: \RR^{1,d}\to M$ from $(d+1)$-dimensional
Minkowski space to a Riemannian manifold $M$ is called a 
\emph{wave map} if it satisfies
\[ \partial_\mu \partial^\mu u^a+\Gamma^a_{bc}(u)\partial_\mu u^b \partial^\mu u^c=0 \]
where $\Gamma^a_{bc}$ are the Christoffel symbols on $M$.
Wave maps are natural geometric nonlinear generalizations of the wave equation and they attracted
a lot of interest in the recent past. 
We refer the reader to the survey article \cite{Kri08} and the references in \cite{Don11, DonSchAic11}
for more background.

In the special case $d=3$ and $M=\mathbb{S}^3$, the three-sphere, the wave maps equation reduces
to the well-known SU(2) sigma model from particle physics \cite{GelLev60, Mis78}.
In hyperspherical coordinates on $\mathbb{S}^3$ and for so-called \emph{co-rotational maps} of the form 
$u(t,r,\theta,\varphi)=(\psi(t,r),\theta,\varphi)$, where $(t,r,\theta,\varphi)$ are the standard
spherical coordinates on Minkowski space, the SU(2) sigma model is described by the single semilinear
wave equation
\begin{equation}
\label{eq:wm}
\psi_{tt}-\psi_{rr}-\frac{2}{r}\psi_r
+\frac{\sin(2\psi)}{r^2}=0. 
\end{equation}
Eq.\ \eqref{eq:wm} is a supercritical wave equation and it exhibits finite-time blowup.
This was demonstrated by Shatah \cite{Sha88, CazShaTah98} 
who proved the existence of a self-similar solution
to Eq.\ \eqref{eq:wm} which was later found in closed form by Turok and Spergel \cite{TurSpe90}, 
\[ \psi^T(t,r)=2\arctan\left (\frac{r}{T-t}\right ), \]
where $T$ is a free parameter (the blowup time).
The relevance of $\psi^T$ for generic time evolutions depends on its stability. 
Numerical investigations by Bizo\'n, Chmaj and Tabor \cite{BizChmTab00} 
suggest that $\psi^T$ describes
a stable blowup regime.

\subsection{Mode stability}
A first, heuristic stability analysis (which was already performed in \cite{BizChmTab00})
addresses the question of \emph{mode stability} and proceeds as follows.
One introduces a new coordinate system $(\tau,\rho)$ which is adapted to self-similarity and given
by
\[ \tau=-\log(T-t),\qquad \rho=\frac{r}{T-t}. \]
It is natural to restrict the analysis to the backward lightcone of the blowup point $(T,0)$
which yields the coordinate domain $\tau\geq -\log T$ and $\rho \in [0,1]$ (note that
$t\to T-$ corresponds to $\tau\to\infty$).
Eq.\ \eqref{eq:wm} transforms into
\begin{equation}
\label{eq:wmcss}
\phi_{\tau\tau}+\phi_\tau+2\rho \phi_{\tau\rho}
-(1-\rho^2)\phi_{\rho\rho}-2\frac{1-\rho^2}{\rho}\phi_\rho
+\frac{\sin(2\phi)}{\rho^2}=0
\end{equation}
where $\psi(t,r)=\phi(\tau,\rho)$.
The point of this transformation is of course that the self-similar Shatah-Turok-Spergel
solution $\psi^T$ becomes static in the new coordinates and is simply given by $2\arctan \rho$.
In order to obtain information on its stability one inserts the mode ansatz
\[ \phi(\tau,\rho)=2 \arctan(\rho)+e^{\lambda \tau}u_\lambda(\rho),\qquad \lambda\in \CC \]
and linearizes 
in $u_\lambda$ which yields the ODE spectral problem
\begin{equation}
\label{eq:ODEspec}
-(1-\rho^2)u_\lambda''-2\frac{1-\rho^2}{\rho}u_\lambda'+2\lambda \rho u_\lambda'
+\lambda(\lambda+1)u_\lambda+\frac{V(\rho)}{\rho^2}u_\lambda=0
\end{equation}
where
\begin{equation}
\label{eq:ODEspecV} V(\rho)=2\cos(4\arctan \rho)=\frac{2(1-6\rho^2+\rho^4)}{(1+\rho^2)^2}. 
\end{equation}
Eq.~\eqref{eq:ODEspec} has the two regular singular points $0$ and $1$ and we call a (nontrivial)
solution $u_\lambda$ of Eq.~\eqref{eq:ODEspec} \emph{admissible} if $u_\lambda \in C^\infty[0,1]$.
It follows from the functional 
framework developed in \cite{DonSchAic11} that one may restrict oneself to smooth functions
here.
Admissible solutions of Eq.\ \eqref{eq:ODEspec} with $\Re \lambda\geq 0$ are called
\emph{unstable modes}.
Obviously, the existence of unstable modes is expected to indicate instabilities of $\psi^T$.
In fact, it is easy to see that there exists an unstable mode 
$u_\lambda(\rho)=2\rho \frac{d}{d\rho}\arctan\rho=\frac{2\rho}{1+\rho^2}$ with $\lambda=1$.
However, this instability is a coordinate effect (a ``symmetry mode'')
and it is related to the freedom of choosing
the blowup time $T$ (see \cite{DonSchAic11} for a thorough discussion on this). 
Consequently, the solution $\psi^T$ is called \emph{mode stable} if there do not exist
unstable modes except for the symmetry mode.

In and by itself, the study of mode solutions does not yield any information
on the \emph{nonlinear} stability of $\psi^T$.
However, the second author proved that {\bf mode stability of $\psi^T$ implies nonlinear
stability} \cite{Don11}.

\subsection{Known results}
There are convincing numerical studies that exclude the existence of unstable
modes \cite{BizChmTab00, BIZON, DonAic10}.
Unfortunately, it seems to be very challenging to prove the mode stability of $\psi^T$ rigorously.
The spectral problem \eqref{eq:ODEspec}
is nonstandard since $\lambda$ appears in the coefficient of $u'_\lambda$. 
Thus, in principle it is possible that
there exist unstable modes for nonreal $\lambda$ which complicates matters tremendously.
Of course, one may remove the first derivative by transforming Eq.~\eqref{eq:ODEspec}
into Liouville normal form.
This transformation, however, involves the spectral parameter and it turns out
that it is only admissible if $\Re\lambda>1$ \cite{DonAic08, DonAic09, DonSchAic11}.
It is thus easy to exclude unstable modes with $\Re\lambda>1$ by Sturm oscillation theory
\cite{DonAic09}
but the domain $0\leq \Re\lambda \leq 1$ remains open.
In \cite{DonAic08} it is proved that there do not exist unstable modes for $\lambda \in (0,1)$.
Furthermore, in \cite{DonSchAic11} the existence of unstable modes for 
$\Re\lambda\geq \frac12$
and $\lambda\not=1$ is excluded.
Finally, also in \cite{DonSchAic11} there is given an argument why there cannot exist unstable
modes for $\lambda$ far away from the real axis.

\subsection{Content of the paper}
As a first step, we quantify the constants in the argument from \cite{DonSchAic11}
to show that there are no 
unstable eigenvalues $\lambda$ with $|\Im \lambda|\geq 380$.
The remaining compact region
$\{\lambda\in \CC: 0\leq \Re\lambda\leq \frac12, |\Im \lambda|\leq 380\}$
is then studied by a new approach which is in some sense ``brute force''.
This means that we work with suitable approximations to solutions of the equation,
e.g.~truncated Chebyshev expansions, and a rigorous version of a method used by Bizo\'n
\cite{BIZON} which relates the existence of eigenvalues to the convergence of a certain
continued fraction. 
We provide rigorous error estimates for all our approximations.
Furthermore, we emphasize that all computations are free of rounding errors since
we work exclusively with rational numbers.
As a consequence, our result is completely rigorous although we have to admit that it seems
unrealistic to perform all our computations without the aid of some computer algebra
system.

In the present paper we prove the following result.

\begin{Theorem}
\label{thm:mainmain}
The Shatah-Turok-Spergel wave map is mode stable, i.e., there does not exist a nontrivial
solution $u_\lambda \in C^\infty[0,1]$ to Eq.~\eqref{eq:ODEspec} if $\Re\lambda\geq 0$
and $\lambda\not=1$.
\end{Theorem}

We reiterate that Theorem \ref{thm:mainmain} in conjunction with the result in \cite{Don11}
provides a completely rigorous proof of the \emph{nonlinear} stability 
of $\psi^T$ in the sense of \cite{Don11}.

In view of our method of proof it is natural to split the problem as follows.

\begin{Theorem}\label{MainTh3}
There are no unstable modes for $\lambda$ with $\Re\lambda \in [0,\frac12]$
and $|\Im \lambda|\geq 380$.
\end{Theorem}

\begin{Theorem}\label{MainTh}
There are no unstable modes for $\lambda$ with $\Re\lambda \in [0,\frac12]$
and $|\Im \lambda|\leq 10$.
\end{Theorem}

\begin{Theorem}\label{MainTh2}
	There are no unstable modes for $\lambda$ with $\Re\lambda \in [0,\frac{1}{2}]$
	and $10\leq |\Im \lambda|\leq 380$.
\end{Theorem}

\section{Proof of Theorem \ref{MainTh3}}

\subsection{Functional setup}
We use the functional framework from \cite{DonSchAic11}.
Recall the following spaces used in \cite{DonSchAic11}.
Let
\begin{align*} 
\tilde H_1&:=\{u\in C^2[0,1]: u(0)=u'(0)=0\} \\
\tilde H_2&:=\{u \in C^1[0,1]: u(0)=0\} 
\end{align*}
and
\[ \|u\|_1^2:=\int_0^1 |u'(\rho)|^2\frac{d\rho}{\rho^2},\qquad
\|u\|_2^2:=\int_0^1 |u'(\rho)|^2 d\rho. \]
We denote by $H_j$ the completion of $\tilde H_j$ with respect to $\|\cdot\|_j$, $j\in \{1,2\}$. 
We define two linear operators (one unbounded, the other one bounded) acting on 
$H:=H_1\times H_2$,
\begin{align*} 
\mb L_0 \mb u(\rho)&:=\left (\begin{array}{c}
-\rho u_1'(\rho)+u_1(\rho)+\rho u_2'(\rho)-u_2(\rho) \\
\frac{1}{\rho}u_1'(\rho)-\rho u_2'(\rho) \end{array} \right ) \\
\mb L' \mb u(\rho)&:=\left (\begin{array}{c}
-V_1(\rho)\int_0^\rho s u_2(s)ds \\ 0 \end{array} \right ),
\end{align*}
where $V_1(\rho)=-\frac{16}{(1+\rho^2)^2}$.
The precise domain of $\mb L_0$ is given in \cite{DonSchAic11} but irrelevant for our
purposes. We only note that $\mb L_0$ is a closed operator.
We also set $\mb L:=\mb L_0+\mb L'$.
A straightforward computation shows that $(\lambda-\mb L)\mb u=\mb 0$ is equivalent to
\[ -(1-\rho^2)u''-2\frac{1-\rho^2}{\rho}u'+2\lambda \rho u'+\lambda(\lambda+1)u+\frac{V(\rho)}{\rho^2}u=0 \]
where 
\begin{align*}
u_1(\rho)&=\rho^2 u_2(\rho)+(\lambda-2)\int_0^\rho s u_2(s)ds \\ 
u(\rho)&=\frac{1}{\rho^2}\int_0^\rho su_2(s)ds 
\end{align*}
and
\[ V(\rho)=\rho^2 V_1(\rho)+2=\frac{2(1-6\rho^2+\rho^4)}{(1+\rho^2)^2}. \]
Thus, our goal is to show that $\lambda$ is not an eigenvalue of $\mb L$
if $\Re\lambda\geq 0$
and $|\Im \lambda|\geq 380$.

\subsection{Resolvent bounds}
In order to show absence of eigenvalues, we construct the resolvent.
By the Birman-Schwinger principle we have
\[ \mb R_{\mb L}(\lambda)=\mb R_{\mb L_0}(\lambda)[\mb I-\mb L'\mb R_{\mb L_0}(\lambda)]^{-1}. \]
Thus, we need to estimate $\|\mb L' \mb R_{\mb L_0}(\lambda)\|$ where 
\[ \|\mb u\|^2:=\|u_1\|_1^2+\|u_2\|_2^2. \]
If $\lambda$ is such that $\|\mb L'\mb R_{\mb L_0}(\lambda)\|<1$ then 
$[\mb I-\mb L'\mb R_{\mb L_0}(\lambda)]^{-1}$ exists (Neumann series) and thus,
$\lambda$ is not an eigenvalue.
By recalling the definition of $\mb L'$ this boils down to estimating
\[ \|V_1 K [\mb R_{\mb L_0}(\lambda)\mb f]_2\|_1 \]
in terms of $\mb f$, where $K u(\rho):=\int_0^\rho s u(s)ds$.
To this end, we use Hardy's inequality,
\[ \left \|\frac{u}{|\cdot|}\right \|_{L^2}\leq 2 \|u'\|_{L^2}, 
\qquad u \in H_1,\]
and the estimate $\|u\|_{L^2}\leq \frac{1}{\sqrt 2}\|u'\|_{L^2}$ for $u\in H_2$, which
is a consequence of Cauchy-Schwarz (all function spaces are on the interval $(0,1)$, i.e.,
$L^2=L^2(0,1)$).
Furthermore, we use the bound $\|\mb R_{\mb L_0}(\lambda)\|\leq \frac{1}{\Re \lambda+\frac12}$
which follows from semigroup theory \cite{DonSchAic11}.
First, we estimate the operator norm of $V_1$, viewed as a map from $H_1$ to $H_1$.
We have
\begin{align*}
\|V_1 u\|_1&=\left \|\frac{(V_1 u)'}{|\cdot|}\right \|_{L^2}
\leq \left \|\frac{V_1 u'}{|\cdot|}\right \|_{L^2}+\left \|\frac{V_1' u}{|\cdot|}\right \|_{L^2} \\
&\leq \|V_1\|_{L^\infty}\|u\|_1+2\|V_1'\|_{L^\infty}\|u'\|_{L^2} \\
&\leq 50 \|u\|_1.
\end{align*}
It remains to estimate $\|K[\mb R_{\mb L_0}(\lambda)\mb f]_2\|_1$.
If we set $\mb u:=\mb R_{\mb L_0}(\lambda)\mb f$, we obtain $(\lambda-\mb L_0)\mb u=\mb f$ which
implies
\[ u_1(\rho)=\rho^2 u_2(\rho)+(\lambda-2)Ku_2(\rho)-Kf_2(\rho) \]
as a straightforward computation shows.
Consequently, we find
\begin{align*} |\lambda-2|\|Ku_2\|_1&\leq \|u_1\|_1+\||\cdot|^2 u_2\|_1+\|Kf_2\|_1 \\
&=\|u_1\|_1+\||\cdot|^{-1}(|\cdot|^2 u_2)'\|_{L^2}+\||\cdot|^{-1}(Kf_2)'\|_{L^2} \\
&\leq \|u_1\|_1+2\|u_2\|_{L^2}+\|u_2'\|_{L^2}+\|f_2\|_{L^2} \\
&\leq \|u_1\|_1+(\tfrac{2}{\sqrt 2}+1)\|u_2\|_2+\tfrac{1}{\sqrt 2}\|f_2\|_2 \\
&\leq (2+\sqrt 2)\|\mb u\|+\tfrac{1}{\sqrt 2}\|\mb f\|
\end{align*}
and this yields
\begin{align*}
|\lambda-2|\|K[\mb R_{\mb L_0}(\lambda)\mb f]_2\|_1\leq
\frac{2+\sqrt 2}{\Re\lambda+\frac12}\|\mb f\|+\tfrac{1}{\sqrt 2}\|\mb f\|
< 7.6 \|\mb f\|
\end{align*}
provided $\Re\lambda \geq 0$.
In summary, we find
\[ \|\mb L'\mb R_{\mb L_0}(\lambda)\|< \frac{50\cdot 7.6}{|\lambda-2|}=\frac{380}{|\lambda-2|} \]
and this shows that there are no eigenvalues with 
$\Re\lambda \geq 0$ and $|\Im\lambda| \geq 380$.

\section{Proof of Theorem\,\ref{MainTh}}

\subsection{Preparatory remarks}

A number $\lambda$ for which equation \eqref{eq:ODEspec}, \eqref{eq:ODEspecV} has a solution $u_\lambda \in C^\infty[0,1]$ will be simply called an {\em eigenvalue} of the equation on $[0,1]$.

Note that $\rho=0$ is a regular singular point with indices $1$ and $-2$, and $\rho=1$ is a regular singular point with indices $0$ and $1-\lambda$. For $\lambda$ noninteger, by Fuchsian theory, a solution is $ C^\infty[0,1]$ if and only if it is analytic on $[0,1]$. 

It is convenient to substitute
\begin{equation}\label{chvar}
F(t) = u_\lambda(\sqrt{t})/\sqrt{t}, \quad  G(t) = F(1-t)
\end{equation}
(where $\sqrt{t}>0$ for $t>0$) which transforms \eqref{eq:ODEspec}, \eqref{eq:ODEspecV} to
\begin{equation}
\label{eq:eqG}
t(1-t)G''(t)+\left[ -\frac{5}{2} t+\lambda(1-t) \right] G'(t)+\left[-\frac{\lambda^2+3\lambda}{4}+\frac{1}{2}+\frac{t(4-t)}{(2-t)^2}\right] G(t)=0
\end{equation}

 {\bf Remark.}  {\em A number $\lambda$ with $0 \leq \Re(\lambda) \leq 1/2 $ is an eigenvalue of \eqref{eq:ODEspec}, \eqref{eq:ODEspecV} on $[0,1]$ if and only if $\lambda$ is an eigenvalue of \eqref{eq:eqG} on $[0,1]$.}

 Indeed, it is easy to check using Taylor series at $0$ that a solution $u_\lambda$ analytic at $\rho=0$ has the form $\rho F(\rho^2)$ with $F$ analytic at $0$, hence $u_\lambda$ is analytic on $[0,1]$ if and only if $F$, and therefore $G$,  are analytic on $[0,1]$.
 
 \

 Denote
$$R = \left\{\lambda\in\CC\, |\, 0\leq \Re\lambda \leq \tfrac12,\  -10 \leq \Im\lambda \leq 10\right\}$$
The proof is slightly different on each of the three subsets $R=R_1\cup R_2\cup R_3$ where
\begin{equation}\label{defR1}
R_1=\{\lambda\,|\, 0 \leq \Re\lambda \leq \tfrac12,\   | \Im\lambda| \leq \tfrac12\}
\end{equation}
\begin{equation}\label{defS2}
R_2 = \{\lambda\,|\,0 \leq \Re\lambda \leq \tfrac12,\   \tfrac12 \leq |\Im\lambda| \leq 4\}
\end{equation}
\begin{equation}\label{defS3}
 R_3 = \{\lambda\, |\,0 \leq \Re\lambda \leq \tfrac12,\   4 \leq |\Im\lambda| \leq 10\}
 \end{equation}

 Note that $\lambda$ is an eigenvalue of \eqref{eq:eqG} if and only if 
 $\bar{\lambda}$ is an eigenvalue. 
 Therefore it suffices to consider $\lambda$ with $\Im\lambda \geq 0$. Denote the upper half of $R_i$ by $S_i$, ($i=1,2,3$).

\subsection{Absence of eigenvalues in $S_1$}
\label{subsec:S1}

Let 
$$S_1 = \{\lambda\in R_1\,|\,  \Im\lambda\geq0\}$$ be the upper half of the rectangle \eqref{defR1}. In this section we show that no $\lambda$ in $S_1$ is an eigenvalue. 

\

\noindent \textbf{Outline of the proof}\\
$t=0$ and $t=1$ are regular singularities of eq. \eqref{eq:eqG} with Frobenius indices $\{ 0,1-\lambda \}$ at $t=0$ and $\{ -3/2, 0\}$ at $t=1$. For a fixed eigenvalue $\lambda$, if $G(t)$ is a solution analytic at $t=0$ and $H(t)$ is a solution analytic at $t=1$, then $G(t)$ and $H(t)$ must be linearly dependent, in other words, their Wronskian $W(G,H) \equiv 0$. Thus it suffices to show that $W(G,H)(1/2) \neq 0$ if $\lambda \in S_1$. 

We will prove this by first constructing approximations $G_a(t)$ of $G(t)$ and $H_a(t)$ of $H(t)$, showing that the Wronskian $W(G_a,H_a)$ of the approximations is nonzero at $t=1/2$ and that the difference $| W(G_a,H_a)- W(G,H) | (1/2)$ is smaller than $|W(G_a,H_a)|$. The differences between actual solutions and the chosen approximations are proved rigorously to be small.

\

\noindent \textbf{Proof}
(I) \textit{Approximation.} 

The way the approximation is obtained is irrelevant for the proof. However we describe the method to motivate the approach.

We would like to construct an approximation of $G(t)$ on $[0,1/4]$ and on $(1/4,1/2]$, and an approximation of $H(t)$ on $[3/4,1]$ and on $[1/2,3/4)$. We explain below the process of constructing an approximation of $G(t)$ on interval $[0,1/4]$.

Let $G(t)$ be the solution which is analytic at $t=0$ satisfying $G(0) = \lambda \neq 0$, and $H(t)$ be the solution which is analytic at $t=1$ satisfying $H(1) = 1$. They have series expansions
\begin{equation}
\label{expansionGH}
G(t) = \sum_{n=0}^{\infty} c_n t^n \quad \quad H(t) = \sum_{n=0}^{\infty} d_n (t-1)^n
\end{equation}
Both series have radii of convergence at least one, since the only singularities in $\mathbb{C}$ of \eqref{eq:eqG} are 0, 1 and 2. It is straight forward to see that each $c_n$ is a rational function of $\lambda$ and analytic in $S_1$, the denominators of $c_n$ vanishes only for $\lambda \in \ZZ^-$, and each $d_n$ is a polynomial in $\lambda$. We obtain an approximation $G_a$ of $G(t)$ as follows.\\

(i) We compute $c_n$ for $n \leq 11$ and differentiate the Taylor polynomial of $G(t)$ twice. We obtain:
\begin{equation}
P_0(t,\lambda)=\sum_{n=0}^{9} (n+2)(n+1) c_{n+2} t^n
\end{equation}

(ii) For each coefficient $(n+2)(n+1)c_{n+2}$, we expand $(n+2)(n+1)c_{n+2}$ as a Taylor series in $\lambda$ at $\lambda_0 = \frac{1+i}{4}$. We then approximate $(n+2)(n+1)c_{n+2}$ by its 6-th order Taylor polynomial $c_{T,n}$. After expanding the polynomial, we have
\begin{equation}
c_{T,n} = \sum_{i=0}^{6} a_{n,i} \lambda^i
\end{equation}
so $P_0(t,\lambda)$ is approximated by
\begin{equation}
P_1(t,\lambda) = \sum_{n=0}^{9} c_{T,n} t^n = \sum_{n=0}^{9} \sum_{i=0}^{6} a_{n,i} \lambda^i t^n = \sum_{i=0}^{6}  \left( \sum_{n=0}^{9} a_{n,i} t^n \right) \lambda^i = \sum_{i=0}^{6}  b_i \lambda^i
\end{equation}

(iii) Next,  we use a Chebyshev polynomial approximation in $t$ on $[0,1/4]$ for $\Re(b_i)$ 
and $\Im(b_i)$, with precision $2 \cdot  10^{-3}$. 
Denote the approximations by   $c^{(r)}_{c,i}$ and $c^{(i)}_{c,i}$, respectively,  $P_1(t,\lambda)$ is approximated by
\begin{equation}
P_2(t,\lambda) = \sum_{i=0}^{6}  (c^{(r)}_{c,i} + i c^{(i)}_{c,i}) \lambda^i = \sum_{i,n} c_{c,n,i} \lambda^i t^n
\end{equation}

(iv) Each $c_{c,n,i}$ is replaced by a rational number $c_{n,i}$ within $e^{-7}$ of it, to allow for rigorous proofs. Finally we have an approximation of $G''(t)$ in the form
\begin{equation}
P_a(t,\lambda) = \sum_{i,n} c_{n,i} \lambda^i t^n
\end{equation}

(v) The approximation $G_a(t)$ of $G(t)$ on $[0,1/4]$ is obtained 
by integrating $P_a$ with respect to $t$ twice and the constants of integration 
are determined by the initial condition: 
\begin{equation}
G_a(0) = G(0) = \lambda, \quad  G'_a(0) = G'(0) = c_1 = \frac{1}{4} \lambda^2+\frac{3}{4} \lambda - \frac{1}{2} 
\end{equation}

To obtain the approximation $G_{a}$ of $G$ on $(1/4, 1/2]$ we work 
with an approximate Taylor polynomial of $G(t)$ at $t=1/4$ obtained by finding a 
series solution of \eqref{eq:eqG} with the initial conditions $G_a(1/4),\, G_a'(1/4)$. 
In order to obtain ``nicer'' expansion coefficients,
we allow for a slight discontinuity of $G_a$ at $t=1/4$.

 Similarly, we obtain a piecewise function $H_a$ as an approximation of $H$ on [1/2,1]. See the appendix \S \ref{subsec:GaHa} for the expressions of $G_a$ and $H_a$. \\

\noindent(II) Estimate of $W(G,H)(1/2)$.

Let $\delta_1(t,\lambda) = G(t,\lambda) - G_a(t,\lambda) $ for $t \in [0,1/2]$, and  $\delta_2(t,\lambda) = H(t,\lambda) - H_a(t,\lambda) $ for $t \in [1/2,1]$. Both $\delta_1$ and $\delta_2$ are analytic at $t=1/2$. In a small neighborhood of $t=1/2$ we have
\begin{equation}
\label{Wro}
W(G,H) = W(G_a,H_a) + G_a \delta'_2+\delta_1 H'_a+\delta_1\delta'_2-H_a \delta'_1-\delta_2 G'_a-\delta_2 \delta'_1
\end{equation}

It is not hard to show that 
\begin{equation}\label{loubou}
|W(G_a,H_a) (1/2)| > 0.46
\end{equation}
(the proof is found in \S\ref{Ploubou}).

We will show that the sum of the absolute values of the other terms on the right hand side of (\ref{Wro}) less than 0.03, implying that $W(G,H)(1/2) $ is not zero.

We detail below the estimates of $\delta_1(1/2)$ and $\delta'_1(1/2)$, obtained in two steps.\\

Denote
\begin{equation}
\epsilon_1 = \mathcal{L} \delta_1=-\mathcal{L} G_a,	\quad t\in(0,1/4)\cup(1/4,1/2)
\end{equation}
where $\mathcal{L}$ is the differential operator in (\ref{eq:eqG}). More precisely,
\begin{equation}\label{oeq}
\mathcal{L} \delta_1= t(1-t)\delta''_1+\left[ -\frac{5}{2} t+\lambda(1-t) \right] \delta'_1+\left[\alpha+g(t)\right] \delta_1 = \epsilon_1
\end{equation}
where $\alpha = -\frac{\lambda^2+3\lambda}{4}+\frac{1}{2}$, $g(t) =\frac{t(4-t)}{(2-t)^2} $ and $(2-t)^2 \epsilon_1$ has an explicit expression as a polynomial. \\

(i) Estimates for $\delta_1(t)$ and $\delta'_1(t)$ on $[0,1/4]$.\\
\indent Note that $G_a (0) = G(0)$, $G'_a (0) = G'(0)$, so $\delta_1(t) = t^2\tilde{ \tilde{\delta_1}}(t)$, where $\tilde{ \tilde{\delta_1}}(t)$ is analytic at $t=0$.
It is also obvious from \eqref{oeq} that $\epsilon_1(t,\lambda) = t \,\tilde{\epsilon}_1(t,\lambda)$ where $(2-t)^2\tilde{\epsilon}_1(t)$ is a polynomial as well. 

\begin{Lemma}\label{estimep} Let $\lambda \in S_1$.

(i) For all $t \in [0,1/4]$, we have the following bounds:
\begin{align}
|\tilde{\epsilon}_1(t)| &\leq e_1 = 0.0015 \\
|\alpha| &\leq A =\sqrt{ \frac{117}{256}}\quad\quad(independent\,\,of\,\, t)\\
 1\leq\frac{4-t}{(2-t)^2}&\leq \frac{60}{49}  
\end{align}

(ii) For $t\in (1/4,1/2]$:
\begin{align}
 \frac{15}{49} \leq \frac{t\,(4-t)}{(2-t)^2}&\leq \frac{7}{9} \\
|\epsilon_1(t)| &\leq e_2 = 2.7 \cdot 10^{-3} 
\end{align}

\end{Lemma}
The proof of \ref{estimep} is given in the Appendix \S\ref{subsec:estipoly}.\\

For $t \in [0,1/4]$, we have
\begin{equation}
\label{eq:preint}
\frac{\mathrm d}{\mathrm d t} \big( (1-t)^{5/2} t^\lambda \delta'_1(t) \big) = (1-t)^{3/2} t^{\lambda} \tilde{\epsilon_1}(t) - (1-t)^{3/2} t^{\lambda+1} \big(\alpha+g(t)\big) \tilde{\tilde{\delta_1}}(t)
\end{equation}
Integrating from $0$ to $t$, we have
\begin{align}
\label{d1p}
\delta'_1(t) &=   (1-t)^{-5/2} t^{-\lambda}  \int_0^t (1-s)^{3/2} s^{\lambda} \tilde{\epsilon}_1(s) \ \mathrm{d}s \nonumber \\
&- (1-t)^{-5/2} t^{-\lambda}  \int_0^t (1-s)^{3/2} s^{\lambda+1}  \big(\alpha+g(s)\big)\tilde{\tilde{\delta}}_1(s) \mathrm{d}s 
\end{align}
One more integration gives
\begin{align}
\label{d1}
\delta_1(t) &= \int_0^t (1-s)^{-5/2} s^{-\lambda}\mathrm{d}s  \int_0^s (1-u)^{3/2} u^{\lambda} \tilde{\epsilon}_1(u) \ \mathrm{d}u \nonumber \\
&-\int_0^t (1-s)^{-5/2} s^{-\lambda}\mathrm{d}s  \int_0^s (1-u)^{3/2} u^{\lambda+1}  \big(\alpha+g(u)\big) \tilde{\tilde{\delta}}_1(u)\ \mathrm{d}u
\end{align}
Let $x=\Re{\lambda}$; \eqref{d1} implies
\begin{align}
\label{ineq1}
|\delta_1(t)| &\leq \int_0^t (1-s)^{-5/2} s^{-x}\mathrm{d}s  \int_0^s (1-u)^{3/2} u^{x} |\tilde{\epsilon}_1(u)| \ \mathrm{d}u \nonumber \\
&+\int_0^t (1-s)^{-5/2} s^{-x}\mathrm{d}s  \int_0^s (1-u)^{3/2} u^{x+1}  \big(|\alpha|+|g(u)|\big) |\tilde{\tilde{\delta}}_1(u)|\ \mathrm{d}u
\end{align}

Denote by $D_1$ an uniform bound of $\tilde{\tilde{\delta}}_1$ for $t\in[0,1/4]$ and $\lambda \in S_1$. From (\ref{ineq1}) we have
\begin{align}
\label{estid1}
t^2 |\tilde{\tilde{\delta}}_1(t)| &\leq e_1\, \int_0^t (1-s)^{-5/2} s^{-x}\mathrm{d}s  \int_0^s (1-u)^{3/2} u^{x}  \ \mathrm{d}u \nonumber \\
&\quad+A \, D_1 \, \int_0^t (1-s)^{-5/2} s^{-x}\mathrm{d}s  \int_0^s (1-u)^{3/2} u^{x+1}  \ \mathrm{d}u \nonumber \\
&\quad+ \frac{60}{49} \, D_1 \, \int_0^t (1-s)^{-5/2} s^{-x}\mathrm{d}s  \int_0^s (1-u)^{3/2} u^{x+2} \ \mathrm{d}u \end{align}
and, using the fact that $\int_0^s(1-u)^{3/2}u^\xi du\leq \frac{s^{\xi+1}}{\xi+1}$ for any $\xi\geq 0$, we obtain further
\begin{align}
\label{estid1c}
&t^2 |\tilde{\tilde{\delta}}_1(t)|\nonumber\\
&\leq \frac{e_1}{x+1} \int_0^t (1-s)^{-5/2} s \,  \mathrm{d}s 
+\frac{A\, D_1}{x+2}  \int_0^t (1-s)^{-5/2} s^{2} \mathrm{d}s \nonumber\\
&\quad+ \frac{60}{49} \, \frac{D_1}{x+3}  \int_0^t (1-s)^{-5/2} s^{3} \mathrm{d}s \nonumber \\
&\leq {e_1}\, \frac{2}{3} \frac{-2+3t +2(1-t)^{3/2}}{(1-t)^{3/2}} 
+\frac{A\, D_1}{2}\, \frac{16}{3} \frac{-1+3t/2 -3t^2/8+(1-t)^{3/2}}{(1-t)^{3/2}} \nonumber \\
&\quad+ \frac{60}{49} \,\frac{D_1}{3}\, \frac{32}{3} \frac{-1+3t/2-3t^2/8-t^3/16 +(1-t)^{3/2}}{(1-t)^{3/2}} \nonumber \\
&\leq {e_1}\,\frac{2}{3} \frac{\sqrt{3} t^2 / 2}{(1-t)^{3/2}} 
+\frac{A\, D_1}{2}\, \frac{16}{3} \frac{t^3/6\sqrt{3}}{(1-t)^{3/2}} + \frac{60}{49} \,\frac{D_1}{3}\, \frac{32}{3} \frac{t^4/12\sqrt{3}}{(1-t)^{3/2}}  \nonumber \\
&\leq \frac{8}{9} e_1 t^2+\frac{32}{81} A D_1 t^3+ \frac{1280}{3969} D_1 t^4  \end{align}
After canceling $t^2$ and taking supremum of $|\tilde{\tilde{\delta}}_1(t)|$ over $[0,1/4]$, we have 
\begin{equation}
D_1 \leq \frac{8}{9} e_1 + \left(\frac{32}{81}\cdot A \cdot \frac{1}{4} +\frac{1280}{3969} \cdot \frac{1}{4^2} \right)\cdot D_1
\end{equation}
Using the bounds of Lemma\,\ref{estimep}, we get
\begin{equation}
D_1 < 0.0015
\end{equation}
and consequently
\begin{equation}
|\delta_1(t) |= t^2 |\tilde{\tilde{\delta}}_1(t)| < \frac{15}{16} \cdot 10^{-4}
\end{equation}
In a similar way,  using \eqref{d1p}  and  we obtain
\begin{equation}
|\delta'_1(t) | < 9 \cdot 10^{-4} \text{ on } [0,1/4)
\end{equation}

\

(ii) Estimate of $G_a$ on $(1/4,1/2]$.\\

\indent For $t\in (1/4,1/2]$, we integrate \eqref{eq:preint} and obtain the following equations:
\begin{align}
\label{d2p}
\delta'_1(t) &= (1-t)^{-5/2} t^{-\lambda} (3/4)^{5/2}\, (1/4)^{\lambda} \delta'_1(1/4+0) \nonumber \\
 &\quad+ (1-t)^{-5/2} t^{-\lambda} \int_{1/4}^t (1-s)^{3/2} s^{\lambda-1} \epsilon_1(s) \mathrm{d}s \nonumber \\
  &\quad-(1-t)^{-5/2} t^{-\lambda} \int_{1/4}^t (1-s)^{3/2} s^{\lambda-1} \left(\alpha+g(s)\right)\,\delta_1(s) \mathrm{d}s
 \end{align}
\begin{align}
\label{d2}
\delta_1(t) &= \delta_1(1/4+0)+\int_{1/4}^t (1-s)^{-5/2} s^{-\lambda} (3/4)^{5/2}\, (1/4)^{\lambda} \delta'_1(1/4+0) \mathrm{d}s \nonumber \\
 &\quad+ \int_{1/4}^t(1-s)^{-5/2} s^{-\lambda} \mathrm{d}s \int_{1/4}^s (1-u)^{3/2} u^{\lambda-1} \epsilon_1(u) \mathrm{d}u \nonumber \\
  &\quad-\int_{1/4}^t (1-s)^{-5/2} s^{-\lambda} \mathrm{d}s \int_{1/4}^s (1-u)^{3/2} u^{\lambda-1} \left(\alpha+g(u)\right)\,\delta_1(u) \mathrm{d}u
\end{align}
Estimating the integrals \eqref{d2p} and \eqref{d2} as in \eqref{estid1} and \eqref{estid1c}, in the end we get
\begin{align}\label{vip1}
&|\delta_1(t) | < 1.6 \cdot 10^{-3} \ \ \  \text{and}\\
&|\delta'_1(t)| <1.3\cdot 10^{-2} \ \ \ \text{for all }t\in (1/4,1/2]
\end{align}

\indent The estimates of $\delta_2(t)$ and $\delta'_2(t)$ on $[3/4,1]$ and then on $[1/2,3/4)$ are obtained similarly. We obtain:
\begin{align}
\label{vip2}
&|\delta_2(1/2) | < 1.3 \cdot 10^{-3} \nonumber\\
&|\delta'_2(1/2)| <8 \cdot 10^{-3}
\end{align}

Also, we have the following estimates for the approximating polynomials using the method in \S\ref{subsec:estipoly}:
\begin{align}
\label{G0H0}
&|G_a(1/2) | < 1.07 \nonumber\\
&|G'_a(1/2) | < 2.63 \nonumber\\
&|H_a(1/2) | < 0.59 \nonumber\\
&|H'_a(1/2) | < 0.74
\end{align}

Finally, using (\ref{vip1}), (\ref{vip2}) and (\ref{G0H0}) in (\ref{Wro}), we get, at $t=1/2$:
\begin{align}
|W(G,H)| &\geq |W(G_a,H_a)| - \left( |G_a \delta'_2|+ |\delta_1 H'_a|+|\delta_1\delta'_2| + |H_a \delta'_1| + |\delta_2 G'_a| + |\delta_2 \delta'_1| \right) \nonumber\\
&> 0.46 - 0.03 = 0.43 > 0
\end{align}

Hence $S_1$ contains no eigenvalues.

\subsection{Absence of eigenvalue in $S_2$}

Next we show that the rectangle \eqref{defS2}, adjacent to $R_1$, does not contain any eigenvalues either.

\begin{Proposition}
\label{PropnS2}
$\lambda$ is not an eigenvalue of (\ref{eq:eqG}) if $\lambda \in S_2= \{z\in\CC\, |\, 0 \leq \Re(z) \leq 1/2, 1/2 \leq \Im(z) \leq 4\}$. 
\end{Proposition}

\noindent \textbf{Proof}\\
The idea is similar to the proof in \S\ref{subsec:S1}: we look for an approximation $G_a(t)$ of a solution $G(t)$ which is real analytic at  $t=0$, and an approximation $H_a(t)$ of a solution $H(t)$ which is real analytic at  $t=1$ then estimate and compare the Wronskians.

Let the solutions have the following expansions:
\begin{equation}
\label{expansionGH2}
G(t) = \sum_{n=0}^{\infty} a_n t^n \quad \quad H(t) = \sum_{n=0}^{\infty} b_n (1-t)^n
\end{equation}
where $G(0)=H(1)=1$.

The coefficients $a_n$ and $b_n$ satisfy the following recurrence relations:
\begin{equation}
\label{rec1}
a_{n+1}=p_n(\lambda)\, a_n - \frac{1}{(n+1)(n+\lambda)}\sum_{k=0}^{n-1}\frac{n-k+1}{2^{n-k}}a_k
\end{equation}

\begin{equation}
\label{rec2}
b_{n+1}=q_n(\lambda) \, b_n + \frac{1}{(n+1)(n+\frac{5}{2})}\sum_{k=0}^{n-1}4(-1)^{n-k+1}(n-k+1)b_k
\end{equation}
where
\begin{equation}
\label{pn}
p_n=\frac{n^2+(\lambda+\frac{3}{2}) n+\frac{\lambda^2+3 \lambda}{4}-\frac{1}{2}}{(n+1)(n+\lambda)}
\end{equation}
\begin{equation}
\label{qn}
q_n=\frac{n^2+(\lambda+\frac{3}{2}) n+\frac{\lambda^2+3 \lambda}{4}-\frac{7}{2}}{(n+1)(n+\frac{5}{2})}
\end{equation}

First we prove the following estimates:

\begin{Lemma}
\label{JnMnLemma}
Assume $\lambda \in S_2$. There are constants $J_n$ and $M_n$ so that

(i) $|a_n|\leq J_n$, for all $n \geq 0$. For $n \geq 30$, $J_n \leq k_1^n$, where $k_1=1.026$.

(ii) $|b_n|\leq M_n$, for all $n \geq 0$. For $n \geq 30$, $M_n \leq k_2^n$, where $k_2=1.2$.\\
The values of $J_n$ and $M_n$ for $0 \leq n \leq 30$ are given in \S\ref{valJnMn}.
\end{Lemma}

The proof is given in \S\ref{JnMnProof}.

Now we obtain estimates of remainders $G(t)-G_a(t)$ and $H(t)-H_a(t)$, where $G_a$ is the type of rational functions discussed in \S\ref{rationalS2} and $H_a$ is a polynomial where we choose:
\begin{equation}
G_a(t)=\sum_{k=0}^{8}a_k t^k,	\quad\quad H_a(t) =\sum_{k=0}^{10}b_k (1-t)^k
\end{equation}
and consider the Wronskians $W(G,H)$ and $W(G_a,H_a)$ at $t_0=0.56$. With the same notations as in the proof for $S_1$, $\delta_1 = G(t)-G_a(t)$ and $\delta_2 = H(t)-H_a(t)$, we have (\ref{Wro}) for $\lambda \in S_2$. Using the method in \S\ref{rationalS2} and \S\ref{polynS2} we have the following estimates :
\begin{equation}
\label{Gat0}
|G_a(t_0)| \leq 2.14
\end{equation}
\begin{equation}
|H_a(t_0)| \leq 0.61
\end{equation}
\begin{equation}
|G'_a(t_0)| \leq 5.83
\end{equation}
\begin{equation}
|H'_a(t_0)| \leq 1.32
\end{equation}
The lower bound of $|W(G_a,H_a)(t_0)|$ on $S_2$ is proven in \S\ref{WaS2}.
\begin{equation}
|W_a(t_0)| =|G_a(t_0)H'_a(t_0)-G'_a(t_0)H_a(t_0)| > 1.06
\end{equation}
\begin{equation}
|\delta_1(t_0)| \leq  \sum_{i=9}^{\infty} |a_i| t_0^i  \leq  \sum_{i=9}^{30} J_i t_0^i + \sum_{i=31}^{\infty} (k_1 t_0)^i <0.014
\end{equation}
\begin{equation}
|\delta'_1(t_0)| \leq  \sum_{i=9}^{\infty}i |a_i| t_0^{i -1} \leq  \sum_{i=9}^{30}i J_i t_0^{i-1} + \sum_{i=31}^{\infty} i k_1^i t_0^{i-1} <0.252
\end{equation}
Similarly:
\begin{equation}
|\delta_2(t_0)| <0.003
\end{equation}
\begin{equation}
\label{d2t0}
|\delta'_2(t_0)| <0.058
\end{equation}
Combining equations \eqref{Gat0} - \eqref{d2t0} we get
\begin{align}
|W(G,H)(t_0)| 
&\geq | W(G_a,H_a)(t_0)| -| G_a \delta'_2|(t_0)-|\delta_1 H'_a|(t_0)-|\delta_1\delta'_2|(t_0)		\nonumber\\
&\quad -|H_a \delta'_1|(t_0)-|\delta_2 G'_a|(t_0)-|\delta_2 \delta'_1|(t_0)		\nonumber\\
&\geq 1.06 - 2.14\cdot0.058-0.014\cdot1.32-0.014\cdot0.073		\nonumber\\
&\quad -0.61\cdot0.252-0.003\cdot5.83-0.003\cdot0.252 > 0.744
\end{align}
Hence no $\lambda \in S_2$ is an eigenvalue. Proposition \ref{PropnS2} is proved.

\subsection{Absence of eigenvalues in $S_3$}\label{S3}

\begin{Proposition}
$\lambda$ is not an eigenvalue of \eqref{eq:ODEspec}, \eqref{eq:ODEspecV} if $\lambda \in S_3 = \{z\in\CC\, |\,0 \leq \Re(z) \leq 1/2, 4 \leq \Im(z) \leq 10\}$. 
\end{Proposition}

\noindent \textbf{Proof}\\
The proof follows the idea of \cite{BIZON}. First we introduce some notation. Make a change of variable
\begin{equation}
y(x) = \left(\frac{2}{1+\rho^2}\right)^{\frac{1-\lambda}{2}}\rho\,u(\rho), \quad\quad x=\frac{2\rho^2}{1+\rho^2}
\end{equation}
and \eqref{eq:ODEspec}, \eqref{eq:ODEspecV} becomes:
\begin{align}
\label{eqbizon}
{x}^{2} \left( 1-x \right)  \left( 2-x \right) {\frac {d^{2}}{d{x}^{2}}}y \left( x \right) &+x \left[ 1- \left( 1+\lambda \right) x \left( 2-x \right)  \right] {\frac {d}{dx}}y \left( x \right) \nonumber\\
&-\frac{1}{4}\,\left[ {\lambda}^{2}x \left( 1-x \right) +9\,{x}^{2}-17\,x+4 \right] y \left( x \right) = 0
\end{align}
We can see from the Frobenius indices of $u(\rho)$ (see the \S 2.1) that the change of variable preserves analyticity of solutions at 0 and 1. So the eigenvalues of $u(\rho)$ are the same as those of $y(x)$.

Let $y$ be the solution of (\ref{eqbizon}) which is analytic at $x=0$ and satisfies $y'(0)=1$. We write 
\begin{equation}
\label{seriesS3}
y(x)=\sum_{n=0}^{\infty}c_n x^n
\end{equation}
If $\lambda$ is an eigenvalue then 
the radius of convergence of the series is strictly larger than 1, therefore, 
it is at least 2. By substituting the series into (\ref{eqbizon}) we find that $c_n$ satisfy the recurrence relation:
\begin{align}\label{rec5}
& p_2(0)c_2+p_1(0)c_1 = 0		\nonumber\\
&p_2(n)c_{n+2}+p_1(n)c_{n+1}+p_0(n)c_n=0,	\quad n=1,2,...
\end{align}
where $c_0=0$ and $c_1=1$ and
\begin{align}\label{coeffps}
& p_2(n) = 8n^2+28n+20,		\nonumber\\
& p_1(n) = -12n^2-(20+8\lambda)n-\lambda^2-8\lambda+9,		\nonumber\\
& p_0(n) = 4n^2+4\lambda n+\lambda^2-9.
\end{align}
Also, let
\begin{equation}
A_n=\frac{p_1(n)}{p_2(n)}, \quad B_n=\frac{p_0(n)}{p_2(n)}, \quad r_n=\frac{c_{n+1}}{c_{n}}.
\end{equation}
Then (\ref{rec5}) is rewritten as
\begin{equation}
r_n = -\frac{B_n}{A_n+r_{n+1}},
\end{equation}

We notice that the recurrence relation (\ref{rec5}) is a second order difference equation of Poincar\'{e} type with characteristic roots $1$ and $\frac{1}{2}$, so by Poincar\'{e}'s theorem (\cite{Elaydi}) either $c_n \equiv 0$ for $n$ large or 
\begin{equation}\label{2lim}
\lim_{n\to \infty} r_n \in\left\{1,\frac{1}{2}\right\} 
\end{equation}

Suppose that $c_n \equiv 0$ for $n$ large for some values of $\lambda$, then there are only finitely many values of $n$ such that $c_n \neq 0$. Assume that $n=N$ is the largest positive integer such that $c_n \neq 0$, then $c_{N+1} =c_{N+2} =0 $. Using recurrence relation \eqref{rec5} and the fact that $p_0(n) \neq 0$ for any values of $\lambda$ with $\Re{\lambda} \in [0,1/2]$ we have $c_N = 0$, contradicting the assumption $c_N \neq 0$. Hence the situation that $c_n \equiv 0$ for $n$ large will not arise for the values of $\lambda$ in concern.

Note that $\lambda$ is an eigenvalue if and only if $\lim_{n\to \infty} r_n = 1/2$ (meaning that the radius convergence of  of \eqref{seriesS3} is 2). In this case $c_n$ is a minimal solution \footnote{A solution $\Phi(n)$ of a second-order difference equation
\begin{equation*}
x(n+2) + p_1(n)\,x(n+1)+p_2(n)\,x(n) =0
\end{equation*}
 is said to be minimal if
 $$\lim_{n \to \infty} \frac{\Phi(n)}{\Psi(n)} = 0$$
 for any other solution $\Psi(n)$ of the equation that is not a multiple of $\Phi(n)$.} (\cite{Elaydi}) of the recurrence (\ref{rec5}). \\

For an arbitrary $\lambda$, let $\psi_n = \phi_{n+1}/\phi_n$ where $(\phi_n)_{n\geq1}$ is a minimal solution of (\ref{rec5}), then $(\psi_n)_{n\geq1}$ is unique. By Pincherle's theorem \cite{Elaydi, Lor92}, 
for each $n\geq 1$:
\begin{equation}\label{eqcf}
\psi_n = -\frac{B_n}{A_n-}\,\frac{B_{n+1}}{A_{n+1}-}\,\frac{B_{n+2}}{A_{n+2}-} \cdot\cdot\cdot
\end{equation}

Therefore we just need to show that for some $n \geq 1 $ the equation
\begin{equation}
\frac{c_{n+1}}{c_n} = -\frac{B_n}{A_n-}\,\frac{B_{n+1}}{A_{n+1}-}\,\frac{B_{n+2}}{A_{n+2}-} \cdot\cdot\cdot
\end{equation}
does not have any roots in $S_3$.

Let $n_0 = 10$. Consider the Banach space of sequences
\begin{equation}
\mathcal{B}:=\{(a_n)_{n \geq n_0} : a_n\in \CC \}
\end{equation}
with sup norm: $\|(a_n)_{n \geq n_0}\|:= \sup_{n \geq n_0} |a_n|$, and the closed ball in $\mathcal{B}$ 
\begin{equation}
S=\{(a_n)_{n \geq n_0} \in \mathcal{B}: \|(a_n)_{n \geq n_0}\| \leq 0.6 \}.
\end{equation}
Define an operator $\mathcal{N}$ in $S$:
\begin{equation}
(\mathcal{N} [(a_n)_{n \geq n_0}] )_n = -\frac{B_n}{A_n+a_{n+1}}
\end{equation}

\begin{Lemma}
(i) The operator $\mathcal{N}$ preserves $S$ and \\
(ii) $\mathcal{N}$ is contractive on $S $ with a contractivity factor less than 0.8.
\end{Lemma}

\noindent \textbf{Proof}\\
(i) $\mathcal{N}$ preserves $S$.\\
Let $(r_n)_{n\geq n_0} \in S$, write
\begin{equation}
(\mathcal{N} [(r_n)_{n \geq n_0}] )_n = -\frac{B_n}{A_n+r_{n+1}}=-\frac{B_n/A_n}{1+r_{n+1}/A_n}
\end{equation}
We obtain estimates of upper bounds of $|B_n/A_n|$ and lower bounds of $|A_n|$ in $S_3$, which are provided in the formula \eqref{Z1}, \eqref{Z2} and \eqref{Z3} in \S\ref{S_3}, and from them we know that for $n\geq n_0$,
\begin{equation}
\Big|-\frac{B_n/A_n}{1+r_{n+1}/A_n}\Big| \leq \frac{|B_n/A_n|}{1-0.6/|A_n|} \leq 0.6
\end{equation}

(ii) $\mathcal{N}$ is contractive on $S $ with a contraction factor less than 0.8.\\
Similarly, let $(r_n)_{n\geq n_0}$,$(s_n)_{n\geq n_0} \in S$, 
\begin{align}
(\mathcal{N} [(r_n)_{n \geq n_0}] )_n & - (\mathcal{N} [(s_n)_{n \geq n_0}] )_n= -\frac{B_n}{A_n+r_{n+1}}+\frac{B_n}{A_n+s_{n+1}}		\nonumber\\
&=\frac{B_n}{(A_n+r_{n+1})(A_n+s_{n+1})}\,(r_{n+1}-s_{n+1})		\nonumber\\
&=\frac{B_n/A^2_n}{(1+r_{n+1}/A_n)(1+s_{n+1}/A_n)}\,(r_{n+1}-s_{n+1})
\end{align}
With estimates of upper bounds of $|B_n/A^2_n|$ (see \eqref{Z3} and \eqref{Z32}), we are able to show that the contraction factor is less than 0.8, so that proves this part of the lemma.

Hence by contractive mapping theorem, $\mathcal{N}^m(\mathbf{0}) \to \mathcal{N}^{\infty}(\mathbf{0})$ as $m \to \infty$, where $\mathcal{N}^{\infty}(\mathbf{0})$ is used to denote the unique fixed point in $S$. On the other hand, by Pincherle's theorem, if $\lambda$ is an eigenvalue then $(r_n)_{n\geq n_0} = (c_{n+1}/c_n)_{n\geq n_0}$ is the fixed point. In particular, eigenvalues $\lambda$ are roots of the equation
\begin{equation}
(\mathcal{N}^{\infty}(\mathbf{0}))_{n_0}=r_{n_0}
\end{equation}
We then observe that over the boundary of $R_3$
\begin{equation}
\label{eqS_34}
\|(\mathcal{N}(\mathbf{0}))-(\mathcal{N}^{2}(\mathbf{0}))\| < 0.106
\end{equation}
\begin{equation}
\label{eqS_35}
|(\mathcal{N}(\mathbf{0}))_{n_0} -r_{n_0}|>0.595
\end{equation}
(see \S\ref{S_34} and \S\ref{S_35} for the proof). So
\begin{align}
\|(\mathcal{N}(\mathbf{0})) -(\mathcal{N}^{\infty}(\mathbf{0}))\| 
& \le \|\mathcal{N}(\mathbf{0})-\mathcal{N}^{2} (\mathbf{0})\| (\mathrm{I} - \|\mathcal{N}\|)^{-1} \nonumber\\
&<|(\mathcal{N}(\mathbf{0}))_{n_0} -r_{n_0}|
\end{align}
Now since over the boundary of $S_3$ we have
\begin{equation}
\begin{split}
&\quad |(\mathcal{N}^{\infty}(\mathbf{0}))_{n_0} -(\mathcal{N}(\mathbf{0}))_{n_0}| \\&= |[ (\mathcal{N}^{\infty}(\mathbf{0}))_{n_0}-r_{n_0} ]- [(\mathcal{N}(\mathbf{0}))_{n_0} -r_{n_0}]|	 \\
&<|(\mathcal{N}(\mathbf{0}))_{n_0} -r_{n_0}|
\end{split}
\end{equation}
By Rouch\'e's theorem, $(\mathcal{N}^{\infty}(\mathbf{0}))_{n_0}-r_{n_0} $ 
and $(\mathcal{N}(\mathbf{0}))_{n_0} -r_{n_0}$ have the same number of roots 
in $S_3$, but the latter is an explicit rational function with no roots 
in $S_3$, hence there are no eigenvalues in $S_3$.

\section{Proof of Theorem \ref{MainTh2}}

\subsection{Introduction of notations}$ $\\
\indent Denote
\begin{equation}\label{S4}
S_4:=\{\lambda \,\,\vert\,\, 0 \leq \Re\lambda \leq\frac{1}{2},\,\,10 \leq \Im\lambda \leq 380\}
\end{equation}
With the notations as in \S\ref{S3}, we rewrite (\ref{rec5}) as
\begin{equation}\label{rec101}
r_{n+1} =-A_n -\frac{B_n}{r_{n}},
\end{equation}
Proving Theorem \ref{MainTh2} is equivalent to proving that for $\lambda \in S_4$,
\begin{equation}
\lim_{n\to\infty} r_n = 1
\end{equation}

The method of obtaining the estimates of $r_n$ is as follows. We first obtain a ``quasi-solution'' $\tilde{r}_n$ of the recurrence \eqref{rec101} (more precisely for $r_n = c_{n+1}/c_{n}$). That is, a function $\tilde{r}_n$ which is close enough to $r_n$. For this purpose the analysis need not be rigorous, since in the sequel we prove that the relative error $|r_n/\tilde{r}_n-1|$ is small. The concrete construction of $\tilde{r}_n$ is as follows.

For $n\le 50$ we  project a suitable branch of $\log{r_n}$ on Chebyshev polynomials in $\Re(\lambda)$ and $\Im(\lambda)$, of degree at most 6. For $n > 50$ we simply define $\tilde{r}_n$ to be one of the exact solutions of the approximate equation
\begin{equation}\label{appeq}
r_n = -A_n-B_n/r_n
\end{equation}

The proof of the estimate of $|r_n/\tilde{r}_n-1|$ is based on calculations with polynomials with rational coefficients and is rigorous.

The boundary of $S_4$ (see \eqref{S4}) is partitioned into 6 line segments,
\begin{align}\label{l16}
\begin{array}{cccc}
& l_1=[10\,i,70\,i] & l_3=[10\,i+1/2,70\,i+1/2] & l_5=[10\,i,10\,i+1/2] \\
& l_2=(70\,i,380\,i] & l_4=(70\,i+1/2,380\,i+1/2] & l_6=[380\,i,380\,i+1/2] \\
\end{array}
\end{align}

Let $\lambda = t\,i+s$. We choose the following quasi-solution 
\begin{equation}\label{r1def}
r^{(1)}_n = e^{W_n(\lambda)}\quad\quad \left( 1\leq n\leq 50\right)
\end{equation}
to $r_n$, where
\begin{align}\label{Wndef}
W_n(\lambda) =\left\{
\begin{array}{lr}
\sum_{i=0}^{5} {\left((1-2s)a^{(1)}_{n,i}+2s\,a^{(2)}_{n,i}\right)} T_i\left(\frac{t-40}{30}\right):&10 \leq t\leq 70\\
\sum_{i=0}^{5} {\left((1-2s)b^{(1)}_{n,i}+2s\,b^{(2)}_{n,i}\right)} T_i\left(\frac{t-225}{155}\right):&70 < t\leq 380
\end{array}
\right.
\end{align}
where $a_{n,i}^{(j)}$ and $b_{n,i}^{(j)}$ $(1\leq j\leq 2,\, 1\leq i\leq 5,\, 1\leq n\leq 50)$ are given in Table \ref{table} in the Appendix. Define also 
\begin{equation}\label{r2def}
r^{(2)}_n = -\frac{A_n}{2}\left(1+\sqrt{1-\frac{4\,B_n}{A^2_n}}\right)
\end{equation}
where the branch of square root is the one which is positive on the positive real axis. We note that $r^{(2)}_n$ is a solution to the quadratic equation \eqref{appeq}. Also write
\begin{align}\label{delta12}
\begin{aligned}
r_n &= r^{(1)}_n\left(1+\delta^{(1)}_n\right) \\
r_n &= r^{(2)}_n\left(1+\delta^{(2)}_n\right) \\
\end{aligned}
\end{align}

\subsection{Proof of the Theorem}
We prove Theorem \ref{MainTh2} in four steps.

\begin{Lemma}\label{L14}
	
	(i)\begin{equation}\label{L1}
	\left\vert \delta^{(1)}_1 (\lambda) \right\vert \leq D_1 = \frac{19}{1000}\quad\quad \forall \lambda \in \partial S_4
	\end{equation}
	(ii)\begin{equation}\label{L2}
	\left\vert \delta^{(1)}_n (\lambda) \right\vert \leq D_2 = \frac{9}{200}\quad\quad \forall \lambda \in \partial S_4,\, 1\leq n\leq 50
	\end{equation}
	(iii)\begin{equation}\label{L3}
	\left\vert \delta^{(2)}_{50} (\lambda) \right\vert \leq D_3 = \frac{3}{20}\quad\quad \forall \lambda \in \partial S_4
	\end{equation}
	(iv)\begin{equation}\label{L4}
	\left\vert \delta^{(2)}_{n} (\lambda) \right\vert \leq D_3 =\frac{3}{20} \quad\quad \forall \lambda \in \partial S_4,\, n\geq 50
	\end{equation}
\end{Lemma}

\begin{Corollary}\label{C5}
	\begin{equation}
	\lim_{n\to\infty} r_{n} = 1\quad\quad \forall \lambda \in S_4
	\end{equation}
\end{Corollary}

\subsubsection{Proof of Lemma \ref{L14} (i)} Direct application of an elementary method described in Section \ref{esmet} in the Appendix.

\subsubsection{Proof of Lemma \ref{L14} (ii)} By the definitions \eqref{r1def} and \eqref{delta12} we have for $1\leq n\leq 50$.
\begin{equation}\label{r1def2}
r_n^{(1)}(\lambda) = e^{W_n (\lambda)}\left(1+\delta^{(1)}_n(\lambda)\right)
\end{equation}
Substituting the equation \eqref{r1def2} into \eqref{rec101} we obtain
\begin{equation}\label{rec103}
\delta^{(1)}_{n+1}(\lambda) = \epsilon^{(1)}_{n+1}+C^{(1)}_n\frac{\delta^{(1)}_{n}(\lambda)}{1+\delta^{(1)}_{n}(\lambda)}
\end{equation}
where
\begin{align}
\epsilon^{(1)}_{n+1}&=-1-A_n e^{-W_{n+1}}-B_n e^{-W_{n+1}-W_n}\\
C^{(1)}_n&=B_n e^{-W_{n+1}-W_n}
\end{align}
We use the method in Section \ref{esmet} to obtain the following upper bounds
\begin{align}
&\left\vert \epsilon^{(1)}_{n+1}\right\vert \le 21/1000	\quad\left(1\le n\le49,\,\lambda\in l_1\cup l_3 \cup l_5\right)\label{es1}\\ 
&\left\vert \epsilon^{(1)}_{n+1}\right\vert \le {27/1000}	\quad\left(1\le n\le49,\,\lambda\in l_2\cup l_4 \cup l_6\right) \label{es2}\\
&\left\vert C^{(1)}_{n}\right\vert \le 1/2	\quad\left(1\le n\le49,\,\lambda\in l_1\cup l_3 \cup l_5\right)\label{es3}\\
&\left\vert C^{(1)}_{n}\right\vert \le 1/3	\quad\left(1\le n\le49,\,\lambda\in l_2\cup l_4 \cup l_6\right)\label{es4}
\end{align}
Hence if $\delta^{(1)}_{n}\le D_2 = 21/500$, it is straightforward to check that $\delta^{(1)}_{n+1}\le D_2$ by \eqref{rec103} and the estimates \eqref{es1} -- \eqref{es4} hold.

\subsubsection{Proof of Lemma \ref{L14} (iii)}
Define $\delta^{(3)}$ by
\begin{equation}\label{d3def}
r^{(2)}_{50}=r^{(1)}_{50}\left(1+\delta^{(3)}\right)
\end{equation} 

By the definition \eqref{delta12} we have
\begin{equation} \label{del1del2}
\left(1+\delta^{(1)}_{50}\right)r^{(1)}_{50} = \left(1+\delta^{(2)}_{50}\right)r^{(2)}_{50}
\end{equation}
Substituting \eqref{d3def} into \eqref{del1del2} we have
\begin{equation}
\delta^{(2)}_{50} = \frac{\delta^{(1)}_{50}-\delta^{(3)}}{1+\delta^{(3)}} 
\end{equation}
Thus
\begin{equation}\label{ubdel2}
\left\vert \delta^{(2)}_{50} \right\vert \leq \frac{\left\vert \delta^{(1)}_{50}\right\vert +\left\vert \delta^{(3)}\right\vert }{1-\left\vert \delta^{(3)}\right\vert} 
\end{equation}

By definitions \eqref{r1def} and \eqref{d3def} we have
\begin{equation}\label{d3def2}
r_{50}^{(2)}(\lambda) = e^{W_{50} (\lambda)}\left(1+\delta^{(3)}(\lambda)\right)
\end{equation}
Substituting \eqref{d3def2} into \eqref{appeq} we get
\begin{equation}\label{rec102}
\delta^{(3)} = \epsilon^{(3)}+C^{(3)}\frac{\delta^{(3)}}{1+\delta^{(3)}}
\end{equation}
where
\begin{align}
\begin{aligned}
\epsilon^{(3)}&=-1-A_{50}\,e^{-W_{50}}-B_{50}\,e^{-2\,W_{50}}\\
C^{(3)}&=B_{50}\, e^{-2\,W_{50}}
\end{aligned}
\end{align}
Using the method in Section \ref{esmet} we obtain
\begin{align}\label{es5}
\begin{aligned}
\left\vert \epsilon^{(3)}\right\vert &<\frac{1}{50}\\
\left\vert C^{(3)}\right\vert &\leq 1/2
\end{aligned}\quad\quad \left(\lambda \in l_1\cup l_3 \cup l_5 \right)
\end{align}
and
\begin{align}\label{es6}
\begin{aligned}
\left\vert \epsilon^{(3)}\right\vert &<\frac{19}{500}\\
\left\vert C^{(3)}\right\vert &\leq 1/3
\end{aligned}\quad\quad \left(\lambda \in l_2\cup l_4 \cup l_6 \right)
\end{align}
Rewriting \eqref{rec102}, we obtain a quadratic equation for $\delta^{(3)}$
\begin{equation}\label{rewrec2}
\left(\delta^{(3)}\right)^2+\left(1-C^{(3)}-\epsilon^{(3)}\right)\delta^{(3)}-\epsilon^{(3)} = 0
\end{equation}
Denoting the two roots of the equation above by $s_1$ and $s_2$ and using the estimates \eqref{es5} and \eqref{es6} we have
\begin{align}\label{ds12}
\left\vert|s_1|-|s_2|\right\vert >19/50\quad\quad \left(\lambda \in l_1\cup l_3 \cup l_5 \right)\\
\left\vert|s_1|-|s_2|\right\vert > {49}/{100}\quad\quad \left(\lambda \in l_2\cup l_4 \cup l_6 \right)
\end{align}
Assume without loss of generality that $s_1(\lambda)$ is the root with the larger modulus. Then \eqref{ds12} implies
\begin{align}
\begin{aligned}\label{ubdel31}
|s_1|&>19/50\\
|s_2|&=|\epsilon^{(3)}/s_1|<\frac{{1}/{50}}{{19}/{50}}=\frac{1}{19}
\end{aligned}\quad\quad \left(\lambda \in l_1\cup l_3 \cup l_5 \right)\\
\begin{aligned}\label{ubdel32}
|s_1|&>49/100\\
|s_2|&=|\epsilon^{(3)}/s_1|<\frac{{19}/{500}}{{49}/{100}}=\frac{19}{245}
\end{aligned}\quad\quad \left(\lambda \in l_2 \cup l_4 \cup l_6 \right)
\end{align}
By definition, $\delta^{(3)}$ is continuous on $l_1\cup l_3\cup l_5$ and on $l_2\cup l_4\cup l_6$. In addition, using the definition of $r_n^{(1)}$ in \eqref{r1def} and \eqref{Wndef}, the definition of $r_n^{(2)}$ in \eqref{r2def} and the definition of $\delta^{(3)}$ in \eqref{d3def} we can check that
\begin{align}
\left\vert\delta^{(3)}{\left(10 \,i\right)}\right\vert<1/100\,\,, \quad\quad
\left\vert\delta^{(3)}{\left(380 \,i\right)}\right\vert<1/25 
\end{align}
This implies that for all $\lambda\in\partial S_4$, $\delta^{(3)}(\lambda)$ is the root $s_2$ of \eqref{rewrec2}. Lemma \ref{L14} (iii) follows from the estimates \eqref{ubdel31}, \eqref{ubdel32}, Lemma \ref{L14} (ii) and \eqref{ubdel2}.

\subsubsection{Proof of Lemma \ref{L14} (iv)}
Substituting the second equation in \eqref{delta12} into \eqref{rec101} and using the fact that $r^{(2)}_n$ is a solution of the equation \eqref{appeq} we obtain
\begin{equation}\label{rec104}
\delta^{(2)}_{n+1} = \epsilon^{(2)}_n+C^{(2)}_n \frac{\delta^{(2)}_n}{1+\delta^{(2)}_n}
\end{equation}
where
\begin{align}
\epsilon^{(2)}_n&=\frac{r^{(2)}_{n}}{r^{(2)}_{n+1}}-1\\
C^{(2)}_n&=\frac{B_n}{r^{(2)}_{n}\, r^{(2)}_{n+1}}
\end{align}
Denoting
\begin{equation}
F_n = 1-\frac{4\,B_n}{A^2_n}
\end{equation}
we get
\begin{align} \label{r2def3}
r^{(2)}_n = -\frac{A_n}{2}\left(1+\sqrt{F_{n}}\right)
\end{align}
and 
\begin{align}\label{es7}
\begin{aligned}
\epsilon^{(2)}_n &=\frac{A_n\left(1+\sqrt{F_n}\right)}{A_{n+1}\left(1+\sqrt{F_{n+1}}\right)}-1
=\frac{A_n}{A_{n+1}} \left(\frac{1+\sqrt{F_n}}{1+\sqrt{F_{n+1}}}-1\right) +\left(\frac{A_n}{A_{n+1}}-1\right)\\
&=\frac{A_n}{A_{n+1}} \frac{{F_n}-{F_{n+1}}}{\left(1+\sqrt{F_{n+1}}\right)\left(\sqrt{F_{n}}+\sqrt{F_{n+1}}\right)} +\left(\frac{A_n}{A_{n+1}}-1\right)\\
\end{aligned}
\end{align}
\begin{align}\label{es8}
\begin{aligned}
C^{(2)}_n&=\frac{4\,B_n}{A_n A_{n+1}\left(1+\sqrt{F_{n}}\right)\left(1+\sqrt{F_{n+1}}\right)}
\end{aligned}
\end{align}

We estimate $\epsilon^{(2)}_n$ and $C^{(2)}_n$ by first showing the following inequalities valid for $\lambda$ on all sides of $\partial S_4$ and $n \ge 50$, proved in Section \ref{L4App}

\begin{equation}\label{U1}
\left\vert \frac{A_n}{A_{n+1}}-1\right\vert \le U_1 := \frac{1}{30}
\end{equation}

\begin{equation}\label{U2}
\left\vert F_n-F_{n+1}\right\vert \le U_2 := \frac{7\sqrt{2}}{500}
\end{equation}

\begin{equation}\label{numberL1}
 \Re{\sqrt{F_{n}}}  \ge L_1 := \frac{8}{25}
\end{equation}

\begin{equation}\label{U3}
\left\vert \frac{B_n}{A_n A_{n+1}}\right\vert \le U_3 := \frac{9}{40}
\end{equation}

\indent Using \eqref{U1} -- \eqref{U3} in \eqref{es7} and \eqref{es8} we have
\begin{align}
\left\vert \epsilon_n^{(2)} \right\vert& \leq \left(1+U_1\right)\frac{U_2}{\left(1+L_1\right)\left(2\,L_1\right)}+U_1 < 29/500\label{es9}\\
\left\vert C_n^{(2)} \right\vert&\leq \frac{4 U_3}{\left(1+L_1\right)^2} < 517/1000\label{es10}
\end{align}
\indent Hence $\left\vert \delta^{(2)}_{n} \right\vert \le D_3 = 3/20$ implies $\left\vert \delta^{(2)}_{n+1} \right\vert \leq D_3$ by \eqref{rec104}, the estimates \eqref{es9} and \eqref{es10}. Also from Lemma \ref{L14} (iii) we know that $\left\vert \delta^{(2)}_{50} \right\vert \leq D_3$. Lemma \ref{L14} (iv) is proved by straightforward induction.\\

\subsubsection{Proof of Corollary \ref{C5}}
First we show that for each $n>0$, $r_n$ is analytic in $S_4$ since we later make use of the maximum principle. We provide below all the details, but conceptually the proof is simple: from the recurrence \eqref{rec101} one can easily check that if $r_n$ has a pole, then $r_{n-1}$ has a zero; we use the argument principle to show that there are no zeros in the strip of interest.\\
\indent Each $r_n$ is meromorphic since each $c_n$ is a polynomial by recurrence \eqref{rec5} and \eqref{coeffps}. On the boundary $\partial S_4$, the boundedness of the quasi-solutions $r^{(1)}_n$ and $r^{(2)}_n$ together with the estimates in Lemma \ref{L14} show that $r_n$ is analytic for each $n$ and therefore also nonzero. Assume that $n=n_0$ is the smallest index such that $r_n$ has a zero in $S_4$. $n_0>1$ obviously.\\
 \indent If $n_0 \geq 50$, $r^{(2)}_{n_0}$ is manifestly analytic. Definition of $\delta^{(2)}$ in \eqref{delta12} and Lemma \ref{L14} (iv) imply that
\begin{equation}
\left\vert r_{n_0}- r^{(2)}_{n_0} \right\vert < \frac{3}{20} \left\vert r^{(2)}_{n_0} \right\vert < \left\vert r^{(2)}_{n_0} \right\vert \quad\quad  (\lambda \in \partial S_4)
\end{equation}
By Rouch\'e's Theorem $r_{n_0}$ has the same number of roots as $r^{(2)}_{n_0}$ in $S_4$, which is zero. \\
\indent If $n_0 < 50$, let $\phi:[0,741] \to \partial S_4$ denote a parameterization of $\partial S_4$ by arc length (the perimeter of $S_4$ is $741$), satisfying $\phi(0) = \phi(741) = 10 i$ and that as t increases $\partial S_4$ is traversed counterclockwise. We observe that the $W_n$ given in \eqref{Wndef} is continuous in $S_4 \backslash \{\lambda: \Im \lambda = 70\}$. It is straightforward to check that for each $n<50$ there exists a suitable integer $k_n$ such that if we redefine $W_n$ as $W_n+2 \pi k_n i$ in $\{\lambda: 70<\Im \lambda \leq 380\}$ we obtain the following estimates at the only two discontinuities $70 i$ and $70 i +1/2$ of $W_n$ on $\partial S_4$
\begin{align}
\begin{aligned}
& \lim_{t\to 0^+}\left\vert W_n\left(70i+t i\right)-W_n\left(70i-t i\right)\right\vert < 2/250\\
& \lim_{t\to 0^+}\left\vert W_n\left(70i+1/2+t i\right)-W_n\left(70i+1/2-t i\right)\right\vert < 7/1000
\end{aligned}
\end{align}
 Thus there exists a function $\widetilde{W}_n$ which is continuous on $\partial S_4 \backslash \{10 i\}$ and such that both $\left(\widetilde{W} \circ \phi\right)(0+)$ and $\left(\widetilde{W}\circ \phi\right)(741-)$ exist, $\widetilde{W}(10i) = \left(\widetilde{W} \circ \phi\right)(0+)$ and the following holds:
 \begin{equation}\label{logdif2}
 \left\vert W_n - \widetilde{W}_n \right\vert <2/250			\quad\quad(\lambda \in \partial S_4)
 \end{equation}
 Consequently, if we let $\tilde{\delta}_n := e^{W_n}/e^{\widetilde{W}_n}-1$ it is straightforward from \eqref{logdif2} that
 \begin{equation}\label{tildel}
 \left\vert \tilde{\delta}_n \right\vert < 13/1000
 \end{equation}
Next we choose a branch $\overline{W}_n$ of $\log r_n$ on $\partial S_4$, which is continuous on $\partial S_4 \backslash \{10 i\}$ with $\overline{W}_n(10i):=\overline{W}_n\left(\phi(0+)\right) $. Denoting $\tilde{\tilde{\delta}}_n = r_n / e^{\widetilde{W}_n}-1$, we have from \eqref{tildel} and Lemma \ref{L14} (ii) that
\begin{equation}
 \left\vert \tilde{\tilde{\delta}}_n \right\vert= \left\vert e^{\overline{W}_n-\widetilde{W}_n}-1\right\vert = \left\vert \left(1+\delta^{(1)}_n\right)\left(1+\tilde{\delta}_n\right)-1 \right\vert < 3/50
\end{equation}
Therefore $\left\vert \mathrm{ln}(1+\tilde{\tilde{\delta}}_n) \right\vert < 1/10$ where by definition $\Im (\ln z) \in (-\pi,\pi]$. For each $\lambda \in \partial S_4$ we have
\begin{equation}
 \overline{W}_n-\widetilde{W}_n  = \ln(1+\tilde{\tilde{\delta}}_n)+2 \pi k i  \quad\quad (k \in \ZZ)
\end{equation}
In other words, for each $t \in [0,741)$ we have
\begin{equation}\label{disjointunion}
\left(\overline{W}_n \circ \phi \right)(t) - \left(\widetilde{W}_n \circ \phi \right)(t) \subset \bigcup_{k\in \ZZ} B(0,1/10)+2 k \pi i
\end{equation} 
where $B(0,1/10) = \{z: |z|<1/10\}$. Since $\overline{W}_n\circ \phi - \widetilde{W}_n \circ \phi$ is continuous on $[0,741)$, its image must be connected, and the right hand side of \eqref{disjointunion} is a disjoint union of open balls, so the image must be contained in only one of them, say $B(0,1/10)+2 k_0 \pi i$. Then we have
\begin{align}
\begin{aligned}
&\left\vert \overline{W}_n \left( \phi (0+)\right) - \overline{W}_n \left( \phi (741-)\right) \right\vert \\
& < \left\vert \left(\overline{W}_n -\widetilde{W}_n\right)\left( \phi (0+)\right)   -  \left(\overline{W}_n -\widetilde{W}_n\right)\left( \phi (741-)\right)  \right\vert + \left\vert \widetilde{W}_n \left( \phi (0+)\right) - \widetilde{W}_n \left( \phi (741-)\right) \right\vert\\
& < 1/10 \cdot 2+2/250 \cdot 2 <3/10
\end{aligned}
\end{align}
By argument principle $\overline{W}_n \left( \phi (0+)\right) = \overline{W}_n \left( \phi (741-)\right)$ and thus $r_{n_0}$ has no zeros in $S_4$. Therefore each $r_n$ is analytic in $S_4$.\\
\indent We know that for all $\lambda \in \CC$ we have 
\begin{align}
& \lim_{n\to\infty}A_n = -\frac{3}{2}\\
& \lim_{n\to\infty}B_n = \frac{1}{2}
\end{align}
so by the definition \eqref{r2def},
\begin{equation}
\lim_{n\to\infty}r^{(2)}_n = 1
\end{equation}
Also, for each $n\geq 50 $, $r_n^{(2)}$ is analytic in $S_4$, thus $\delta_n^{(2)} = r_n/r^{(2)}_n-1$ is also analytic in $S_4$. Therefore the estimate \eqref{L4} in Lemma \ref{L14} (iv) and the maximum modulus principle imply that
\begin{equation}\label{d2bd}
\left\vert \delta_n^{(2)}\right\vert \le 3/20 \quad\quad (\lambda \in S_4)
\end{equation}
In addition \eqref{2lim} implies that either $\lim_{n \to \infty} \delta_n^{(2)} = 0$, which means $r_n \to 1$ or $\lim_{n \to \infty} \delta_n^{(2)} = -\frac{1}{2}$, which means $r_n \to 1/2$. Since \eqref{d2bd} implies $\lim_{n \to \infty} \delta_n^{(2)} \neq -\frac{1}{2}$, Corollary \ref{C5} is proved.

\section{Appendix}

\subsection{Bounds of polynomials and Proof of Lemma \ref{estimep}}
\label{subsec:estipoly}
First we describe the method used  to estimate upper bounds of absolute values of polynomials of the form
\begin{equation}
P(x,\lambda)=\sum_{i,j} a_{i,j} x^i \lambda^j
\end{equation}
where $x \in [-1,1]$, $\lambda \in S_1 \subseteq B(\frac{1+i}{4},\frac{\sqrt{2}}{4})$, where $B(\frac{1+i}{4},\frac{\sqrt{2}}{4})$ denotes the closed disk centered at $\frac{1+i}{4}$ with radius $\frac{\sqrt{2}}{4}$.\\
\indent Reexpand the polynomial $P(x,\lambda)$ at $\lambda_0 = \frac{1+i}{4}$.
\begin{equation}
P(x,\lambda)=\sum_{i,j} \tilde{a}_{i,j} x^i (\lambda-\lambda_0)^j
\end{equation}
Then we express each $x^i$ as a linear combination of Chebyshev polynomials of the first kind.
\begin{equation}
x^i = \sum_{l=0}^i c_l T_l (x)
\end{equation}
Rewrite the polynomial again
\begin{align}
P(x,\lambda)	&=\sum_{i,j} \tilde{a}_{i,j} \left(\sum_{l=0}^i c_l T_l (x)\right) (\lambda-\lambda_0)^j \\
					&=\sum_{i,j}\sum_{l=0}^i \tilde{a}_{i,j}c_l T_l (x) (\lambda-\lambda_0)^j
\end{align}
and use the facts that for each $n$, $|T_n(x)| \leq 1$ on $[0,1]$ and $|\lambda-\lambda_0| \leq \frac{\sqrt{2}}{4}$ to obtain the following estimate:
\begin{equation}
| P(x,\lambda)|	
					\leq \sum_{i,j}\sum_{l=0}^i |\tilde{a}_{i,j}| |c_l |\left(\sqrt{2}/4\right)^j
\end{equation}
\\ \textbf{Proof of Lemma \ref{estimep}} 
As explained before Lemma \ref{estimep}, 

$$\epsilon_1(t,\lambda) = \mathcal{L} \delta_1=-\mathcal{L} G_a = t \, \tilde{\epsilon}_1(t,\lambda)$$

For $t\in (0,1/4)$ and $\lambda \in S_1$ we are interested in an upperbound of $|\tilde{\epsilon}_1(t,\lambda)|$. For $t\in (1/4,1/2)$ and $\lambda \in S_1$ we are interested in an upperbound of $|\epsilon_1(t,\lambda) |=|\mathcal{L} G_a |$.

Let 
\begin{align}
R_1(t,\lambda) 	&:=(2-t)^2 \tilde{\epsilon}_1(t,\lambda) \quad\quad t\in(0,1/4)\\
R_2(t,\lambda) 	&:=(2-t)^2 \epsilon_1(t,\lambda)\quad\quad t\in(1/4,1/2)
\end{align}

$R_1$ and $R_2$ are both polynomials in $t$ and $\lambda$. Make a change of variable $t=(x+1)/8$, then an upper bound of  $|R_1((x+1)/8,\lambda)|$ with $x \in [-1,1]$, $\lambda\in S_1$ is an upper bound of $|R_1(t,\lambda)|$ with $t \in [0,1/4]$,$\lambda\in S_1$, which is what we seek. Using the method mentioned in the beginning of the section, we obtain the following:
\begin{align}
|R_1((x+1)/8,\lambda)| &\leq 0.00441 \quad\quad x\in[-1,1] \quad \lambda\in S_1 \\
i.e.\quad\quad|R_1(t,\lambda)| &\leq 0.00441 \quad\quad t\in[0,1/4]\quad \lambda\in S_1 
\end{align}
Hence
\begin{equation}
|\tilde{\epsilon}_1 (t,\lambda)| = \left\vert \frac{R_1(t,\lambda)}{(2-t)^2}\right\vert \leq \frac{0.00441}{(2-1/4)^2} = 0.00144
\end{equation}

For $R_2(t,\lambda)$, let $t=(x+3)/8$, then we an upper bound of $|R_2((x+3)/8,\lambda)|$ with $x \in [-1,1]$, $\lambda\in S_1$ is what we seek. Apply the method to $|R_2((x+3)/8,\lambda)|$ we get
\begin{align}
|R_2((x+3)/8,\lambda)| &\leq 0.00595 \quad\quad x\in[-1,1] \quad \lambda\in S_1 \\
i.e.\quad\quad|R_2(t,\lambda)| &\leq 0.00595 \quad\quad t\in[1/4,1/2]\quad \lambda\in S_1 
\end{align}
Hence
\begin{equation}
|\epsilon_1 (t,\lambda)| = \left\vert \frac{R_2(t,\lambda)}{(2-t)^2}\right\vert \leq \frac{0.00595}{(2-1/2)^2} < 0.00265
\end{equation}
The other estimates in Lemma \ref{estimep} are trivial.

\subsection{{Proof} of the estimate \ref{loubou}}\label{Ploubou}

Let $W(\lambda)= W(G_a,H_a) (1/2,\lambda)$. $W(\lambda)$ is a polynomial in $\lambda$, the region of interest is $\lambda\in S_1$.\\
\indent First take $W_0 = W(1/2,\frac{1+i}{4})$, which is just a constant with $|W_0| > 0.79$. Then use the method described in the previous section to acquire an upper bound of $|W(\lambda) -W_0 |$. Notice that $W(\lambda) - W_0$ is just a special case of the class of polynomial discussed, so we obtain:
\begin{equation}
|W(\lambda) -W_0 | \leq 0.33 		\quad \lambda\in S_1
\end{equation}
Thus
\begin{equation}
|W(\lambda)| \geq |W_0| - |W(\lambda) -W_0 | > 0.79-0.33>0.46
\end{equation}

Sharper estimates can be achieved at the cost of simplicity of the proof.

\subsection{{Proof} of Lemma \ref{JnMnLemma}}\label{JnMnProof}

For $n \leq 10$ we  obtain $J_n$ and $M_n$ using the method in \S\ref{rationalS2} and \S\ref{polynS2}. Using (\ref{rec1}) and (\ref{rec2}), we obtain:
 \begin{equation}
\label{rec3}
|a_{n+1}| \leq |p_n(\lambda)| |a_n| + \frac{1}{(n+1)|n+\lambda|}\sum_{k=0}^{n-1}\frac{n-k+1}{2^{n-k}}|a_k|
\end{equation}
\begin{equation}
\label{rec4}
|b_{n+1}| \leq |q_n(\lambda)| |b_n| + \frac{1}{(n+1)(n+\frac{5}{2})}\sum_{k=0}^{n-1}4(n-k+1)|b_k|
\end{equation}
which give upper bounds $J_n$ of $|a_n|$ and $M_n$ of $|b_n|$ for $11\leq n \leq 30$, $\lambda \in S_2$. The details are as follows.

We first estimate the difference between $p_n(\lambda)$ and its leading behavior for $n\geq 30$:

\begin{equation}
\label{en}
e_n(\lambda) := p_n(\lambda)-\left(1+\frac{1}{2n}\right)=\frac{1}{4}\,{\frac {{\lambda}^{2}n-3\,\lambda\,n-2\,\lambda-4\,n}{ \left( n+1
 \right)  \left( n+\lambda \right) n}}
\end{equation}
With the notation 
\begin{align}
\lambda&= r+i\, s \quad\quad (0\leq r \leq 1/2, 1/2 \leq s \leq 4)\\
z&=1/n \quad\quad\,\,\,\,\,\,(0 < z \leq 1/30)
\end{align}
we have
\begin{align}
\label{enr}
\Re(e_n) &= {\frac {{z}^{2} \left( {r}^{3}z-2\,{r}^{2}{z}^{2}+r{s}^{2}z-2\,{s
}^{2}{z}^{2}-3\,{r}^{2}z-3\,{s}^{2}z+{r}^{2}-6\,rz-{s}^{2}-3\,r-4
 \right) }{4\, \left( 1+z \right)  \left( {r}^{2}{z}^{2}+{s}^{2}{z}^{2}+2
\,rz+1 \right) }}\nonumber\\
              &\leq {\frac {{z}^{2} \left( {r}^{3}z+r{s}^{2}z+{r}^{2}-4
 \right) }{4\, \left( 1+z \right)  \left( {r}^{2}{z}^{2}+{s}^{2}{z}^{2}+2
\,rz+1 \right) }}\nonumber\\
             &\leq {\frac {{z}^{2} \left( {0.5}^{3}/30+(0.5)\,{4}^{2}/30+{0.5}^{2}-4
 \right) }{4\, \left( 1+z \right)  \left( {r}^{2}{z}^{2}+{s}^{2}{z}^{2}+2
\,rz+1 \right) }} <0
\end{align}
and
\begin{align}
\Re(e_n) &= {\frac {{z}^{2} \left( {r}^{3}z-2\,{r}^{2}{z}^{2}+r{s}^{2}z-2\,{s
}^{2}{z}^{2}-3\,{r}^{2}z-3\,{s}^{2}z+{r}^{2}-6\,rz-{s}^{2}-3\,r-4
 \right) }{4\, \left( 1+z \right)  \left( {r}^{2}{z}^{2}+{s}^{2}{z}^{2}+2
\,rz+1 \right) }}\nonumber\\
              &\geq {\frac {{z}^{2} \left( -2\,{r}^{2}{z}^{2}-2\,{s
}^{2}{z}^{2}-3\,{r}^{2}z-3\,{s}^{2}z-6\,rz-{s}^{2}-3\,r-4
 \right) }{4\, \left( 1+z \right)  \left( {r}^{2}{z}^{2}+{s}^{2}{z}^{2}+2
\,rz+1 \right) }} \geq -\frac{4187}{648000}
\end{align}

\begin{equation}
\label{eni}
\Im(e_n) = \frac{1}{4}\,{\frac {{z}^{2}s \left( {r}^{2}z+{s}^{2}z+2\,r+2\,z-3 \right) }{
 \left( 1+z \right)  \left( {r}^{2}{z}^{2}+{s}^{2}{z}^{2}+2\,rz+1
 \right) }}
\end{equation}
Hence
\begin{equation}
1/100 < \Re(e_n) < 0 
\end{equation}
\begin{equation}
0> \Im(e_n)> -\frac{3 z^2}{z^2(r^2+s^2)+1} \geq -\frac{3}{r^2+s^2+1/z^2} \geq -\frac{3}{1/4+n^2}
\end{equation}
for $\lambda \in S_2$ and $n \geq 30$. So
\begin{equation}
|p_n|^2 = \left(1+\frac{1}{2n} + \Re(e_n)\right)^2+\Big(\Im(e_n)\Big)^2 \leq \left(1+\frac{1}{2n}\right)^2+\left(\frac{3}{1/4+n^2}\right)^2
\end{equation}
Let 
\begin{equation}
P_{n} =\left( \left(1+\frac{1}{2n}\right)^2+\left(\frac{3}{1/4+n^2}\right)^2\right)^{1/2}
\end{equation} 
then $|p_n| \leq P_{30}$ for $n \geq 30$. For $10\leq n \leq 29$, we obtain $P_n$ using the method in \S\ref{rationalS2}, see \S\ref{valJnMn} for these values of $P_n$. From \eqref{rec3} we have:
\begin{equation}
|a_{n+1}| \leq P_n  |a_n| + \frac{1}{(n+1)|n+i/2|}\sum_{k=0}^{n-1}\frac{n-k+1}{2^{n-k}}|a_k| \quad (n\geq 10)
\end{equation}
For $10 \leq n\leq 29$ if we choose $J_{n+1}$ such that
\begin{equation}
J_{n+1} > P_n \,J_n + \frac{1}{(n+1)|n+i/2|}\sum_{k=0}^{n-1}\frac{n-k+1}{2^{n-k}}J_k
\end{equation}
then
\begin{equation}
|a_{n+1}|  < J_{n+1} 	\quad\quad (10 \leq n\leq 29)
\end{equation}

We use induction to take care of larger $n$. Explicit calculation shows that for our choice of $k_1$ and $J_n$, $J_{30}<k_1^{30}$ (see \S\ref{valJnMn}). For $n>30$, assuming for $30 \leq m \leq n$, $|a_m| \leq k_1^m$,  we have:
\begin{align}
|a_{n+1}| & < P_n |a_n| + \frac{1}{(n+1)n} \sum_{i=0}^{n-1}\frac{n-i+1}{2^{n-i}} |a_i|	\nonumber\\
			& < P_n k_1^n + \frac{1}{(n+1)n}\sum_{i=0}^{n-1}\frac{n-i+1}{2^{n-i}} k_1^i +\frac{1}{(n+1)n}\sum_{i=0}^{n-1}\frac{n-i+1}{2^{n-i}} (|a_i|-k_1^i)\nonumber\\
               & < P_{30} k_1^n+{\frac {4\, \left( 2\,k_{{1}} \right) ^{n}k_{{1}}-2\,nk_{{1}}- \left( 
2\,k_{{1}} \right) ^{n}-4\,k_{{1}}+n+1}{{2}^{n} \left(2\,k_{{1}}-1\right) ^{2} \left( n+1 \right) n}}	\nonumber\\
			&\qquad +\frac{1}{(n+1)n}\sum_{i=0}^{30}\frac{n-i+1}{2^{n-i}} (J_i-k_1^i)		\nonumber\\
			& < P_{30} k_1^n+k_1^n\,{\frac {4k_1-1}{{(2\,k_1-1)^2 \left( n+1 \right) n}}} +\frac{(1-2k_1)n+1-4k_1}{2^n(n+1)n(2k_1-1)^2}	\nonumber\\
			&\qquad +\frac{1}{(n+1)n}\,(\tilde{V}_1\,n +\tilde{V}_0)		\nonumber\\
			& <\left(P_{30}+\frac {4k_1-1}{{(2\,k_1-1)^2 \left( n+1 \right) n}}\right) k_1^n+ \frac{\left(\frac{-1}{2k_1-1}+\tilde{V}_1\right)n+\frac{1-4k_1}{(2k-1)^2}+\tilde{V}_0}{2^n n(n+1)}		\nonumber\\
						& <(1.01668+0.00301)\, k_1^n		\nonumber\\
			&<1.0197\, k_1^n <k_1^{n+1}
\end{align}

where $\tilde{V}_1=-196921202$ and $\tilde{V}_0=5563721416$ are chosen to satisfy 
\begin{equation}
\sum_{i=0}^{30}\frac{n-i+1}{2^{-i}} (J_i-k_1^i)\,<\,\tilde{V}_1\,n +\tilde{V}_0
\end{equation}
and explicit calculation shows that
\begin{equation}
\left(\frac{-1}{2k_1-1}+\tilde{V}_1\right)n+\frac{1-4k_1}{(2k-1)^2}+\tilde{V}_0 <0  \quad (n \geq 30)
\end{equation}

So Lemma \ref{JnMnLemma} (i) is proved.\\

Part (ii) is similar. It is straightforward to check that for each $n\geq 2$, $|q_n|$ is monotonic on each line segment of $\partial S_2$ and the maximum is attained at $\lambda= 4\,i+1/2$. Let
\begin{equation}
Q_n := |q_n(4\,i+1/2)|	\quad\quad (n\geq 10)
\end{equation}
Then $|q_n(\lambda)| \leq Q_n$ for all $\lambda \in S_2$ and $n \geq 10$. It is also straightforward to check that $|Q_n|<1$ for all $n\geq 2$.\\

From \eqref{rec4} and $M_n$ ($n\leq 10$) we choose $M_n$ recursively for $11 \leq n \leq 30$ such that
\begin{equation}
M_{n+1} >  Q_n M_n + \frac{1}{(n+1)(n+\frac{5}{2})}\sum_{k=0}^{n-1}4(n-k+1)M_k
\end{equation}
so
\begin{equation}
|b_{n+1}| \leq |q_n(\lambda)| |b_n| + \frac{1}{(n+1)(n+\frac{5}{2})}\sum_{k=0}^{n-1}4(n-k+1)|b_k| \, < M_{n+1}
\end{equation}

It is easy to check (see \S\ref{valJnMn}) that $M_{30}<k_2^{30}$. For $n>30$, assuming for all $30 \leq m \leq n$, $|b_m| \leq k_2^m$,  we have:
\begin{align}
|b_{n+1}| & < |b_n| + \frac{4}{(n+1)(n+5/2)} \sum_{i=0}^{n-1}(n-i+1) |b_i|	\nonumber\\
			& < k_2^n +  \frac{4}{(n+1)(n+5/2)} \sum_{i=0}^{n-1}(n-i+1) k_2^i \nonumber\\
			&\qquad+\frac{4}{(n+1)(n+5/2)} \sum_{i=0}^{30}(n-i+1)\,(M_i-k_2^i)\nonumber\\
               & < k_2^n+\frac{4}{(n+1)(n+5/2)}\left[k_2^n\frac{2k_2-1}{(k_2-1)^2}+\frac{n(1-k_2)+1-2k_2}{(k_2-1)^2}+\tilde{W}_1n +\tilde{W}_0\right]	\nonumber\\
			& <  k_2^n \left(1+\frac{4\,(2k_2-1)}{(k_2-1)^2(n+1)(n+5/2)}\right)	\nonumber\\
			&\qquad + k_2^n\left(4\,k_2^{-n}\frac{(\tilde{W}_1-\frac{1}{k_2-1})n+\left(\frac{1-2k_2}{(k_2-1)^2}+\tilde{W}_0\right)}{(n+1)(n+5/2)}\right)		\nonumber\\
			&<k_2^n\left(1+\frac{4\cdot1.4}{0.2^2\cdot 30\cdot 32.5}+\frac{4}{1.2^{30}}\,2.89\right)  \nonumber\\
			&<k_2^n (1+0.144+0.049) \nonumber\\
			&<k_2^{n+1}
\end{align}
where $\tilde{W}_1=271$ and $\tilde{W}_0=-5622$ are chosen such that 
\begin{equation}
\sum_{i=0}^{30}(n-i+1) (M_i-k_2^i)\,<\,\tilde{W}_1\,n +\tilde{W}_0
\end{equation}
and explicit calculation shows that for $n \geq 30$
\begin{equation}
\frac{(\tilde{W}_1-\frac{1}{k_2-1})n+\left(\frac{1-2k_2}{(k_2-1)^2}+\tilde{W}_0\right)}{(n+1)(n+5/2)} = \frac{266 n -5657}{(n+1)(n+5/2)} \leq 2.89
\end{equation}
ending the proof of Lemma \ref{JnMnLemma}.
\newpage

\subsection{Approximations $G_a$ and $H_a$}
\label{subsec:GaHa}
(1). $\lambda \in S_1$, expression of $G_a(t)$, $H_a(t)$.\\

For $t\in [0,1/4)$:
\begin{align*}
&G_a(t)  := \Bigg( {\frac {106}{111}}\,{\lambda}^{4}-{\frac {89}{86}}+{\frac {52}{113}}\,
i{\lambda}^{5}-{\frac {82}{89}}\,{\lambda}^{3}+{\frac {61}{36}}\,{
\lambda}^{2}
 \Bigg) {t}^{6}+ \Bigg(  {\frac {18}{181}}\,i{\lambda}^{3}-{\frac {79}{214}}\,{\lambda}^{5} \nonumber\\
 &-{
\frac {11}{75}}-{\frac {26}{259}}\,{\lambda}^{3}-{\frac {59}{362}}\,i{
\lambda}^{6}-{\frac {19}{88}}\,\lambda-{\frac {29}{93}}\,i{\lambda}^{4
}+{\frac {44}{205}}\,{\lambda}^{2}+{\frac {28}{283}}\,{\lambda}^{4}+{
\frac {58}{1195}}\,i{\lambda}^{5}
 \Bigg) {t}^{5}\nonumber\\
&+ \Bigg(-{\frac {5}{574}}\,i{\lambda}^{2}+{\frac {35}{96}}\,{\lambda}^{4}+{
\frac {23}{873}}\,{\lambda}^{6}+{\frac {259}{368}}\,{\lambda}^{2}-{
\frac {20}{43}}-{\frac {37}{450}}\,{\lambda}^{5}-{\frac {19}{527}}\,i{
\lambda}^{6}+{\frac {19}{864}}\,i{\lambda}^{3}\nonumber\\
&-{\frac {58}{161}}\,{
\lambda}^{3}-{\frac {2}{29}}\,i{\lambda}^{4}+{\frac {43}{243}}\,i{
\lambda}^{5}+{\frac {17}{749}}\,\lambda
 \Bigg) {t}^{4}
+ \Bigg( -{\frac {41}{90}}-{\frac {31}{257}}\,i{\lambda}^{4}+{\frac {22}{75}}\,
{\lambda}^{4}+{\frac {78}{127}}\,{\lambda}^{2}\nonumber\\
&+{\frac {271}{7045}}\,i{
\lambda}^{3}+{\frac {45}{314}}\,i{\lambda}^{5}-{\frac {33}{230}}\,{
\lambda}^{5}-{\frac {87}{293}}\,{\lambda}^{3}+{\frac {6}{713}}\,{
\lambda}^{6}-{\frac {3}{1081}}\,i{\lambda}^{2}-{\frac {7}{111}}\,i{
\lambda}^{6}+{\frac {43}{519}}\,\lambda
 \Bigg) {t}^{3}\nonumber\\
&+ \Bigg(  {\frac {79}{359}}\,{\lambda}^{4}+{\frac {8}{917}}\,{\lambda}^{6}+{
\frac {9}{83}}\,i{\lambda}^{5}-{\frac {4}{1389}}\,i{\lambda}^{2}-{
\frac {41}{402}}\,{\lambda}^{5}-\frac{1}{2}+{\frac {39}{125}}\,\lambda-{\frac 
{19}{424}}\,i{\lambda}^{6}\nonumber\\
&+{\frac {12}{439}}\,i{\lambda}^{3}+{\frac {
59}{110}}\,{\lambda}^{2}-{\frac {2}{7677}}\,i\lambda-{\frac {272}{3173
}}\,i{\lambda}^{4}-{\frac {93}{407}}\,{\lambda}^{3}
 \Bigg) {t}^{2}\nonumber\\
&+ \left( \frac{1}{4}\,{\lambda}^{2}+\frac{3}{4}\,\lambda-\frac{1}{2} \right) t+\lambda \nonumber\\
\end{align*}

For $t\in (1/4,1/2)$
\begin{align*}
&G_a(t):= \Bigg( -{\frac {121}{18}}\,i{\lambda}^{4}-{\frac {421}{77}}\,\lambda-{\frac {
183}{23}}\,{\lambda}^{5}-{\frac {1045}{58}}-17\,{\lambda}^{3}
 \Bigg) {t}^{7}+ \Bigg(  {\frac {361}{27}}\,i{\lambda}^{4}-{\frac {38}{17}}\,i{\lambda}^{6}+{\frac {76}{7}}\,\lambda \nonumber\\
 &+{\frac {657}{122}}\,i{\lambda}^{5}+{\frac {287}{15}}\,{\lambda}^{2}+{
\frac {236}{21}}\,{\lambda}^{4}+{\frac {491}{31}}\,{\lambda}^{5}+{
\frac {250}{7}}+{\frac {175}{129}}\,i{\lambda}^{3}+{\frac {72}{163}}\,
{\lambda}^{6}+{\frac {1387}{41}}\,{\lambda}^{3}
 \Bigg) {t}^{6}\nonumber\\
&+ \Bigg( -{\frac {313}{21}}\,{\lambda}^{5}-{\frac {29}{330}}\,i{\lambda}^{2}-{
\frac {327}{26}}\,i{\lambda}^{4}-{\frac {571}{31}}\,{\lambda}^{4}-{
\frac {1305}{41}}\,{\lambda}^{3}-{\frac {193}{19}}\,\lambda+{\frac {
216}{59}}\,i{\lambda}^{6}-{\frac {218}{7}}\,{\lambda}^{2}\nonumber\\
&-{\frac {238}{27}}\,i{\lambda}^{5}-{\frac {1178}{35}}-{\frac {62}{85}}\,{\lambda}^{
6}-{\frac {71}{32}}\,i{\lambda}^{3}
 \Bigg) {t}^{5}+ \Bigg({\frac {153}{22}}\,i{\lambda}^{5}+7/4\,i{\lambda}^{3}+{\frac {335}{64}
}\,\lambda+{\frac {423}{25}}\nonumber\\
\end{align*}
\begin{align*}
&+{\frac {145}{9}}\,{\lambda}^{3}+{\frac {
1381}{56}}\,{\lambda}^{2}+{\frac {76}{677}}\,i{\lambda}^{2}+{\frac {
151}{20}}\,{\lambda}^{5}+{\frac {116}{201}}\,{\lambda}^{6}+{\frac {247
}{17}}\,{\lambda}^{4}+{\frac {312}{49}}\,i{\lambda}^{4}-{\frac {26}{9}
}\,i{\lambda}^{6}
 \Bigg) {t}^{4}\nonumber\\
&+ \Bigg( -{\frac {181}{33}}\,{\lambda}^{3}-{\frac {93}{43}}\,i{\lambda}^{4}-{
\frac {113}{72}}\,\lambda-{\frac {136}{53}}\,{\lambda}^{5}+{\frac {65}
{59}}\,i{\lambda}^{6}-{\frac {61}{839}}\,i{\lambda}^{2}-{\frac {220}{
83}}\,i{\lambda}^{5}-{\frac {39}{176}}\,{\lambda}^{6}-{\frac {504}{85}
}\nonumber\\
&-{\frac {172}{31}}\,{\lambda}^{4}-{\frac {5}{3333}}\,i\lambda-{\frac 
{324}{35}}\,{\lambda}^{2}-{\frac {2543}{3815}}\,i{\lambda}^{3}
 \Bigg) {t}^{3}+ \Bigg({\frac {4}{223}}\,i{\lambda}^{2}+{\frac {530}{187}}\,{\lambda}^{2}+{
\frac {67}{107}}\,\lambda+{\frac {63}{83}}\,i{\lambda}^{5} \nonumber\\
 &+{\frac {17}
{89}}\,i{\lambda}^{3}+{\frac {34}{45}}\,{\lambda}^{3}+{\frac {71}{132}
}+{\frac {1}{9639}}\,i-{\frac {35}{111}}\,i{\lambda}^{6}+{\frac {264}{
167}}\,{\lambda}^{4}+{\frac {14}{39}}\,{\lambda}^{5}+{\frac {8}{127}}
\,{\lambda}^{6}+{\frac {56}{185}}\,i{\lambda}^{4}\nonumber\\
 &+{\frac {1}{1180}}\,i\lambda
 \Bigg) {t}^{2}+ \Bigg( {\frac {32}{953}}\,i{\lambda}^{6}-{\frac {18}{511}}\,{\lambda}^{2}-{
\frac {2}{7333}}\,i\lambda-{\frac {23}{222}}\,{\lambda}^{3}-{\frac {17
}{837}}\,i{\lambda}^{3}-{\frac {2}{637}}\,i{\lambda}^{2}-{\frac {133}{
787}}\,{\lambda}^{4}\nonumber\\
 &-{\frac {5}{739}}\,{\lambda}^{6}-{\frac {9}{220}}
\,i{\lambda}^{4}-{\frac {12}{247}}\,{\lambda}^{5}-{\frac {53}{656}}\,i
{\lambda}^{5}+{\frac {38}{53}}\,\lambda-{\frac {92}{151}}-{\frac {1}{
19278}}\,i
 \Bigg) t- {\frac {673}{672}}\,\lambda+{\frac {3}{1631}}\,i{\lambda}^{4}\nonumber\\
 &+{\frac {
8}{1717}}\,{\lambda}^{3}-{\frac {5}{2888}}\,i{\lambda}^{6}+{\frac {1}{
955}}\,i{\lambda}^{3}+{\frac {1}{44699}}\,i\lambda+{\frac {1}{240}}\,i
{\lambda}^{5}+{\frac {15}{3052}}+{\frac {1}{154231}}\,i+{\frac {1}{
5265}}\,i{\lambda}^{2}\nonumber\\
 &+{\frac {11}{1262}}\,{\lambda}^{4}+{\frac {3}{
1372}}\,{\lambda}^{5}+{\frac {2}{5759}}\,{\lambda}^{6}+{\frac {110}{
7481}}\,{\lambda}^{2}
\nonumber\\
\end{align*}

For $t \in (1/2,3/4)$,let $s=1-t$
\begin{align*}
&H_a(1-s) :=\Bigg({\frac {64}{67}}\,i{\lambda}^{5}-{\frac {95}{158}}\,i{\lambda}^{6}-{
\frac {205}{461}}\,i{\lambda}^{4}-{\frac {17}{13}}
\Bigg){s}^{6}+ \Bigg( {\frac {148}{61}}+{\frac {19}{305}}\,{\lambda}^{4}+{\frac {10}{199}}\,
\lambda-{\frac {69}{35}}\,i{\lambda}^{5}\nonumber\\
		&+{\frac {98}{107}}\,i{\lambda}
^{4}+{\frac {19}{206}}\,{\lambda}^{5}+{\frac {23}{348}}\,{\lambda}^{6}
+{\frac {41}{33}}\,i{\lambda}^{6}-{\frac {28}{485}}\,{\lambda}^{2}
 \Bigg) {s}^{5}+ \Bigg({\frac {182}{109}}\,i{\lambda}^{5}+{\frac {24}{583}}\,{\lambda}^{3}+{
\frac {31}{408}}\,{\lambda}^{2}\nonumber\\
		&-{\frac {159}{151}}\,i{\lambda}^{6}-{
\frac {14}{135}}\,{\lambda}^{6}-{\frac {17}{74}}\,\lambda-{\frac {118}
{77}}-{\frac {42}{605}}\,{\lambda}^{4}-{\frac {31}{40}}\,i{\lambda}^{4
}-{\frac {41}{302}}\,{\lambda}^{5}
 \Bigg) {s}^{4}+ \Bigg( {\frac {9}{134}}\,{\lambda}^{6}+{\frac {61}{130}}\,i{\lambda}^{6}\nonumber\\
		&-{
\frac {33}{236}}+{\frac {60}{193}}\,\lambda-{\frac {24}{593}}\,{
\lambda}^{2}-{\frac {119}{160}}\,i{\lambda}^{5}+{\frac {30}{569}}\,{
\lambda}^{4}-{\frac {11}{602}}\,{\lambda}^{3}+{\frac {99}{287}}\,i{
\lambda}^{4}+{\frac {31}{350}}\,{\lambda}^{5}
 \Bigg) {s}^{3}+ \Bigg( -{\frac {233}{599}}\,\lambda\nonumber\\
		&-{\frac {83}{975}}\,i{\lambda}^{4}-{\frac 
{15}{697}}\,{\lambda}^{6}-{\frac {11}{399}}\,{\lambda}^{5}+{\frac {47}
{256}}\,i{\lambda}^{5}+{\frac {24}{535}}\,{\lambda}^{3}+{\frac {9}{359
}}\,{\lambda}^{2}+{\frac {68}{65}}-{\frac {49}{423}}\,i{\lambda}^{6}-{
\frac {1}{102}}\,{\lambda}^{4}
 \Bigg) {s}^{2}\nonumber\\
		&+ \Bigg( {\frac {193}{643}}\,\lambda-{\frac {5}{2957}}\,{\lambda}^{3}-{\frac {
79}{58}}+{\frac {9}{92}}\,{\lambda}^{2}+{\frac {11}{733}}\,i{\lambda}^
{6}+{\frac {5}{1156}}\,{\lambda}^{5}-{\frac {75}{3151}}\,i{\lambda}^{5
}+{\frac {4}{1967}}\,{\lambda}^{4}+{\frac {5}{453}}\,i{\lambda}^{4}+\nonumber\\
		&{
\frac {5}{1458}}\,{\lambda}^{6}
 \Bigg) s-{\frac {4}{6811}}\,i{\lambda}^{4}-{\frac {1}{24270}}\,\lambda-{\frac 
{5}{6274}}\,i{\lambda}^{6}+{\frac {10}{7901}}\,i{\lambda}^{5}-{\frac {
2}{7433}}\,{\lambda}^{5}-{\frac {1}{4713}}\,{\lambda}^{6}+{\frac {498}
{499}}\nonumber\\
		&+{\frac {1}{8478}}\,{\lambda}^{3}+{\frac {1}{7548}}\,{\lambda}^{
2}-{\frac {1}{8040}}\,{\lambda}^{4}
\nonumber\\
\end{align*}

For $t\in(3/4,1]$, let $s=1-t$,
\begin{align*}
&H_a(1-s) := \Bigg( {\frac {13}{76}}\,\lambda-{\frac {201}{286}}
 \Bigg) {s}^{5}+ \Bigg({\frac {11}{2267}}\,{\lambda}^{5}+{\frac {13}{1097}}\,{\lambda}^{4}+{
\frac {21}{1000}}\,{\lambda}^{3}-{\frac {17}{3153}}\,{\lambda}^{2}-{
\frac {21}{62}}\,\lambda \nonumber\\
&+{\frac {75}{62}}
 \Bigg) {s}^{4}+ \Bigg( {\frac {1}{1267}}\,{\lambda}^{5}+{\frac {11}{1538}}\,{\lambda}^{4}+{
\frac {5}{1326}}\,{\lambda}^{3}+{\frac {7}{953}}\,{\lambda}^{2}+{
\frac {85}{246}}\,\lambda-{\frac {346}{255}}
 \Bigg) {s}^{3}+ \Bigg( {\frac {1}{6961}}\,{\lambda}^{6}\nonumber\\
 &+{\frac {1}{37023}}\,{\lambda}^{5}+{
\frac {4}{1103}}\,{\lambda}^{4}+{\frac {60}{1679}}\,{\lambda}^{3}+{
\frac {17}{1590}}\,{\lambda}^{2}-{\frac {152}{387}}\,\lambda+{\frac {
98}{73}} \Bigg) {s}^{2}\nonumber\\
		&+ \Bigg( 1/10\,{\lambda}^{2}+3/10\,\lambda-7/5
 \Bigg) s+1\nonumber\\
\end{align*}

\subsection {Estimate of modulus of some rational functions on the boundary of  a rectangle.}
\label{rationalS2} 
\subsubsection {A description of the method}
\label{rationalme}
  Here we describe a method to obtain an upper bound and a lower bound of $|F(\lambda)|$ on the boundary of a rectangle $$R:=\{z:\Re(z) \in[0,1/2], \Im(z) \in [a,b]\}$$ where
\begin{equation}
\label{rationalfn}
F(\lambda)=P(\lambda)+\sum_{j=0}^{n-1}\, \frac{\tilde{a}_j}{\lambda-s_j}
\end{equation} 
where $P(\lambda)$ is a polynomial of degree $m$, ($m\leq 12$), $s_j$ are distinct complex numbers outside of R, and $\tilde{a}_j$'s are constants. The estimates of $F(\lambda)$ are obtained in 4 steps.\\

(i) Define a partition $\mathcal{P}$ on the boundary $\partial S_2$, given by $\{d_i\}_{0\leq i \leq N-1}$, satisfying
\begin{equation}
d_0 \in \{a\,i,\,a\,i+1/2,\,b\, i+1/2,\,b\,i\} \subseteq \mathcal{P}
\end{equation}
$d_i$'s are ordered counterclockwise. Let $d_N = d_0$.\\
For each $i \in \{1,2,...,N\}$, denote the disk that has $[d_{i-1},d_{i}]$ as a diameter by $D_i$, the center of $D_i$ by $c_i$. Also let $r_i = \frac{d_i - d_{i-1}}{2}$. Since all vertices of $S_2$ are in $\mathcal{P}$, each $[d_{i-1},d_{i}]$ either lies on a horizontal side or a vertical side of $\partial S_2$. In the following steps we will find an upper bound and a lower bound of $|F(\lambda)|$ on each $[d_{i-1},d_{i}]$.\\

Now fix $i \in \{1,2,...,N\}$.\\

(ii) Make a change of variable 
\begin{equation}
\label{cov}
\lambda_i(x) = c_i + r_i \,x
\end{equation}
$\lambda_i$ is a bijection that maps the unit disk $\{|x|\leq 1\}$ to $D_i$, and the image of $[-1,1]$ under $\lambda_i$ is $[d_{i-1},d_{i}]$. Consider the Taylor expansion of $F(\lambda_i(x))$:
\begin{align}
&\,\,F(\lambda_i(x)) \nonumber\\
&= P(c_i+r_i\,x) +\sum_{j=0}^{n-1}\, \frac{\tilde{a}_j}{c_i+r_i\,x-s_j}\nonumber\\
&= P(c_i+r_i\,x) +\sum_{j=0}^{n-1}\, \frac{\tilde{a}_j}{c_i-s_j}\,\frac{1}{1+\frac{r_i\,x}{c_i-s_j}}\nonumber\\
&= P(c_i+r_i\,x) +\sum_{j=0}^{n-1}\, \frac{\tilde{a}_j}{c_i-s_j}\,\sum_{k=0}^{m_0}\left(-\frac{r_i\,x}{c_i-s_j}\right)^k+\sum_{j=0}^{n-1}\, \frac{\tilde{a}_j}{c_i-s_j}\,\frac{\left(-\frac{r_i\,x}{c_i-s_j}\right)^{m_0+1}}{1+\frac{r_i\,x}{c_i-s_j}}
\end{align}

Let
\begin{equation}
\tilde{F}_{i,1}(x) = P(c_i+r_i\,x) +\sum_{j=0}^{n-1}\, \frac{\tilde{a}_j}{c_i-s_j}\,\sum_{k=0}^{m_0}\left(-\frac{r_i\,x}{c_i-s_j}\right)^k = \sum_{k=0}^M \tilde{b}_{k,i}\,x^k
\end{equation}

\begin{equation}
E_{i,1} > \sum_{j=0}^{n-1}\, \left\vert \frac{\tilde{a}_j}{c_i-s_j}\right\vert\,\frac{\left\vert\frac{r_i}{c_i-s_j}\right\vert^{m_0+1}}{1-\left\vert\frac{r_i\,x}{c_i-s_j}\right\vert}
\end{equation}
where $M = \max(m,m_0)$, then
\begin{equation}
\label{ratlink1}
\left\vert F(\lambda_i(x))-\tilde{F}_{i,1}(x)\right\vert < E_{i,1} 	\quad( |x| \leq 1)
\end{equation}
Repacing each coefficient $\tilde{b}_{k,i}$ by a rational number $\tilde{c}_{k,i}$ within $e^{-9}$ of it we obtain
\begin{equation}
F_{i,1}(x) = \sum_{k=0}^{M} \tilde{c}_{k,i}\,x^k
\end{equation}
Let
\begin{equation}
E_{i,2}>\sum_{k=0}^M\left\vert \tilde{b}_{k,i}-\tilde{c}_{k,i}\right\vert
\end{equation}
Then
\begin{equation}
\label{ratlink2}
\left\vert F_{i,1}(x)-\tilde{F}_{i,1}(x) \right\vert < E_{i,2} \quad\quad( |x| \leq 1)
\end{equation}

(iii) Consider $x\in [-1,1]$. We estimate $\left\vert F_{i,1}(x) \right\vert^2$ using Chebyshev polynomials of the first kind.
\begin{equation}
\left\vert F_{i,1}(x) \right\vert^2 = \left[ \Re(F_{i,1}(x))\right]^2 +\left[\Im(F_{i,1}(x))\right]^2
\end{equation}
$\left\vert F_{i,1}(x) \right\vert^2$ is a polynomial of degree $2M$. Express each power of $x$ by a linear combination of Chebyshev polynomials $T_k(x)$ then we have a nonnegative polynomial:
\begin{equation}
\left\vert F_{i,1}(x) \right\vert^2 = \sum_{k=0}^{2M} \tilde{e}_{k,i} \,T_k(x)
\end{equation}
For $k=0,1,2,3$ let $\tilde{f}_{k,i}$ be a rational number within $e^{-9}$ of it and let
\begin{equation}
A_i(x) = \sum_{k=0}^3 \tilde{f}_{k,i} T_k(x)
\end{equation}
\begin{equation}
E_{i,3} > \sum_{k=0}^3 \left\vert \tilde{e}_{k,i}-\tilde{f}_{k,i} \right\vert 
+\sum_{k=4}^{2M} \left\vert\tilde{e}_{k,i} \right\vert
\end{equation}
Then
\begin{equation}
\label{ratlink3}
\left\vert |F_{i,1}(x)|^2-A_i(x) \right\vert < E_{i,3} \quad\quad x\in[-1,1]
\end{equation}
Now $A_i(x)$ is a cubic polynomial, so we can calculate its maximum $U_{i,1}$ and minimum $L_{i,1}$ on $[-1,1]$,
\begin{equation}
\label{ratlink4}
L_{i,1}\leq A_i (x) \leq U_{i,1}
\end{equation}
In cases we deal with in this paper, we can always choose $E_{i,l}$ ($l=1,2,3$), $U_i$ and $L_i$ such that
$$\sqrt{L_{i,1}-E_{i,3}} >>E_{i,1}+E_{i,2}$$
so from \eqref{ratlink3} and \eqref{ratlink4} we get:
\begin{equation}
\sqrt{L_{i,1}-E_{i,3}} \leq \left\vert F_{i,1}(x) \right\vert \leq \sqrt{U_{i,1}+E_{i,3}}
\end{equation}
Furthermore, from \eqref{ratlink1} and \eqref{ratlink2} we get:
\begin{equation}
\left\vert F_{i,1}(x)-F(\lambda_i(x))\right\vert \leq E_{i,1}+E_{i,2}
\end{equation}
Thus if for each i we choose
\begin{align}
U_i &> \sqrt{U_{i,1}+E_{i,3}} + E_{i,1}+E_{i,2} \\
L_i &< \sqrt{L_{i,1}-E_{i,3}} - E_{i,1} - E_{i,2}
\end{align}
then
\begin{equation}
L_i<\left\vert F(\lambda_i(x))\right\vert <U_i 	\quad\quad x\in[-1,1]
\end{equation}
i.e.
\begin{equation}
\quad L_i<\left\vert F(\lambda)\right\vert <U_i 	\quad\quad \lambda \in[d_{i-1},d_{i}]
\end{equation}
Finally let
\begin{align}
U&= \max_{1\leq i \leq N}{U_i}\\
L&= \min_{1\leq i \leq N}{L_i}
\end{align}
Then
\begin{equation}
L<\left\vert F(\lambda)\right\vert <U\quad\quad(\lambda \in \partial{R})
\end{equation}

\subsubsection{Derivation of $J_n$, $P_n$, upper bounds of $|G_a(t_0)|$ and $|G'_a(t_0)|$ on $S_2$ and a lower bound of $|W_a(t_0)|$ on $\partial S_2$.} Each one of $a_n$, ($1\leq n\leq 10$), $G_a(t_0)$, $G'_a(t_0)$ and $W_a(t_0)$ is a rational function of the form in \eqref{rationalfn} with $R = S_2$. We use the method described in \S\ref{rationalme} to obtain upper bounds of $|a_n|$, ($1\leq n\leq 10$), $|G_a(t_0)|$, $|G'_a(t_0)|$ over $\partial S_2$ and a lower bound of $|W_a(t_0)|$ on $\partial S_2$.

Two partitions are used
\begin{align}
\mathcal{P}_1 = \{d_i^{(1)}\}_{0\leq i \leq 8} \\
\mathcal{P}_2 = \{d_i^{(2)}\}_{0\leq i \leq 8}
\end{align}
where 
$$d_0^{(1)} = \frac{1}{2}+4\, i,\,\,\, d_1^{(1)} = \frac{1}{4}+4\, i,\,\,\, d_2^{(1)} = 4\, i,\,\,\, d_3^{(1)} = \frac{17}{6}\, i,\,\,\,d_4^{(1)} = \frac{5}{3}\, i,$$
$$d_5^{(1)} = \frac{1}{2}\, i,\,\,\, d_6^{(1)} = \frac{1}{2}+\frac{1}{2}\, i,\,\,\, d_8^{(1)} = \frac{1}{2}+\frac{17}{8}\, i,\,\,\,d_9^{(1)} = \frac{1}{2}+\frac{15}{4}\, i$$
and $d_i^{(1)}+d_i^{(2)} = \frac{9}{2}\,i+\frac{1}{2}$ for each $i$.

To obtain upper bounds of $|a_n|$, ($1\leq n\leq 6$) we used partition $\mathcal{P}_1$. To obtain upper bounds of $|a_n|$, ($7\leq n\leq 10$), $|p_n|$ ($11\leq n\leq 29$), $|G_a(t_0)|$ and $|G'_a(t_0)|$ we used partition $\mathcal{P}_2$.

In all these calculations we chose to keep 11 terms in the Taylor expansions in Step (i), i.e. $m_0 = 10$.

Note that each one of $a_n$'s, $G_a(t_0)$ and $G'_a(t_0)$ is analytic in $S_2$, so an upper bound of its modulus on the boundary $\partial S_2$ is also an upper bound on whole $S_2$.\\

To obtain a lower bound of $W_a(t_0)$ on $S_2$ we added a couple of steps.\\
 \indent First approximate $W_a$ by a rational function $F(\lambda)$ which can be estimated by the method in \S\ref{rationalme}: 
\begin{equation}
W_a(t_0) = F(\lambda) + \lambda^{13} \, \sum_{k=0}^{15} \tilde{g}_{i} \lambda^i
\end{equation}
where $F(\lambda)$ is a rational function as described in \S\ref{rationalme}, and fortunately by explicit calculation
\begin{equation}
\label{polye}
|\lambda|^{13}\,\sum_{k=0}^{15} |\tilde{g}_{i}|\, |\lambda|^i < 0.079 	\quad\quad (|\lambda| \leq |4\,i+1/2|)
\end{equation}
Use the method in \S\ref{rationalme} and choose the partition to be $\mathcal{P}_2$, $m_0=10$, we obtain:
\begin{equation}
|F(\lambda) | > 1.14 \quad\quad (\lambda \in \partial S_2)
\end{equation}
Thus
\begin{equation}
\label{Washarp}
|W_a(t_0) | > 1.14-0.08=1.06 \quad\quad (\lambda \in \partial S_2)
\end{equation}
Finally, to show a lower bound of $|W_a(t_0)|$ on the boundary $\partial S_2$ is also a lower bound on $S_2$ it suffices to show that $W_a(t_0)$ has no roots in $S_2$, see \S\ref{WaS2} for the proof.

\subsection {Estimate of modulus of some polynomials on the boundary of  $S_2$}
\label{polynS2}
\subsubsection{A description of the method}
In this subsection we describe the method used to estimate $|P(\lambda)|$ on the boundary of $S_2$, where $P(\lambda)$ can be any one of $b_n$ ($1\leq n\leq 10$), $H_a(t_0)$ or $H'_a(t_0)$. Denote the degree of $P(\lambda)$ by $m$. The method is similar to the one in \S\ref{rationalS2}. Slight modifications are made.\\
(i) Let  $\mathcal{P}$, $ \{d_i\}_{0\leq i\leq N}$ and $ \{c_i\}_{1\leq i\leq N}$  all be the same as in \S\ref{rationalme} \\

Fix $i \in \{1,2,...,N\}$.\\

(ii) Make a change of variable \eqref{cov}, then
\begin{equation}
P(\lambda_i(x)) = \sum_{k=0}^{m} \tilde{a}_{k,i}\,x^k
\end{equation}
For each $k$ express $x^k$ by a linear combination of Chebyshev polynomials $T_k(x)$ then we have
\begin{equation}
P(\lambda_i(x)) = \sum_{k=0}^{m} \tilde{b}_{k,i}\,T_k(x)
\end{equation}
 Let
 \begin{align}
P_{i,1}(x) &= \sum_{k=0}^{5} \tilde{b}_{k,i}\,T_k(x)\\
E_{i,1} &> \sum_{k=6}^{m} \left\vert\tilde{b}_{k,i}\right\vert
\end{align}
Then
 \begin{equation}
 \label{polylink1}
\left\vert P(\lambda_i(x)) - P_{i,1}(x)\right\vert < E_{i,1} \quad\quad x\in [-1,1]
\end{equation}

(iii) Estimate the polynomial $\left\vert P_{i,1}(x)\right\vert^2 $ which is of degree at most 10.
\begin{equation}
\left\vert P_{i,1}(x)\right\vert^2 = \sum_{i=0}^{10} \tilde{c}_{k,i}\, x^k
\end{equation}
Express each power of $x$ by a linear combination of Chebyshev polynomials $T_k(x)$ then we have
\begin{equation}
\left\vert P_{i,1}(x)\right\vert^2 = \sum_{i=0}^{10} \tilde{d}_{k,i}\, T_k(x)
\end{equation}
Let
\begin{align}
A_i(x) &= \sum_{k=0}^3 \tilde{e}_{k,i} \, T_k(x)\\
\left\vert E_{i,2}\right\vert &> \sum_{k=0}^3\left\vert \tilde{d}_{k,i}-\tilde{e}_{k,i}\right\vert +\sum_{k=4}^{10}\left\vert \tilde{d}_{k,i}\right\vert
\end{align}
where each $\tilde{e}_{k,i}$ is a rational number within $e^{-7}$ of $\tilde{d}_{k,i}$, then
\begin{equation}
 \label{polylink2}
\left\vert |P_{i,1}(x)|^2- A_i(x)\right\vert < E_{i,2} \quad\quad x\in [-1,1]
\end{equation}
Now we can determine the maximum $U_{i,1}$ and the minimum $L_{i,1}$ of the cubic polynomial $A_i$ on $[-1,1]$. We obtain
\begin{equation}
 \label{polylink3}
 L_{i,1} \leq A_i(x)\leq U_{i,1} \quad\quad x\in [-1,1]
\end{equation}
From \eqref{polylink2} and \eqref{polylink3} we get
\begin{equation}
\sqrt{L_{i,1}-E_{i,2}} \leq \left\vert P_{i,1}(x) \right\vert \leq \sqrt{U_{i,1}+E_{i,2}}
\end{equation}
and in view of \eqref{polylink1} if we choose for each i:
\begin{align}
U_i &> \sqrt{U_{i,1}+E_{i,2}} +E_{i,1} \\
L_i &< \sqrt{L_{i,1}-E_{i,2}} - E_{i,1} 
\end{align}
then
\begin{equation}
L_i<\left\vert P(\lambda_i(x))\right\vert <U_i 	\quad\quad x\in[-1,1]
\end{equation}
i.e.
\begin{equation}
L_i<\left\vert P(\lambda)\right\vert <U_i 	\quad\quad \lambda \in[d_{i-1},d_{i}]
\end{equation}
Let
\begin{align}
U&= \max_{1\leq i \leq N}{U_i}\\
L&= \min_{1\leq i \leq N}{L_i}
\end{align}
Then
\begin{equation}
L<\left\vert P(\lambda)\right\vert <U \quad\quad \lambda \in \partial S_2
\end{equation}

\subsubsection{Derivation of $M_n$, upper bounds of $H_a(t_0)$ and $H'_a(t_0)$}
Each one of $b_n$, $H_a(t_0)$ and $H'_a(t_0)$ is a polynomial. To obtain the upper bounds of their moduli, 2 partitions are used
\begin{align}
\mathcal{P}_3 = \{d_i^{(3)}\}_{0\leq i \leq 6} \\
\mathcal{P}_4 = \{d_i^{(4)}\}_{0\leq i \leq 6}
\end{align}
where 
$$d_0^{(3)} = \frac{i}{2},\,\,\, d_1^{(3)} = \frac{i}{2}+\frac{1}{2},\,\,\,d_2^{(3)} = \frac{7\,i}{2}+\frac{1}{2} $$
$$d_3^{(3)} = \frac{1}{2}+4\, i,\,\,\,d_4^{(3)} = \frac{1}{4}+4\, i,\,\,\, d_5^{(3)} = 4\, i,\,\,\, d_6^{(3)} = \frac{7}{2}\, i,$$
and $d_i^{(3)}+d_i^{(4)} = \frac{9}{2}\,i+\frac{1}{2}$ for each $i$.

To obtain upper bounds of $|b_n|$ and $|H'_a(t_0)|$, ($1\leq n\leq 10$) we used partition $\mathcal{P}_3$. To obtain upper bounds of $|H_a(t_0)|$ we used partition $\mathcal{P}_4$.

\subsection{Estimate of $W(G_a,H_a)(t_0)$ in $S_2$}
\label{WaS2} 
Equally divide $S_2$ into 7 squares:
\begin{equation}
S_{2,k} := \{z: 0\leq \Re(z) \leq \frac{1}{2},\,  \Im(C_k)-\frac{1}{4} \leq \Im(z) \leq \Im(C_k)+\frac{1}{4} \}
\end{equation}
where for each $k \in \{1, 2, ...,7\}$
$$C_k = \frac{3\,i}{4}+\frac{1}{4}+(k-1)\,\frac{i}{2}$$ 
is the center of the square $S_{2,k}$ .
$W_a(t_0):=W(G_a,H_a)(t_0)$ is of the form
\begin{equation}
W_a(t_0) = \sum_{j=0}^7\frac{\tilde{a}_j}{\lambda+j} +  \sum_{k=0}^{12} \tilde{b}_k\, \lambda^k + \lambda^{13} \, \sum_{k=0}^{15} \tilde{g}_{k} \lambda^k
\end{equation}
In view of \eqref{polye}, if we let 
\begin{equation}
W_1(\lambda) = \sum_{j=0}^7\frac{\tilde{a}_j}{\lambda+j} +  \sum_{k=0}^{12} \tilde{b}_k\, \lambda^k
\end{equation}
then
\begin{equation}
\label{Walink1}
 |W_a(t_0)-W_1(\lambda)| <0.079 \quad\quad (\lambda \in S_2)
\end{equation}
For each i, let
\begin{equation}
\bar{A}_i = W_1(C_i)
\end{equation}
then for $\lambda \in S_{2,i}$:
\begin{align}
\label{Walink2}
&\left\vert W_1(\lambda) - \bar{A}_i \right\vert = \left\vert \sum_{j=0}^7\frac{\tilde{a}_j}{\lambda+j} - \sum_{j=0}^7\frac{\tilde{a}_j}{C_i+j} +  \sum_{k=0}^{12} \tilde{b}_j\, \lambda^k - \sum_{k=0}^{12} \tilde{b}_k\, C_i^k\right\vert \nonumber\\
&\leq \sum_{j=0}^7\frac{|\tilde{a}_j|\,|C_i - \lambda|}{|\lambda+j |\,|C_i+j|} + \left\vert \sum_{k=0}^{12} \left(\tilde{b}_k - C_i^k\right)\right\vert
\end{align}
For $\lambda \in S_{2,i}$, 
\begin{align}
\left\vert C_i-\lambda \right\vert &\leq \frac{\sqrt{2}\,\,}{4}; \\
\left\vert \lambda+j \right\vert &\geq \left\vert C_i - \frac{1}{4} -\frac{i}{4}+j \right\vert
\end{align}
so
\begin{equation}
\label{Walink3}
\sum_{j=0}^7\frac{|\tilde{a}_j|\,|C_i - \lambda|}{|\lambda+j |\,|C_i+j|} \leq \sum_{j=0}^7\frac{|\tilde{a}_j|\,\sqrt{2}/4}{|C_i - \frac{1}{4} -\frac{i}{4}+j |\,|C_i+j|} :=R_{i,1}
\end{equation}
On the other hand
\begin{equation}
\sum_{k=0}^{12} \left(\tilde{b}_k - C_i^k\right) = \sum_{k=1}^{12} \tilde{c}_{k,i}\left(\lambda - C_i\right)^k
\end{equation}
so by explicit calculation,
\begin{align}
\label{Walink4}
\left\vert \sum_{k=0}^{12} \left(\tilde{b}_k - C_i^k\right) \right\vert \leq \sum_{k=1}^{12} |\tilde{c}_{k,i}|\left\vert \lambda - C_i\right\vert ^k
\leq \sum_{k=1}^{12} |\tilde{c}_{k,i}| \left(\frac{\sqrt{2}}{4}\right)^k := R_{i,2}
\end{align}

Now choose $E_i$ such that
\begin{equation}
\label{Walink5}
E_i > R_{i,1}+R_{i,2}
\end{equation}
Combine \eqref{Walink2}, \eqref{Walink3}, \eqref{Walink4} and \eqref{Walink5} we have
\begin{equation}
\label{Walink6}
\left\vert W_1(\lambda) - \bar{A}_i\right\vert < E_i \quad\quad (\lambda \in S_{2,i})
\end{equation}
Combine \eqref{Walink6} with \eqref{Walink1} we obtain:
\begin{equation}
\label{Walink7}
|W_a(t_0)| > \bar{A}_i - (E_i +0.079)
\end{equation}
From straightforward calculation we have the following:
\begin{align}
\label{barA}
&|\bar{A}_1| > 1.95, \quad|\bar{A}_2| > 1.80,\quad |\bar{A}_3| > 1.74,\quad |\bar{A}_4| > 1.68, \nonumber\\
&|\bar{A}_5| > 1.61,\quad |\bar{A}_6| > 1.54,\quad |\bar{A}_7| > 1.48,      
\end{align}
According to calculations for $R_{i,1}$ and $R_{i,2}$, we can choose valid $E_i$ as follows:
\begin{align}
\label{errWa}
&|E_1| =1.28, \quad|E_2| > 0.77,\quad |E_3| =0.67,\quad |E_4| = 0.65, \nonumber\\
&|E_5| =0.66,\quad |E_6| =0.67,\quad |E_7| =0.65,      
\end{align}
Hence from \eqref{Walink7}, \eqref{barA} and \eqref{errWa} we have
\begin{equation}
|W_a(t_0)| > \min_{1\leq i \leq 7} \{\bar{A}_i - (E_i + 0.079)\} > 0.59
\end{equation}
This crude estimate is enough for the proof of Proposition \ref{PropnS2}, but we are able to obtain a better estimate once we know that that $W_a(t_0)$ has no roots in $S_2$ and thus the minimum of $|W_a(t_0)|$ on $S_2$ is attained on the boundary. The method in \S\ref{rationalS2} provides us with a sharper result \eqref{Washarp}.

\newpage
\subsection {Values of $J_n$, $P_n$ and $M_n$} 
\label{valJnMn}
\begin{center}
    \begin{tabular}{| l ||  l  |  l  |  l  |}
    \hline
    n	&		$10^7\,J_n$			&	$P_n$	&		$10^7\,M_n$		 \\ \hline
    0 	& $10000000$		 					&	  				& 10000000 						\\ \hline
   1	     &$14142986$	               &				&		32466329						\\\hline
    2	&$14385562$	 &	$ $			&	 43849481				\\\hline
    3	&$13584388$	&	$ $			&		25844183						\\\hline
    4	&$12609618$   &	$ $			&		22641590								\\\hline
    5	&$11709327$	&	$ $			&		34192044							\\\hline
    6	&$10953644$	&	$ $			&		18884907						\\\hline
    7	&$10383616$	&	$ $			&	29748748						\\\hline
    8	&10372084	&	$ $			&		25563511				\\\hline
    9	&10432493		&	$ $			&		22348202						\\\hline
    10	&10540865	&	$ 2847/2737$			&		27433266				\\\hline
    11	&11256315	&	$1975/1904$			&		74641000						\\\hline
    12	&11917388	&	$5417/5235$			&		117161600								\\\hline
    13	&12539695	&	$793/768$			&	158112600						\\\hline
    14	&13132397	&	$1381/1340$			&	200013200				\\\hline
    15	&13700995	&	$1533/1490$			&		244928600			\\\hline
    16	&14248989		&	$2897/2820$			&		294605900						\\\hline
    17	&14778796	&	$20077/19570$			&	350585000					\\\hline
    18	&15292210	&	$1871/1826$			&	414285300				\\\hline
    19	&15790675	&	$6316/6171$			&	487071000				\\\hline
    20	&16275383	&	$8243/8062$			&	570301600		\\\hline
    21	&16747362	&	$5323/5211$			&	665370000			\\\hline
    22	&17207507	&	$3763/3687$			&	773733100			\\\hline
    23	&17656621	&	$103/101$			&	896936200					\\\hline
    24	&18095424	&	$2193/2152$			&	1036633000 				\\\hline
    25	&18524561	&	$1276/1253$			&	1194603700					\\\hline
    26	&18944602	&	$1839/1807$			&	1372770700 				\\\hline
    27	&19356079	&	$3627/3566$			&	1573213300 			\\\hline
    28	&19759479	&	$2581/2539$			&	1798182200 					\\\hline
    29	&20155227	&	$2157/2123$			&	2050113400	 			\\\hline
    30	&20543752	&	$1047/1031$			&	2331642500 		\\
    \hline
    \end{tabular}\\
\end{center}

\subsection{Estimates of $|B_n/A_n|$, $|B_n/A_n|$ and $|B_n/A^2_n|$.} \label{S_3}
\subsubsection{Upper bound of $|B_n/A_n|$}\label{S_31}
Let $\lambda= x+iy$, then 
\begin{equation}
\label{BA}
F_1(x,y,n):= \frac{B_n(x+iy)}{A_n(x+iy)}
\end{equation}
is a rational function in $x$, $y$ and $n$. Let $z=1/n$, then
\begin{align}
M_1(x,y,z):=|F_1(x,y,1/z)|^2 = \frac{P_1(x,y,z)}{Q_1(x,y,z)}
\end{align}
is a real valued rational function, where
\begin{align}
P_1(x,y,z)&=\left( {x}^{2}{z}^{2}+{z}^{2}{y}^{2}-6\,x{z}^{2}+4\,xz+9\,{z}^{2}-12
\,z+4 \right) \cdot\nonumber\\
&\left( {x}^{2}{z}^{2}+{z}^{2}{y}^{2}+6\,x{z}^{2}+4\,xz+9\,{z}^{2}+12\,z+4 \right) \\
Q_1(x,y,z)&={x}^{4}{z}^{4}+2\,{x}^{2}{y}^{2}{z}^{4}+{y}^{4}{z}^{4}+16\,{x}^{3}{z}^
{4}+16\,x{y}^{2}{z}^{4}+16\,{x}^{3}{z}^{3}+46\,{x}^{2}{z}^{4}\nonumber\\
&+16\,x{y}^{2}{z}^{3}+82\,{y}^{2}{z}^{4}+168\,{x}^{2}{z}^{3}-144\,x{z}^{4}+88\,{
y}^{2}{z}^{3}+88\,{x}^{2}{z}^{2}+176\,x{z}^{3}\nonumber\\
&+40\,{z}^{2}{y}^{2}+81\,{z}^{4}+512\,x{z}^{2}-360\,{z}^{3}+192\,xz+184\,{z}^{2}+480\,z+144\\
\end{align}
Take partial derivative of $M_1(x,y,z)$ with respect to x:
\begin{equation}
\frac{\partial}{\partial x}M_1(x,y,z) = z \frac{P_2(x,y,z)}{Q_2(x,y,z)}
\end{equation}
where
\begin{align}
P_2(x,y,z)&=1536+ \left( 4096\,x+7168 \right) z+ \left( 4480\,{x}^{2}+640\,{y}^{2}
+17664\,x+6528 \right) {z}^{2}\nonumber\\
                &+ \left( 2560\,{x}^{3}+1024\,x{y}^{2}+
17664\,{x}^{2}+2304\,{y}^{2}+14848\,x-6912 \right) {z}^{3}\nonumber\\
                &+ \big( 800\,{x}^{4}+576\,{x}^{2}{y}^{2}+32\,{y}^{4}+9088\,{x}^{3}+3456\,x{y}^{2}
+13504\,{x}^{2}-2496\,{y}^{2}\nonumber\\
                &\quad\quad-10368\,x-13536 \big) {z}^{4}\nonumber\\
                &+ \big( 
128\,{x}^{5}+128\,{x}^{3}{y}^{2}+2496\,{x}^{4}+1792\,{x}^{2}{y}^{2}+64
\,{y}^{4}+6144\,{x}^{3}-384\,x{y}^{2}\nonumber\\
                &\quad\quad-5760\,{x}^{2}-16128\,{y}^{2}-
10368\,x-25920 \big) {z}^{5}\nonumber\\
                &+ \big( 8\,{x}^{6}+8\,{x}^{4}{y}^{2}-8
\,{x}^{2}{y}^{4}-8\,{y}^{6}+336\,{x}^{5}+352\,{x}^{3}{y}^{2}+16\,x{y}^
{4}+1400\,{x}^{4}\nonumber\\
                &\quad\quad+944\,{x}^{2}{y}^{2}+120\,{y}^{4}-1440\,{x}^{3}-10656\,x{y}^{2}-1800\,{x}^{2}-9720\,{y}^{2}\nonumber\\
                &\quad\quad-14256\,x-20088 \big) {z}^{6}\nonumber\\
                &+ \big( 16\,{x}^{6}+16\,{x}^{4}{y}^{2}-16\,{x}^{2}{y}^{4}-16\,{y}^{6}+
128\,{x}^{5}+256\,{x}^{3}{y}^{2}+128\,x{y}^{4}-144\,{x}^{4}\nonumber\\
                &-1440\,{x}^
{2}{y}^{2}-144\,{y}^{4}-4608\,x{y}^{2}-1296\,{x}^{2}+1296\,{y}^{2}-
10368\,x+11664 \big) {z}^{7}\\
Q_2(x,y,z)&=\big( {x}^{4}{z}^{4}+2\,{x}^{2}{y}^{2}{z}^{4}+{y}^{4}{z}^{4}+16\,{x}
^{3}{z}^{4}+16\,x{y}^{2}{z}^{4}+16\,{x}^{3}{z}^{3}+46\,{x}^{2}{z}^{4}\nonumber\\
                &+16\,x{y}^{2}{z}^{3}+82\,{y}^{2}{z}^{4}+168\,{x}^{2}{z}^{3}-144\,x{z}^{
4}+88\,{y}^{2}{z}^{3}+88\,{x}^{2}{z}^{2}+176\,x{z}^{3}\nonumber\\
                &+40\,{z}^{2}{y}^
{2}+81\,{z}^{4}+512\,x{z}^{2}-360\,{z}^{3}+192\,xz+184\,{z}^{2}+480\,z
+144 \big) ^{2}
\end{align}
It is straightforward to check that for $x\in[0,1/2]$, $y\in[4,10]$ and $z\in[0,1/10]$, both $P_2(x,y,z)$ and $Q_2(x,y,z)$ are positive, so
\begin{equation}
\frac{\partial}{\partial x}M_1(x,y,z) \geq 0
\end{equation}
so
\begin{equation}
M_1(x,y,z) \leq M_1(1/2,y,z)
\end{equation}
Use similar method we also have:
\begin{equation}
\frac{\partial}{\partial y}M_1(x,y,z) \geq 0
\end{equation}
so
\begin{equation}
M_1(1/2,y,z) \leq M_1(1/2,10,z)
\end{equation}
To summarize, for any fixed $z\in[0,1/10]$ the maximum of $M_1$ is attained at $x=1/2$, $y=10$. In other words, given fixed $n \geq 10$, the maximum of $|B_n(\lambda)/A_n(\lambda)|$ on $S_3$ is attained at $\lambda=1/2+10 i$. Let $Z^{(1)}_n$ denote the maximum of $|B_n(\lambda)/A_n(\lambda)|$ on $S_3$. Then
\begin{align}
\label{Z1}
|B_n(\lambda)/A_n(\lambda)| &\leq Z^{(1)}_n \nonumber\\
&=\sqrt{M_1(1/2,10,1/n)}\nonumber\\
                &= \left\vert\frac{B_n(1/2+10\,i)}{A_n(1/2+10\,i)}\right\vert\ \nonumber\\
                &=\left({\frac { \left( 16\,{n}^{2}-40\,n+425 \right)  \left( 16\,{n}^{2}+56\,
n+449 \right) }{2304\,{n}^{4}+9216\,{n}^{3}+71392\,{n}^{2}+149952\,n+
305161}} \right)^{1/2}
\end{align}
If we take derivative of $M_1(1/2,10,z)$ with respect to z, then
\begin{equation}
\frac{\partial}{\partial z}M_1(1/2,10,z) = \frac{P_3(z)}{Q_3(z)}
\end{equation}
where
\begin{align}
P_3(z)&=26832450160\,{z}^{6}+20079135232\,{z}^{5}+3697458944\,{z}^{4}+
1477025792\,{z}^{3}\nonumber\\
          &+15159296\,{z}^{2}+17563648\,z-1769472\\
 Q_3(z)&=\left( 305161\,{z}^{4}+149952\,{z}^{3}+71392\,{z}^{2}+9216\,z+2304
 \right) ^{2}
\end{align}
It is obvious that $Q_3$ is positive and $P_3$ is increasing on $[0,1/10]$. We notice that
\begin{align}
P_3(1/15)>0\\
P_3(1/16)<0
\end{align}
Hence for $z\in [0,1/16]$, $M_1$ decreases as $z$ increases. In other words, for $n\geq 16$, $Z^{(1)}_n$ increases as $n$ increases. Since
\begin{equation}
\frac{B_n}{A_n} = {\frac {{\lambda}^{2}+4\,\lambda\,n+4\,{n}^{2}-9}{-12\,{n}^{2}-
 \left( 20+8\,\lambda \right) n-{\lambda}^{2}-8\,\lambda+9}}
\end{equation} 
we can see that for each $\lambda$, as $n\to\infty$, $|B_n/A_n| \to 1/3$, so $Z^{(1)}_n \uparrow 1/3$ as $n \to \infty$.

\subsubsection{Lower bound of $|A_n|$}\label{S_32}
The study of $A_n$ is similar.
\begin{equation}
A_n(\lambda) = {\frac {-12\,{n}^{2}- \left( 20+8\,\lambda \right) n-{\lambda}^{2}-8\,
\lambda+9}{8\,{n}^{2}+28\,n+20}}
\end{equation}
Let $\lambda= x+iy$, then 
\begin{equation}
F_2(x,y,n):= {A_n(x+iy)}
\end{equation}
is a rational function in $x$, $y$ and $n$. Let $z=1/n$, then
\begin{align}
M_2(x,y,z):=|F_2(x,y,1/z)|^2 = \frac{P_4(x,y,z)}{Q_4(x,y,z)}
\end{align}
is a real valued rational function, where
\begin{align}
P_4(x,y,z)&={x}^{4}{z}^{4}+2\,{x}^{2}{y}^{2}{z}^{4}+{y}^{4}{z}^{4}+16\,{x}^{3}{z}^
{4}+16\,x{y}^{2}{z}^{4}+16\,{x}^{3}{z}^{3}+46\,{x}^{2}{z}^{4}\nonumber\\
                &+16\,x{y}
^{2}{z}^{3}+82\,{y}^{2}{z}^{4}+168\,{x}^{2}{z}^{3}-144\,x{z}^{4}+88\,{
y}^{2}{z}^{3}+88\,{x}^{2}{z}^{2}\nonumber\\
                &+176\,x{z}^{3}+40\,{z}^{2}{y}^{2}+81\,
{z}^{4}+512\,x{z}^{2}-360\,{z}^{3}+192\,xz+184\,{z}^{2}\nonumber\\
                &+480\,z+144 \\
Q_4(x,y,z)&=16\, \left( z+1 \right) ^{2} \left( 5\,z+2 \right) ^{2}
\end{align}
Take partial derivative of $M_2(x,y,z)$ with respect to x:
\begin{equation}
\frac{\partial}{\partial x}M_2(x,y,z) = z \frac{P_5(x,y,z)}{Q_5(x,y,z)}
\end{equation}
where
\begin{align}
P_5(x,y,z)&=\left( xz+4\,z+4 \right)  \left( {x}^{2}{z}^{2}+{z}^{2}{y}^{2}+8\,x{z
}^{2}+8\,xz-9\,{z}^{2}+20\,z+12 \right) \\
Q_5(x,y,z)&=4\, \left( z+1 \right) ^{2} \left( 5\,z+2 \right) ^{2}\\
\end{align}
By inspection $P_6$ and $Q_6$ are both positive, so
\begin{equation}
\frac{\partial}{\partial x}M_2(x,y,z) \geq 0
\end{equation}
\begin{equation}
M_2(x,y,z) \geq M_2(0,y,z)
\end{equation}
Now consider the partial derivative of $M_2(0,y,z)$with respect to y:
\begin{equation}
\frac{\partial}{\partial y}M_2(0,y,z) = {\frac {4\,{y}^{3}{z}^{4}+164\,y{z}^{4}+176\,y{z}^{3}+80\,{z}^{2}y}{ 16\,\left( z+1 \right) ^{2} \left( 5\,z+2 \right) ^{2}}}
\end{equation}
It is obviously nonnegative for $y\geq0$ and $z\geq0$, so
\begin{equation}
M_2(0,y,z) \geq M_2(0,4,z)
\end{equation}
So given any $z\in[0,1/10]$ the minimum of $M_2(x,y,z)$ is attained at $x=0$, $y=4$. Equivalently, for each $n\geq 10$, the minimum of $|A_n|$ on $S_3$ is attained at $\lambda = 4\,i$. Let $Z^{(2)}_n$ denote the minimum of $|A_n(\lambda)|$ on $S_3$. Then
\begin{align}\label{Z2}
|A_n(\lambda)| &\geq Z^{(2)}_n \nonumber\\
&=\sqrt{M_2(0,4,1/n)}\nonumber\\
                &= \left\vert{A_n(4\,i)}\right\vert\ \nonumber\\
                &=\left({\frac {144\,{n}^{4}+480\,{n}^{3}+824\,{n}^{2}+1048\,n+1649}{16\,
 \left( 1+n \right) ^{2} \left( 5+2\,n \right) ^{2}}}
 \right)^{1/2}
\end{align}
Notice that
\begin{align}
\frac{\partial}{\partial z}M_2(0,4,z)={\frac {8923\,{z}^{4}+6144\,{z}^{3}-456\,{z}^{2}-1472\,z-528}{8\,\left( z+1 \right) ^{3} \left( 5\,z+2 \right) ^{3}}}
\end{align}
By inspection the derivative is negative for $z\in[0,1/10]$, so for $n\geq 10$, $Z^{(2)}_n$ is monotonically increasing. Since for each $n$, $\lim_{n\to\infty}|A_n| = 3/2$, $Z^{(2)}_n \uparrow 3/2$ as $n \to \infty$. 

\subsubsection{Upper bound of $|B_n/A_n^2|$.}\label{S_33}
\begin{equation}
\frac{B_n(\lambda)}{A_n^2(\lambda)} = {\frac { \left( {\lambda}^{2}+4\,\lambda\,n+4\,{n}^{2}-9 \right)  \left( 8\,{n}^{2}+28\,n+20 \right) }{ \left( -12\,{n}^{2}- \left( 20+
8\,\lambda \right) n-{\lambda}^{2}-8\,\lambda+9 \right) ^{2}}}
\end{equation}
Let $\lambda= x+iy$, then 
\begin{equation}
F_3(x,y,n):= \frac{B_n(x+i y)}{A_n^2(x+iy)}
\end{equation}
is a rational function in $x$, $y$ and $n$. Let $z=1/n$, then
\begin{align}
M_3(x,y,z):=|F_3(x,y,1/z)|^2
\end{align}
is a real valued rational function. Consider the partial derivative of $M_3(x,y,z)$ with respect to $x$.
\begin{equation}
\frac{\partial}{\partial x}M_3(x,y,z) = 64\,z\,(z+1)^2\,(5\,z+2)^2 \, \frac{P_6(x,y,z)}{Q_6(x,y,z)}
\end{equation}
where
\begin{align}
P_6(x,y,z)&=\sum_{i=0}^7 f_i(x,y) z^i\\
Q_6(x,y,z)&= \Big( {x}^{4}{z}^{4}+2\,{x}^{2}{y}^{2}{z}^{4}+{y}^{4}{z}^{4}+16\,{x}
^{3}{z}^{4}+16\,x{y}^{2}{z}^{4}+16\,{x}^{3}{z}^{3}+46\,{x}^{2}{z}^{4}\nonumber\\
                &+16\,x{y}^{2}{z}^{3}+82\,{y}^{2}{z}^{4}+168\,{x}^{2}{z}^{3}-144\,x{z}^{
4}+88\,{y}^{2}{z}^{3}+88\,{x}^{2}{z}^{2}+176\,x{z}^{3}\nonumber\\
                &+40\,{z}^{2}{y}^
{2}+81\,{z}^{4}+512\,x{z}^{2}-360\,{z}^{3}+192\,xz+184\,{z}^{2}+480\,z
+144 \Big) ^{3}
\end{align}
where
\begin{align}
f_0(x,y)=&-384\\
f_1(x,y)=&-1216\,x-256\\
f_2(x,y)=&-1632\,{x}^{2}-288\,{y}^{2}-1024\,x+4384\\
f_3(x,y)=&-1200\,{x}^{3}-624\,x{y}^{2}-1536\,{x}^{2}-512\,{y}^{2}+8560\,x+8064\\
f_4(x,y)=&-520\,{x}^{4}-528\,{x}^{2}{y}^{2}-72\,{y}^{4}-1152\,{x}^{3}-960\,x{y}^{2}+6480\,{x}^{2}\nonumber\\
&\quad-1552\,{y}^{2}+13824\,x-4104\\
f_5(x,y)=&-132\,{x}^{5}-216\,{x}^{3}{y}^{2}-84\,x{y}^{4}-464\,{x}^{4}-672\,{x}^{
2}{y}^{2}-144\,{y}^{4}\nonumber\\
&\quad+2360\,{x}^{3}-1064\,x{y}^{2}+8640\,{x}^{2}-5760
\,{y}^{2}-1332\,x-19440\\
f_6(x,y)=&-18\,{x}^{6}-42\,{x}^{4}{y}^{2}-30\,{x}^{2}{y}^{4}-6\,{y}^{6}-96\,{x}^
{5}-208\,{x}^{3}{y}^{2}\nonumber\\
&\quad-112\,x{y}^{4}+410\,{x}^{4}-108\,{x}^{2}{y}^{2}
-86\,{y}^{4}+2304\,{x}^{3}-3600\,x{y}^{2}\nonumber\\
&\quad+1026\,{x}^{2}-3546\,{y}^{2}-
12960\,x-8586\\
f_7(x,y)=&-{x}^{7}-3\,{x}^{5}{y}^{2}-3\,{x}^{3}{y}^{4}-x{y}^{6}-8\,{x}^{6}-24\,{
x}^{4}{y}^{2}-24\,{x}^{2}{y}^{4}\-8\,{y}^{6}\nonumber\\
&\quad+27\,{x}^{5}+18\,{x}^{3}{y}
^{2}-9\,x{y}^{4}+216\,{x}^{4}-432\,{x}^{2}{y}^{2}-72\,{y}^{4}\nonumber\\
&\quad+333\,{x}
^{3}-1647\,x{y}^{2}-1944\,{x}^{2}+648\,{y}^{2}-4455\,x+5832
\end{align}
It is obvious that $Q_6(x,y,z)>0$, $f_1(x,y)\leq0$ and $f_2(x,y)\leq0$ for $x\in[0,1/2]$, $y\in[4,10]$ and $z\in=[0,1/10]$. By taking absolute value of every term and use the condition that $|x| \leq 1/2 $ and $|y| \leq 10$ we obtain crude estimates for upper bounds of $|f_i(x,y)|$. Thus
\begin{align}
P_6(x,y,z)&\leq \sum_{i=0}^7 f_i(x,y) z^i \nonumber\\
                 &\leq -384 + \sum_{i=3}^7 |f_i(x,y)| (1/10)^i\nonumber\\
                 &<-384+96+95+26+9<0
\end{align}
We get
\begin{equation}
\frac{\partial}{\partial x}M_3(x,y,z) \leq 0
\end{equation}
\begin{equation}\label{monoM3}
M_3(x,y,z) \leq M_3(0,y,z)
\end{equation}

From \S\ref{S_31} and \S\ref{S_32} we see that given $n$ both $|B_n(i y)/A_n(iy)|$ and $|A_n(iy)|$ are monotonically increasing for $y\in[4,10]$, so for $\lambda \in [j,j+1]$, $j\in\{4,5,..,9\}$ we have
\begin{equation}
\left\vert\frac{B_n(\lambda)}{A_n^2(\lambda)}\right\vert^2 =
\left\vert \frac{\frac{B_n(\lambda)}{A_n(\lambda)}}{A_n(\lambda)}\right\vert
\leq \frac{M_1(0,j+1,1/n)}{M_2(0,j,1/n)}
\end{equation}

For $j\in\{4,5,..,9\}$ let 
\begin{equation}
M_{4,j}(z)= \frac{M_1(0,j+1,z)}{M_2(0,j,z)}
\end{equation}
\begin{equation}\label{Z3}
Z^{(3)}(n)= \max_{j=4,5..,9}\left\{\left(M_{4,j}(z)\right)^{1/2}\right\}
\end{equation}
then for all $\lambda \in [4i,10i]$
\begin{equation}
\left\vert\frac{B_n(\lambda)}{A_n^2(\lambda)}\right\vert^2 \leq Z^{(3)}(n)^2
\end{equation}

In view of \eqref{monoM3}, for each $n$ and $\lambda \in S_3$:
\begin{equation}\label{Z32}
\left\vert\frac{B_n(\lambda)}{A_n^2(\lambda)}\right\vert \leq Z^{(3)}(n)
\end{equation}

The expressions of $M_{4,j}$ are:
\begin{align}
&\quad\quad M_{4,4}(z) = \nonumber\\
&{\frac {16\, \left( 34\,{z}^{2}+12\,z+4 \right)  \left( 34\,{z}^{2}-12
\,z+4 \right)  \left( z+1 \right) ^{2} \left( 5\,z+2 \right) ^{2}}{
 \left( 2756\,{z}^{4}+1840\,{z}^{3}+1184\,{z}^{2}+480\,z+144 \right) 
 \left( 1649\,{z}^{4}+1048\,{z}^{3}+824\,{z}^{2}+480\,z+144 \right) }}\nonumber\\
 &\quad\quad M_{4,5}(z) = \nonumber\\
&{\frac { 16\,\left( 45\,{z}^{2}+12\,z+4 \right)  \left( 45\,{z}^{2}-12
\,z+4 \right)  \left( z+1 \right) ^{2} \left( 5\,z+2 \right) ^{2}}{
 \left( 4329\,{z}^{4}+2808\,{z}^{3}+1624\,{z}^{2}+480\,z+144 \right) 
 \left( 2756\,{z}^{4}+1840\,{z}^{3}+1184\,{z}^{2}+480\,z+144 \right) }
}\nonumber\\
&\quad\quad M_{4,6}(z) = \nonumber\\
&{\frac { 16\,\left( 73\,{z}^{2}+12\,z+4 \right)  \left( 73\,{z}^{2}-12
\,z+4 \right)  \left( z+1 \right) ^{2} \left( 5\,z+2 \right) ^{2}}{
 \left( 9425\,{z}^{4}+5272\,{z}^{3}+2744\,{z}^{2}+480\,z+144 \right) 
 \left( 6500\,{z}^{4}+3952\,{z}^{3}+2144\,{z}^{2}+480\,z+144 \right) }}\nonumber\\
 &\quad\quad M_{4,7}(z) = \nonumber\\
 &{\frac {16\, \left( 73\,{z}^{2}+12\,z+4 \right)  \left( 73\,{z}^{2}-12
\,z+4 \right)  \left( z+1 \right) ^{2} \left( 5\,z+2 \right) ^{2}}{
 \left( 9425\,{z}^{4}+5272\,{z}^{3}+2744\,{z}^{2}+480\,z+144 \right) 
 \left( 6500\,{z}^{4}+3952\,{z}^{3}+2144\,{z}^{2}+480\,z+144 \right) }}\nonumber\\
 &\quad\quad M_{4,8}(z) = \nonumber\\
 &{\frac { 16\,\left( 90\,{z}^{2}+12\,z+4 \right)  \left( 90\,{z}^{2}-12
\,z+4 \right)  \left( z+1 \right) ^{2} \left( 5\,z+2 \right) ^{2}}{
 \left( 13284\,{z}^{4}+6768\,{z}^{3}+3424\,{z}^{2}+480\,z+144 \right) 
 \left( 9425\,{z}^{4}+5272\,{z}^{3}+2744\,{z}^{2}+480\,z+144 \right) }}\nonumber\\
 &\quad\quad M_{4,9}(z) = \nonumber\\
 &{\frac { 16\,\left( 109\,{z}^{2}+12\,z+4 \right)  \left( 109\,{z}^{2}-
12\,z+4 \right)  \left( z+1 \right) ^{2} \left( 5\,z+2 \right) ^{2}}{
 \left( 18281\,{z}^{4}+8440\,{z}^{3}+4184\,{z}^{2}+480\,z+144 \right) 
 \left( 13284\,{z}^{4}+6768\,{z}^{3}+3424\,{z}^{2}+480\,z+144 \right) }}\nonumber\\
\end{align}

It is straightforward to check that each $M_{4,j}(z)$ has a positive derivative with respect to $z$ on $[0,1/10]$. Hence as $n$ increases, $Z^{(3)}_n$ decreases monotonically.


\subsection{Proof of \eqref{eqS_34} and \eqref{eqS_35}}
\subsubsection{Proof of \eqref{eqS_34}}
\label{S_34}
\begin{align}
\left\vert\left[(\mathcal{N}(\mathbf{0}))-(\mathcal{N}^{2}(\mathbf{0})) \right]_n\right\vert
&= \left\vert-\frac{B_n}{A_n} - \left( -\frac{B_n}{A_n-\frac{B_{n+1}}{A_{n+1}}}\right)\right\vert\nonumber\\
&= \left\vert\frac{B_n}{A_n}\right\vert \left\vert\frac{\frac{B_{n+1}}{A_{n+1}}}{A_n - \frac{B_{n+1}}{A_{n+1}}}\right\vert\nonumber\\
&\leq \frac{Z^{(1)}_n\,Z^{(1)}_{n+1}}{Z^{(2)}_n-Z^{(1)}_{n+1}}
\end{align}

For $10\leq n\leq 20$ we can check directly from \eqref{Z1} and \eqref{Z2} that the value is less than $0.106$. For $n>20$ we have 
\begin{equation}
\frac{Z^{(1)}_n\,Z^{(1)}_{n+1}}{Z^{(2)}_n-Z^{(1)}_{n+1}} \leq \frac{(1/3)^2}{Z^{(2)}_{21} - 1/3} <0.106
\end{equation}

\subsubsection{Proof of \eqref{eqS_35}}
\label{S_35}
Note that we can write
\begin{align}
\left[\mathcal{N}(\mathbf{0})\right]_{n_0} &= \left[\mathcal{N}(\mathbf{0})\right]_{10} \nonumber\\
&=1-{\frac {4}{545}}\,{\frac { \left( -139+6\,\sqrt {545} \right) \sqrt 
{545}}{\lambda+44-\sqrt {545}}}-{\frac {4}{545}}\,{\frac { \left( 139+
6\,\sqrt {545} \right) \sqrt {545}}{\lambda+44+\sqrt {545}}}
\end{align}
so we have unique numbers $\tilde{e}_i$, ($i= 18,19$), and $s_{18}=-44+\sqrt{545}$, $s_{19}=-44-\sqrt{545}$, such that
\begin{align}
\left[\mathcal{N}(\mathbf{0})\right]_{10} =1+ \sum_{i=18}^{19}\frac{\tilde{e}_i}{\lambda-s_i} 
\end{align}

To simplify calculations, we approximate $\left[\mathcal{N}(\mathbf{0})\right]_{10}$ by $Q(\lambda)$, where
\begin{equation}
Q(\lambda)=1+\frac{\tilde{d}_{18}}{\lambda-\tilde{s}_{18}}+\frac{\tilde{d}_{19}}{\lambda-\tilde{s}_{19}} \nonumber\\
\end{equation}
where
\begin{align}
&\tilde{s}_{18} = -1735/84,\quad \tilde{s}_{19} = -1953/29,\nonumber\\
&\tilde{d}_{18} = -38/207, \quad\, \tilde{d}_{19} = -2343/49
\end{align}
By straightforward calculation:
\begin{align}
\label{Nrlink1}
\left\vert \left[\mathcal{N}(\mathbf{0})\right]_{10} - Q(\lambda)\right\vert
    &\leq \sum_{i=18}^{19}\left\vert \frac{\tilde{e}_i}{\lambda-s_i} -\frac{\tilde{d}_{i}}{\lambda-\tilde{s}_{i}}\right\vert \nonumber\\
    &=\sum_{i=18}^{19} \frac{|\lambda|\,|\tilde{e}_i-\tilde{d}_i |+|\tilde{e}_i \tilde{s_i} - \tilde{d_i} s_i |}{|\lambda-s_i| |\lambda- \tilde{s_i}|} \nonumber\\
    &<\frac{1}{312500}
\end{align}

$c_{10}$ is a polynomial of degree $18$ and it has $18$ distinct real roots $\{s_j\}_{0\leq i\leq {17}}$ satisfying
\begin{equation}
1 = s_0 >-0.54>s_1>s_2>...>s_{17}
\end{equation}
$c_{11}$ is a polynomial of degree $20$ and it has $20$ distinct real roots $\{t_j\}_{0 \leq i \leq 19}$ satisfying
\begin{equation}
1 = t_0 >-0.54>t_1>t_2>...>t_{19}
\end{equation}
The only common root of $c_{10}$ and $c_{11}$ is $\lambda = 1$. We can write:
\begin{equation}
c_{10} = l_1 \prod_{i=0}^{17} (\lambda-s_i)
\end{equation}
\begin{equation}
c_{11} = l_2 \prod_{i=0}^{19} (\lambda-t_i)
\end{equation}
where $l_1$ is the leading coefficient of $c_{10}$ and $l_2$ is the leading coefficient of $c_{11}$, so 
\begin{equation}
r_{10} = \frac{c_{11}}{c_{10}} = P_1(\lambda) + \sum_{i=1}^{17} \frac{\tilde{a}_i}{\lambda-s_i}
\end{equation}
where $\tilde{a}_i$'s are constants and $P_1(\lambda)$ is a quadratic polynomial.\\

Now we approximate $c_{10}$ by $\tilde{c}_{10}$
\begin{equation}
\tilde{c}_{10} = l_1 \prod_{i=0}^{17} (\lambda-\tilde{s}_i)
\end{equation}
where for each $1 \leq i \leq 17$, $\tilde{s}_i$ is a rational number within $e^{-9}$ of $s_i$, and $\tilde{s}_0 = s_0 = 1$.
Now we crudely estimate $c_{11}$, $c_{10} - \tilde{c}_{10}$ and $\tilde{c}_{10}$.
Rewrite
\begin{equation}
c_{11} = \sum_{j=0}^{20}\tilde{b}_j (\lambda-7\,i-1/4)^j
\end{equation}
then
\begin{equation}
\label{Nrlink2}
|c_{11}| \leq \sum_{j=0}^{20} |\tilde{b}_j | \left(\frac{\sqrt{145}}{4}\right)^j < 37209/50		\quad(\lambda \in S_3)
\end{equation}
Rewrite
\begin{equation}
c_{10} - \tilde{c}_{10} = \sum_{j=0}^{18}\tilde{c}_j (\lambda-7\,i-1/4)^j
\end{equation}
then
\begin{equation}
\label{Nrlink3}
|c_{10} - \tilde{c}_{10}| \leq \sum_{j=0}^{18} |\tilde{c}_j | \left(\frac{\sqrt{145}}{4}\right)^j < 3/250000		\quad(\lambda \in S_3)
\end{equation}
By definition
\begin{equation}
\tilde{c}_{10} = l_1 \prod_{i=0}^{17}(\lambda-\tilde{s}_i) = l_1 (\lambda - 1) \prod_{i=1}^{17}(\lambda-\tilde{s}_i)
\end{equation}

Since for each $1 \leq i \leq 17$, $\tilde{s}_i <0$, we have
\begin{equation}
\label{Nrlink4}
|\tilde{c}_{10}|>|l_1| |4\, i - 1/2| \prod_{i=1}^{17}|4\,i-\tilde{s}_i| > 786187/1000000 
\end{equation}

Hence from \eqref{Nrlink2}, \eqref{Nrlink3} and \eqref{Nrlink4} we have
\begin{align}\label{Nrlink5}
\left\vert \frac{c_{11}}{c_{10}} - \frac{c_{11}}{\tilde{c}_{10}} \right\vert 
    &= \frac{|c_{11}||\tilde{c}_{10}-c_{10}|}{|c_{10} \tilde{c}_{10}|}\nonumber\\
    &=\frac{\frac{37209}{50}\frac{3}{250000}}{\left(\frac{786187}{1000000}-\frac{3}{250000}\right)\left(\frac{786187}{1000000}\right)}\nonumber\\
    &<0.0145
\end{align}

Now we have explicit numbers $\tilde{e}_i$'s and a quadratic polynomial $\tilde{P}_1(\lambda)$ such that
\begin{equation}
\frac{c_{11}}{\tilde{c}_{10}} = \tilde{P}_1(\lambda) + \sum_{i=1}^{17} \frac{-\tilde{e}_i}{\lambda - \tilde{s}_i}
\end{equation}

To simplify the calculation, replace each $\tilde{e}_i$ by a rational number $\tilde{d}_i$ with $e^{-7}$ of it and each coefficient in $\tilde{P}_1(\lambda)$ by a rational number within $e^{-7}$ of it we obtain:
\begin{equation}
\tilde{r}_{10} := \tilde{\tilde{P}}_1(\lambda) + \sum_{i=1}^{17}  \frac{-\tilde{d}_i}{\lambda - \tilde{s}_i}
\end{equation}
By straightforward calculation
\begin{equation}\label{Nrlink6}
\left\vert \tilde{r}_{10} - \frac{c_{11}}{\tilde{c}_{10}} \right\vert \leq \left\vert\tilde{\tilde{P}}_1(\lambda)- \tilde{P}_1(\lambda)\right\vert + \sum_{i=1}^{17}  \frac{|\tilde{d}_i-\tilde{e}_i|}{|4\,i - \tilde{s}_i|} <3\cdot 10^{-7}
\end{equation}
Finally, consider 
\begin{align}
F(\lambda):=Q(\lambda)-\tilde{r}_{10} = 1-\tilde{\tilde{P}}_1(\lambda) + \sum_{i=1}^{19}\frac{ \tilde{d}_i}{\lambda-\tilde{s_i}}
\end{align}
(see \eqref{Flambda} for explicit expression of $F(\lambda)$). Use the method in \S\ref{rationalme} with partition $\mathcal{P} = \{4\,i,4\,i+1/2,7\,i+1/2,10\,i+1/2, 7\,i\}$ and $m_0 =2$, we have
\begin{align}
\label{Nrlink7}
|F(\lambda)| > 0.6098
\end{align}
In view of \eqref{Nrlink1}, \eqref{Nrlink5} \eqref{Nrlink6} and \eqref{Nrlink7} we have:
\begin{align}
&\left\vert \left[\mathcal{N}(\mathbf{0})\right]_{10} - \frac{c_{11}}{c_{10}}\right\vert  \nonumber\\
    &\geq |Q(\lambda)-\tilde{r}_{10}| - \left\vert \left[\mathcal{N}(\mathbf{0})\right]_{10} - Q(\lambda)\right\vert - \left\vert \frac{c_{11}}{c_{10}} - \frac{c_{11}}{\tilde{c}_{10}} \right\vert -\left\vert \tilde{r}_{10} - \frac{c_{11}}{\tilde{c}_{10}} \right\vert \nonumber\\
    &>0.595
\end{align}
Thus we have \eqref{eqS_35}.

The expression of $F(\lambda)$:
\begin{align}\label{Flambda}
F(\lambda)
     &={\frac {1}{920}}\,{\lambda}^{2}+\frac{2}{23}\,\lambda+{\frac {387}{920}}-{\frac {38}{207}}\, \left( \lambda+{\frac {1735}{84
}} \right) ^{-1}-{\frac {2343}{49}}\, \left( \lambda+{\frac {1953}{29}
} \right) ^{-1}\nonumber\\
     &-{\frac {
1704}{101}}\, \left( \lambda+{\frac {108406}{1283}} \right) ^{-1}-{
\frac {4297}{537}}\, \left( \lambda+{\frac {73708}{1189}} \right) ^{-1
}-{\frac {2765}{768}}\, \left( \lambda+{\frac {59605}{1299}} \right) ^
{-1}\nonumber\\
     &-{\frac {1588}{1195}}\, \left( \lambda+{\frac {53402}{1589}}
 \right) ^{-1}-{\frac {763}{3245}}\, \left( \lambda+{\frac {31328}{
1279}} \right) ^{-1}+{\frac {53}{10426}}\, \left( \lambda+{\frac {
19885}{1029}} \right) ^{-1}\nonumber\\
     &+{\frac {23}{2173}}\, \left( \lambda+{
\frac {52237}{3423}} \right) ^{-1}+{\frac {367}{7339}}\, \left( 
\lambda+{\frac {21932}{1491}} \right) ^{-1}+{\frac {160}{4019}}\,
 \left( \lambda+{\frac {38219}{3029}} \right) ^{-1}\nonumber\\
     &+{\frac {71}{9409}}
\, \left( \lambda+{\frac {1316}{117}} \right) ^{-1}+{\frac {358}{3591}
}\, \left( \lambda+{\frac {28039}{2844}} \right) ^{-1}+{\frac {717}{
10073}}\, \left( \lambda+{\frac {11626}{1435}} \right) ^{-1}\nonumber\\
     &+{\frac {
274}{14543}}\, \left( \lambda+{\frac {36794}{5013}} \right) ^{-1}+{
\frac {280}{3639}}\, \left( \lambda+{\frac {121869}{21650}} \right) ^{
-1}+{\frac {179}{3360}}\, \left( \lambda+{\frac {15715}{4138}}
 \right) ^{-1}\nonumber\\
     &+{\frac {159}{8224}}\, \left( \lambda+{\frac {19660}{
9407}} \right) ^{-1}+{\frac {43}{17095}}\, \left( \lambda+{\frac {2595
}{4718}} \right) ^{-1}\nonumber\\
\end{align}

\subsection{The method to prove \eqref{L1}, \eqref{es1} -- \eqref{es4}, and \eqref{es5} -- \eqref{es6}}\label{esmet}
The functions we estimate on $\partial R_2$ are of two types:\\

Group (I)
\begin{align}
& \delta^{(1)}_{1} = r_1 e^{-W_1} -1\\
& \epsilon^{(1)}_{n}=-1-A_{n-1} e^{-W_{n}}-B_{n-1} e^{-W_{n}-W_{n-1}} \quad (2\leq n\leq 50)\\
& \epsilon^{(3)}=-1-A_{50}\,e^{-W_{50}}-B_{50}\,e^{-2\,W_{50}}
\end{align}

Group (II)
\begin{align}
C^{(1)}_n&=B_n e^{-W_{n+1}-W_n}\\
C^{(3)} &=B_{50}\, e^{-2\,W_{50}}
\end{align}

We obtain estimates for each function on each one of the following subintervals of $l_i$, ($1\leq i \leq 6$)
\begin{align}
&\begin{array}{llll}
&{l}_{1,1}=[10\,i,18\,i] & {l}_{1,2}=[18\,i,30\,i] & {l}_{1,3}=[30\,i,48\,i] \\
& {l}_{1,4}=[48\,i,70\,i] &  &  \\
\end{array}\\
&\begin{array}{llll}
&{l}_{3,j}=l_{1,j}+1/2 \quad\quad(j=1,2,3,4)& &\\
\end{array}\\
&\begin{array}{llll}
&{l}_{2,1}=(70\,i,110\,i] & {l}_{2,2}=[110\,i,175\,i] & {l}_{2,3}=[175\,i,275\,i] \\
& {l}_{2,4}=[275\,i,380\,i] &  &  \\
\end{array}\\
&\begin{array}{llll}
&{l}_{4,j}=l_{2,j}+1/2 \quad\quad(j=1,2,3,4)& &\\
\end{array}\\
&\begin{array}{llll}
&{l}_{5,1}=l_{5}=[10\,i,10\,i+1/2]& &\\
\end{array}\\
&\begin{array}{llll}
&{l}_{6,1}=l_{6}=[380\,i,380\,i+1/2]& &\\
\end{array}
\end{align}

For each function in Group (II), the estimate follows easily from intermediate results obtained in the process of estimating the corresponding function in Group (I). See \eqref{esgp2}. We shall describe the method we use to estimate a function $\tilde{G}(\lambda)$ of Group (I) on an interval $[a_0,b_0]$ in detail.\\

We see that $\tilde{G}(\lambda)$ is of the form
\begin{equation}
\tilde{G}(\lambda) = \tilde{q}_1 e^{\tilde{f}_1}+\tilde{q}_2 e^{\tilde{f}_2}-1
\end{equation}
where $\tilde{q}_1$ and $\tilde{q}_2$ are polynomials of degree 2 in $\lambda$ and $\tilde{f}_1$ and $\tilde{f}_2$ are polynomials of degree {at most} 5 in $\Im(\lambda)$ and of degree 1 in $\Re(\lambda)$. First let $\phi$ be a polynomial of degree 1 which maps $[-1,1]$ bijectively to $[a,b]$ and let
\begin{align}
G:=\tilde{G}\circ \phi \quad\quad q_i:=\tilde{q}_i \circ \phi \quad\quad f_i:=\tilde{f}_i \quad\quad (i=1,2)
\end{align}
The problem thus transforms into finding an upper bound for the function $|G|$ of a single real variable $x$
\begin{equation}
\left\vert {G}\right\vert = \left\vert {q}_1 e^{{f}_1}+{q}_2 e^{{f}_2}-1 \right\vert \quad\quad x\in [-1,1]
\end{equation}
where $q_{1,2}$ are a polynomials of degree 2 in $x$ and ${f}_{1,2}$ are polynomials of degree 5 or 1 in $x$.\\
The estimate of $|G|$ is obtained in 3 steps:\\
Step 1: For each $i$, express $f_i$ as a linear combination of Chebyshev polynomials $T_j(x)$, $0\le j \le 5$.
\begin{align}
f_i(x) &= \sum_{j=0}^5 g_{i,j} T_j(x)
\end{align}
Denote
\begin{align}\label{tilder}
h_i(x) = \sum_{j=1}^2 g_{i,j} T_j(x)\,\, , \quad\quad
\tilde{r}_{i}(x)  = \sum_{j=3}^5 g_{i,j} T_j(x)\,\, , \quad\quad
\overline{r}_i = \sum_{j=3}^5 \left\vert g_{i,j}\right\vert
\end{align}
then
\begin{equation}
f_i = g_{i,0} + h_i + \tilde{r}_{i}
\end{equation}

Step 2: For each $i$, approximate $e^{f_i}$ by a polynomial
\begin{align}
\begin{aligned}
& e^{f_i}  = e^{g_{i,0}+h_i+\tilde{r}_i} \\
&=e^{g_{i,0}}\left(1+h_i+\frac{h_i^2}{2}+\frac{h_i^3}{6}+\int_{0}^{h_i}\int_{0}^{s_1}\int_{0}^{s_2}\int_{0}^{s_3}e^{s_4} \mathrm{d}s_4 \, \mathrm{d}s_3 \, \mathrm{d}s_2 \, \mathrm{d}s_1 \right) \left(1+\int_0^{\tilde{r}_i} e^s\mathrm{d} s\right)\\
& =H_i + R_{i,1} +R_{i,2}
\end{aligned}
\end{align}
where 
\begin{align}
& H_i = e^{g_{i,0}} \left(1+h_i+\frac{h_i^2}{2}+\frac{h_i^3}{6}\right)\\
& R_{i,1} = e^{g_{i,0}} e^{h_i} \int_0^{\tilde{r}_i} e^s\mathrm{d} s\\
& R_{i,2} = e^{g_{i,0}} \left(\int_{0}^{h_i}\int_{0}^{s_1}\int_{0}^{s_2}\int_{0}^{s_3}e^{s_4} \mathrm{d}s_4 \, \mathrm{d}s_3 \, \mathrm{d}s_2 \, \mathrm{d}s_1\right)
\end{align}

Step 3: We have
\begin{align}
\begin{aligned}
\left\vert G(x) \right\vert & \le \left\vert H(x)\right\vert + \left\vert q_1 R_{1,1}\right\vert + \left\vert q_1 R_{1,2}\right\vert + \left\vert q_2 R_{2,1}\right\vert + \left\vert q_2 R_{2,2}\right\vert
\end{aligned}
\end{align}
where
\begin{equation}
H(x) = q_1(x) H_1(x)+q_2(x) H_2(x)-1
\end{equation}
is a polynomial of degree $8$.

To estimate $H(x)$ we simply express $H(x)$ as a linear combination of Chebyshev polynomials $T_l(x)$, $0 \le l \le 8$.
\begin{align}
\left\vert H(x)\right\vert  = \left\vert \sum_{l=0}^{8} \tilde{h}_l T_l(x) \right\vert \le \sum_{l=0}^{8}\left\vert \tilde{h}_l \right\vert 
\end{align}
since
\begin{equation}\label{Tle1}
\sup_{x\in[-1,1]}|T_l(x)| \le 1 \quad\quad \forall l \in \NN
\end{equation}

For $i\in\{1,2\}$, to estimate $q_i R_{i,1}$ we first notice that by \eqref{tilder} and \eqref{Tle1} we have
\begin{equation}
\left\vert \tilde{r}_i(x)\right\vert \leq \overline{r}_i
\end{equation}
We choose $E_{i,1}$ such that
\begin{equation}
E_{i,1}\geq \left\|q_i e^{h_i}\right\| =\sup_{x\in[-1,1]} \left\vert q_i(x) e^{h_i(x)} \right\vert
\end{equation}

$E_{i,1}$ is chosen rigorously in the following way. For each $i \in \{1,2\}$, denote $Q_i = |q_i|^2$ and consider the function
\begin{equation}
S(x) := \left\vert q_i(x) e^{h_i(x)} \right\vert^2 = Q_i(x) e^{2 \Re{h_i(x)}}
\end{equation}
We approximate $S(x)$ by approximating first its derivative:
\begin{equation}
S'(x) = e^{2 \Re{h_i(x)}}\left(2 \Re\left(h'_i(x)\right) Q_i(x)+Q'_i(x)\right)
\end{equation}
Let $\tilde{\epsilon}$ be small enough. We construct a polynomial $\tilde{P}(x)$ with real coefficients and known roots $r_i$ such that explicitly
\begin{equation}
\left\vert \tilde{P}(x) - \left(2 \Re\left(h'_i(x)\right) Q_i(x)+Q'_i(x)\right) \right\vert < \tilde{\epsilon} \quad\quad x\in [-1,1]
\end{equation}
By construction, the maximum value of the approximating function
\begin{equation}
\tilde{S}(x) := S(0) + \int_{0}^{x} e^{2 \Re{h_i(s)}} \tilde{P}(s) \mathrm{d}s
\end{equation}
on $[-1,1]$ is known exactly. Also
\begin{align}
\left\vert S(x) - \tilde{S}(x)   \right\vert &\leq \left\vert \int_0^x  e^{2 \Re{h_i(s)}}\left[\tilde{P}(s) - \left(2 \Re\left(h'_i(s)\right) Q_i(s)+Q'_i(x)\right)\right]  \right\vert \leq \exp\left(\sup_{x\in[-1,1]} 2 \Re(h_i(x)) \right) \tilde{\epsilon}
\end{align}
Thus
\begin{equation}
\left\vert q_i(x)e^{h_i(x)}  \right\vert \leq \sqrt{S(x)} \leq \left( \sup_{-1\leq x \leq 1} \tilde{S}(x) + \exp\left(\sup_{x\in[-1,1]} 2 \Re(h_i(x)) \right) \tilde{\epsilon} \right)^{1/2} := E_{i,1}
\end{equation}
In particular for $i =2$ we obtain upper bounds for the corresponding functions in  Group (II):
\begin{equation}\label{esgp2}
\left\vert q_i(x)e^{f_i(x)}  \right\vert \leq \left\vert q_i(x)e^{h_i(x)}  \right\vert \left\vert \exp \left(g_{i,0}\right)  \right\vert \exp \left(\overline{r}_i\right) \leq E_{i,1}\left\vert \exp \left(g_{i,0}\right)  \right\vert \exp \left(\overline{r}_i\right)
\end{equation}
Now let us get back to estimating $q_i R_{i,1}$, 
\begin{align}
\left\vert q_i R_{i,1}\right\vert < \left\vert e^{g_{i,0}}\right\vert E_{i,1} e^{\overline{r}_i}\overline{r}_i
\end{align}
To estimate $q_i R_{i,2}$ we write
\begin{align}
\begin{aligned}
q_i R_{i,2} = e^{g_{i,0}} q_i \left(\int_{0}^{h_i}\int_{0}^{s_1}\cdots \int_{0}^{s_4}e^{s_5} \mathrm{d}s_5 \cdots \mathrm{d}s_1\right) + e^{g_{i,0}} q_i\frac{h_i^4}{24}
\end{aligned}
\end{align}
and let
\begin{align}
& E_{i,3} :=  \left\vert e^{g_{i,0}} \right\vert \left\|q_i\right\|\frac{\left\|h_i\right\|^5}{120} \exp\left(\sup_{x\in[-1,1]}  \Re(h_i(x)) \right) \\
& e^{g_{i,0}} q_i\frac{h_i^4}{24} = \sum_{j} \tilde{l}_{j} T_j(x)\\
& E_{i,4} :=  \sum_{j}\left\vert \tilde{l}_{j} \right\vert
\end{align}
then
\begin{align}
\left\vert q_i R_{i,2}\right\vert \leq E_{i,3}+E_{i,4}
\end{align}
and
\begin{equation}
|G(x)|\leq \sum_{l=0}^{8}\left\vert \tilde{h}_l \right\vert +\sum_{i=1}^2 \left(\left\vert e^{g_{i,0}}\right\vert E_{i,1} e^{\overline{r}_i}\overline{r}_i+E_{i,3}+E_{i,4}\right)
\end{equation}
implying an upper bound of $|\tilde{G}|$ on $[a,b]$ in terms of $\tilde{h}_l$, $g_{i,0}$, $\overline{r}_i$ and $E_{i,j}$ ($0\le l \le 8$, $i=1,2$, $j = 1,3,4$).

\subsection{Proof of the estimates in Lemma \ref{L14} (iv)}\label{L4App}
\subsubsection{$U_1=1/30$}$ $\\
\indent On $[10\,i,380\,i]$, write $\lambda = t\,i$. Then $t$ is real and positive. We have
\begin{align}
\begin{aligned}
\left\vert \frac{A_n(t\,i)}{A_{n+1}(t\,i)}-1\right\vert^2 &=\frac{P_1(n,t)}{Q_1(n,t)}\\
\end{aligned}
\end{align}
where $P_1$, $Q_1$ are real-valued polynomials in $n$ and $t$, where
\begin{flalign}\label{P1Q1}
\begin{aligned}
& P_1(n,t) = 256\,{n}^{4}{t}^{2}+16\,{n}^{2}{t}^{4}+1888\,{n}^{3}{t}^{2}+72\,n{t}^{
	4}+1936\,{n}^{4}+5720\,{n}^{2}{t}^{2}\\
&\quad +81\,{t}^{4}+17600\,{n}^{3}+8600\,n{t}^{2}+61208\,{n}^{2}+5362\,{t}^{2}+96400\,n+58081
\end{aligned}\\
\begin{aligned}
& Q_1(n,t) = \left( 2\,n+5\right) ^{2} \left( n+1 \right) ^{2}\big( 144\,{n}^{4}+40\,{n}^{2}{t}^{2}+{t}^{4} +1056\,{n}^{3}+168\,n{t}^{2}\\
&\quad +2488\,{n}^{2} +210\,{t}^{2}+2024\,n+529 \big) 
\end{aligned}
\end{flalign}
Now, $P_1/Q_1$ is increasing in $t$ since,
\begin{align}
\begin{aligned}
\frac{\partial}{\partial t} \left\vert \frac{A_n(t\,i)}{A_{n+1}(t\,i)}-1\right\vert^2 = \frac{\partial}{\partial t} \frac{P_1(n,t)}{Q_1(n,t)} = \frac{P_2(n,t)}{Q_2(n,t)}
\end{aligned}
\end{align}
where, after simplification $P_2$ and $Q_2$ are manifestly positive polynomials whose explicit forms are not relevant.
Thus
\begin{align}\label{P1Q1380}
\begin{aligned}
&\frac{P_1(n,t)}{Q_1(n,t)}\leq \frac{P_1(n,380)}{Q_1(n,380)}\\
&={\frac {36968336\,{n}^{4}+272644800\,{n}^{3}+334447789208\,{n}^{2}+
		1502539856400\,n+1689734490881}{ \left( 144\,{n}^{4}+1056\,{n}^{3}+
		5778488\,{n}^{2}+24261224\,n+20881684529 \right)  \left( 2\,n+5
		\right) ^{2} \left( n+1 \right) ^{2}}}
\end{aligned}
\end{align}
The right side of \eqref{P1Q1380} is decreasing in $n$. Indeed
\begin{align}
\begin{aligned}
\frac{\mathrm{d}}{\mathrm{d} n} \frac{P_1(n,380)}{Q_1(n,380)} = -\frac{P_3(n)}{Q_3(n)}
\end{aligned}
\end{align}
where $P_3$ and $Q_3$ are manifestly positive for positive values of $n$. Hence for $n\geq 50$ and $10\leq t\leq 380$ we have
\begin{equation}
\left\vert \frac{A_n(t\,i)}{A_{n+1}(t\,i)}-1\right\vert^2 \leq \frac{P_1(n,t)}{Q_1(n,t)}\leq \frac{P_1(50,380)}{Q_1(50,380)} < U^2_1
\end{equation}
\\
\indent Similarly, on $[10\,i+1/2,380\,i+1/2]$, writing $\lambda = t\,i+1/2$, $\left\vert {A_n}/{A_{n+1}}-1\right\vert^2$ is monotonically increasing in $t$ and thus
\begin{equation}\label{U1v2}
\left\vert \frac{A_n(t\,i+1/2)}{A_{n+1}(t\,i+1/2)}-1\right\vert^2 \le \left\vert \frac{A_n(380\,i+1/2)}{A_{n+1}(380\,i+1/2)}-1\right\vert^2
\end{equation}
The right hand side of the equation \eqref{U1v2} is monotonically decreasing in $n$, so
\begin{equation}
\left\vert \frac{A_n(t\,i+1/2)}{A_{n+1}(t\,i+1/2)}-1\right\vert^2 \le \left\vert \frac{A_{50}(380\,i+1/2)}{A_{51}(380\,i+1/2)}-1\right\vert^2\leq U^2_1
\end{equation}
\\
\indent On $[10\,i,10\,i+1/2]$, write $\lambda = t+10\,i$, $\left\vert {A_n}/{A_{n+1}}-1\right\vert^2$ is monotonically decreasing in $t$ and thus
\begin{equation}\label{U1h1}
\left\vert \frac{A_n(t+10\,i)}{A_{n+1}(t+10\,i)}-1\right\vert^2 \le \left\vert \frac{A_n(10\,i)}{A_{n+1}(10\,i)}-1\right\vert^2
\end{equation}
The right hand side of the equation \eqref{U1h1} is monotonically decreasing in $n$, so
\begin{equation}
\left\vert \frac{A_n(t+10\,i)}{A_{n+1}(t+10\,i)}-1\right\vert^2 \le \left\vert \frac{A_{50}(10\,i)}{A_{51}(10\,i)}-1\right\vert^2\leq U^2_1
\end{equation}
\\
\indent Finally, on $[380\,i,380\,i+1/2]$, write $\lambda = t+380\,i$, $\left\vert {A_n}/{A_{n+1}}-1\right\vert^2$ is monotonically decreasing in $t$ and thus
\begin{equation}\label{U1h2}
\left\vert \frac{A_n(t+380\,i)}{A_{n+1}(t+380\,i)}-1\right\vert^2 \le \left\vert \frac{A_n(380\,i)}{A_{n+1}(380\,i)}-1\right\vert^2
\end{equation}
The right hand side of the equation \eqref{U1h2} is monotonically decreasing in $n$, so
\begin{equation}
\left\vert \frac{A_n(t+380\,i)}{A_{n+1}(t+380\,i)}-1\right\vert^2 \le \left\vert \frac{A_{50}(380\,i)}{A_{51}(380\,i)}-1\right\vert^2\leq U^2_1
\end{equation}

\subsubsection{$U_2=7 \sqrt{2}/500$}$ $\\
We now obtain an upper bound of $\left\vert F_{n}-F_{n+1}\right\vert$ for $n\ge 50$ and $\lambda \in \partial R_2$. Let $N = n-50$.\\
Consider $\lambda \in [10\,i,380\,i]$, write $\lambda = t\,i$, and denote $\Delta_N = F_{N+50}-F_{N+51}=R(N,t) + i S(N,t)$, where
\begin{align}
R(N,t) = \frac{P_4(N,t)}{Q_4(N,t)}\quad\quad\quad t\in[10,380],N \ge 0\\
S(N,t) = \frac{P_5(N,t)}{Q_4(N,t)} \quad\quad\quad t\in[10,380],N \ge 0
\end{align}
where $P_4$, $Q_4$ and $P_5$ are all real polynomials in $N$ and $t$. Explicitly
\begin{align}
\begin{aligned}
&Q_4(N,t)= \big( 144\,{N}^{4}+40\,{N}^{2}{t}^{2}+{t}^{4}+29280\,{N}^{3}+4088\,{
	t}^{2}N+2232184\,{N}^{2}\\
&+104482\,{t}^{2}+75618040\,N+960442081
\big) ^{2} \big( 144\,{N}^{4}+40\,{N}^{2}{t}^{2}+{t}^{4}+29856\,{N
}^{3}\\
&+4168\,{t}^{2}N+2320888\,{N}^{2}+108610\,{t}^{2}+80170824\,N+1038321729 \big) ^{2}
\end{aligned}
\end{align} 
is positive. We shall obtain bounds for $R(N,t)$ and $S(N,t)$ separately. Let $U_{2,R}=7/500$. Since $Q_4$ is obviously positive, proving $\left\vert R(N,t)\right\vert\leq U_{2,R}$ and $\left\vert S(N,t)\right\vert\leq U_{2,I}$ is equivalent to proving the following inequalities
\begin{align}
& -U_{2,R}Q_4(N,t) \le P_4(N,t)\le U_{2,R}Q_4(N,t)\leq 0\label{U3R}\\
& -U_{2,I}Q_4(N,t) \le P_5(N,t)\le U_{2,I}Q_4(N,t)\leq 0\label{U3I}
\end{align}
To prove \eqref{U3R} and \eqref{U3I} we write
\begin{align}
T_1(N,t)& =P_4(N,t)-U_{2,R}Q_4(N,t) = \sum_{i=0}^{16}w^{(1)}_i(t) N^i\\
T_2(N,t)& =P_4(N,t)+U_{2,R}Q_4(N,t) = \sum_{i=0}^{16}w^{(2)}_i(t) N^i\\
T_3(N,t)& =P_5(N,t)-U_{2,I}Q_4(N,t)=\sum_{i=0}^{16} w^{(3)}_i(t) N^i\\
T_4(N,t)&=P_5(N,t)+U_{2,I}Q_4(N,t)=\sum_{i=0}^{16} w^{(4)}_i(t) N^i
\end{align}

By our choice of $U_{2,R}$ and $U_{2,i}$, we can show the following,
\begin{align}
& w^{(j)}_i(t)\leq 0 \quad\quad (j=1,3 ,\, 0\le i\le 16, \, t\in [10,380])\label{w13}\\
& w^{(j)}_i(t)\geq 0 \quad\quad (j=2,4 ,\, 0\le i\le 16, \, t\in [10,380])\label{w24}
\end{align}
and consequently \eqref{U3R} and \eqref{U3I}. The details are in Section \ref{mono} in the Appendix.

Use the same method and still let $U_{2,R} = 7/500$ be an upper bound of $\left\vert \Re(\Delta_n) \right\vert$ and $U_{2,I}=7/500$ as an upper bound of $\left\vert \Im(\Delta_n) \right\vert$, then we are able to show that on the other three sides of $\partial R_2$
\begin{equation}
\left\vert F_{n}-F_{n+1}\right\vert \leq U_2
\end{equation}

\subsubsection{$L_1 = 8/25$}$ $\\
We would like to show now for $n\geq 50$ and $\lambda \in \partial R_2$
\begin{equation}
\Re\sqrt{F_{n}}> L_1:=8/25
\end{equation}
Let $u =\Re\sqrt{F_{n}}$ and $v =\Im\sqrt{F_{n}}$ we have
\begin{align}
\Re(F_{n}) & =u^2-v^2\\
\Im(F_{n}) & =2 u v
\end{align}
We will show that $\Re(F_{n}) \geq L_1^2$ and $\Im(F_n)>0$. Since we chose the branch of square root to be the one that is positive on the positive real axis, we have $$u>\sqrt{u^2-v^2} \geq L_1$$\\
\indent For $\lambda \in [10\,i,380\,i]$, let $\lambda=t\,i$ we can check straightforwardly that for all nonnegative $t$ and $n\geq 50$
\begin{align}\label{mono1}
& \frac{\partial}{\partial t} \Re\left(F_n(t\,i)\right) > 0\\
& \Im\left(F_n(t\,i)\right)>0
\end{align}
So we have
\begin{equation}
\Re\left(F_n(t\,i)\right) \geq \Re\left(F_n(0)\right) = {\frac {16\,{n}^{4}+32\,{n}^{3}+152\,{n}^{2}+648\,n+801}{ \left( 12\,{
			n}^{2}+20\,n-9 \right) ^{2}}}
\end{equation}
Then we can check that
\begin{equation}\label{mono2}
\frac{\mathrm{d}}{\mathrm{d} n} \Re\left(F_n(0)\right) >0
\end{equation}
and thus
\begin{equation}
\Re\left(F_n(t\,i)\right)\geq \Re\left(F_n(0)\right) \geq \Re\left(F_{50}(0)\right)>L_1^2
\end{equation}

\indent For $\lambda \in [10\,i+1/2,380\,i+1/2]$, let $\lambda=t\,i+1/2$ we can check straightforwardly that for all nonnegative $t$ and $n\geq 50$,
\begin{align}
& \frac{\partial}{\partial t} \Re\left(F_n(t\,i+1/2)\right) > 0\\
& \Im\left(F_n(t\,i+1/2)\right)>0
\end{align}
So we have
\begin{equation}
\Re\left(F_n(t\,i+1/2)\right) \geq \Re\left(F_n(1/2)\right) = {\frac {256\,{n}^{4}+1024\,{n}^{3}+3168\,{n}^{2}+9472\,n+11561}{
		\left( 48\,{n}^{2}+96\,n-19 \right) ^{2}}}
\end{equation}
Then we can check that
\begin{equation}
\frac{\mathrm{d}}{\mathrm{d} n} \Re\left(F_n(1/2)\right) < 0
\end{equation}
and thus
\begin{align}
& \Re\left(F_n(t\,i+1/2)\right)\geq \Re\left(F_n(1/2)\right) \geq \lim_{n \to \infty} \Re\left(F_{n}(1/2)\right)=1/9>L_1^2
\end{align}

\indent For $\lambda \in [10\,i,10\,i+1/2]$, let $\lambda=t+10\,i$ we can check straightforwardly that for all nonnegative $t$ and $n\geq 50$,
\begin{align}
& \frac{\partial}{\partial t} \Re\left(F_n(10\,i+t)\right) > 0\\
& \Im\left(F_n(10\,i+t)\right)>0
\end{align}
Thus, we have
\begin{equation}
\Re\left(F_n(10\,i+t)\right) \geq \Re\left(F_n(10\,i)\right)
\end{equation}
Then by \eqref{mono1} and \eqref{mono2} we obtain
\begin{equation}
\Re\left(F_n(10\,i)\right) \geq \Re\left(F_n(0)\right) \geq \Re\left(F_{50}(0)\right)>L_1^2
\end{equation}

\indent For $\lambda \in [380\,i,380\,i+1/2]$, let $\lambda=t+380\,i$ and write
\begin{align}
\Re(F_n) = \frac{P_6(n,t)}{Q_6(n,t)}
\end{align}
where
\begin{align}
\begin{aligned}
& Q_6(n,t) = \big( 144\,{n}^{4}+192\,{n}^{3}t+88\,{n}^{2}{t}^{2}+16\,n{t}^{3}+480\,{n}^{3}+512\,{n}^{2}t+168\,n{t}^{2}+16\,{t}^{3}\nonumber\\
& +{t}
^{4}+5776184\,{n}^{2}+2310576\,nt+288846\,{t}^{2}+12706840\,n+2310256\,t+20863200881
\big) ^{2} 
\end{aligned}
\end{align}
Defining $N=n-50$ then we can check directly that
\begin{align}
T_5(N,t) = P_6(N+50,t)-L_1^2 Q_6(N+50,t)
\end{align}
is positive for all nonnegative $N$ and $t$. Hence $\Re({F_n}) > L_1^2$ for  $\lambda \in [380\,i,380\,i+1/2]$ and all $n\ge 50$.

\subsubsection{$U_3=9/40$}$ $\\
We show that $U_3=9/40$ is an upper bound for $\left\vert B_n/\left(A_n A_{n+1}\right) \right\vert$ for $\lambda \in \partial R_2$. Denote
\begin{equation}
V_n(\lambda) = \left\vert\frac{B_n(\lambda)}{A_n(\lambda)A_{n+1}(\lambda)}\right\vert^2
\end{equation}
\indent For $\lambda \in [10\,i,380\,i]$, let $\lambda=t\,i$ then
\begin{equation}
V_n(t\,i) = \frac{P_7(n,t)}{Q_7(n,t)}
\end{equation}
where 
\begin{align}
\begin{aligned}
& Q_7(n,t) = \big( 144\,{n}^{4}+40\,{n}^{2}{t}^{2}+{t}^{4}+1056\,{n}^{3}+168\,n{t
}^{2}+2488\,{n}^{2}+210\,{t}^{2}+2024\,n\\
&+529 \big) \left( 144\,
{n}^{4}+40\,{n}^{2}{t}^{2}+{t}^{4}+480\,{n}^{3}+88\,n{t}^{2}+184\,{n}^
{2}+82\,{t}^{2}-360\,n+81 \right)
\end{aligned}
\end{align}
which is manifestly positive for all $t\in[10,380]$ and $n\ge 50$. Next we just check that for all nonnegative $N$ and $t$
\begin{equation}
P_7(N+50,t)-U_3^2 Q_7(N+50,t) \leq 0
\end{equation}
which is also obvious.

\indent For $\lambda \in [10\,i+1/2,380\,i+1/2]$, let $\lambda=t\,i+1/2$. Then
\begin{equation}
V_n(t\,i+1/2) = \frac{P_8(n,t)}{Q_8(n,t)}
\end{equation}
where 
\begin{align}
\begin{aligned}
& Q_8(n,t) = \big( 2304\,{n}^{4}+640\,{n}^{2}{t}^{2}+16\,{t}^{4}+18432\,{n}^{3}+
2816\,n{t}^{2}+48864\,{n}^{2}+3624\,{t}^{2}\\
&+48000\,n+15625 \big) 
\big( 2304\,{n}^{4}+640\,{n}^{2}{t}^{2}+16\,{t}^{4}+9216\,{n}^{3}+
1536\,n{t}^{2}+7392\,{n}^{2}\\
&+1448\,{t}^{2}-3648\,n+361 \big)
\end{aligned}
\end{align}
which is manifestly positive for all $t\in[10,380]$ and $n\ge 50$. Next we check that for all nonnegative $N$ and $t$
\begin{equation}
P_8(N+50,t)-U_3^2 Q_8(N+50,t) \leq 0
\end{equation}

\indent For $\lambda \in [10\,i, 10\,i+1/2]$ and $\lambda \in [380\,i, 380\,i+1/2]$, using a similar method we obtain that for $t\in [0,1/2]$ and $n\ge 50$
\begin{align}
& V_n(10\,i+t)\le U_3^2\\
& V_n(380\,i+t)\le U_3^2
\end{align}

\subsection{Proof of \eqref{w13} and \eqref{w24}}\label{mono}$ $\\
\indent We give the proof for $\lambda\in [10\,i,380\,i]$. The proof is similar for $\lambda\in [10\,i+1/2,380\,i+1/2]$. On the other two sides of $\partial R_2$ the sign of $w^{(j)}_i(t)$ is clear by inspection.\\
\subsubsection{$t\in[10,60]$}$ $\\
For each $0\le i\le 16$ and $j\in\{1,2,3,4\}$, assume that $w^{(j)}_i(t)$ is a polynomial of degree $\alpha$. If $\alpha \le 3$ it is trivial to calculate the maximum or minimum of $w^{(j)}_i(t)$ on $[10,60]$. If $\alpha > 3$, we express $w^{(j)}_i(t)$ as a linear combination of Chebyshev polynomials:
\begin{equation}
w^{(j)}_i(t) = \sum_{l=0}^{\alpha}\tilde{a}_l T_l \left(\frac{t-35}{25}\right)
\end{equation}
Using \eqref{Tle1} we have lower bounds and upper bounds for $w^{(j)}_i(t)$
\begin{align}\label{LwU}
\begin{aligned}
& \inf_{ t\in[10,60]}\left\{\sum_{l=0}^{3}\tilde{a}_l T_l\left(\frac{t-35}{25}\right)\right\} - \sum_{l=4}^{\alpha}\left\vert \tilde{a}_l \right\vert  \le w^{(j)}_i(t) \\
& \le  \sup_{t\in[10,60]}\left\{\sum_{l=0}^{3}\tilde{a}_l T_l\left(\frac{t-35}{25}\right)\right\} + \sum_{l=4}^{\alpha}\left\vert \tilde{a}_l \right\vert  
\end{aligned}
\end{align}
We can check that for $j=1,3$ the upper bound of $w^{(j)}_i(t)$ given by \eqref{LwU} is negative and for $j=2,4$ the lower bound. of $w^{(j)}_i(t)$ given by \eqref{LwU} is positive.\\
\subsubsection{$t\in[60,380]$}$ $\\
For each $j\in\{1,2,3,4\}$, assume that $K(t) = w^{(j)}_i(t)$ is a polynomial of degree $\alpha$, then $K(t)$ has the property that the ${\beta}$-th derivative of $K(t)$ is monotonic and does not change sign on $[60,380]$ for each $0\le \beta \le \alpha$. We can show this by induction.

It is obvious that the $\alpha$-th derivative of $K(t)$, which is a constant, does not change sign on any interval. Since the ${\beta}$-th derivative $K^{(\beta)}(t)$ does not change sign on $[60,380]$ we have for each $\beta$

\begin{equation}\label{endsign}
\frac{\mathrm{d}^{\beta}K(60)}{\mathrm{d}t^{\beta}}\frac{\mathrm{d}^{\beta}K(380)}{\mathrm{d}t^{\beta}} \ge 0
\end{equation}

Assume that $j\le 0$. Suppose the $(\alpha-j)$-th derivative does not change sign on $[60,380]$ and \eqref{endsign} holds for $\beta = \alpha-j-1$. Then the $(\alpha-j-1)$-th derivative is monotonic and does not change sign on $[60,380]$. By induction on $j$ we see that $K(t)$ is monotonic and does not change sign on $[60,380]$. Now it is trivial to check that

\begin{align}
w_i^{(j)}(60)<0 \quad \quad w_i^{(j)}(380)<0 \quad \quad(j=1,3)\\
w_i^{(j)}(60)>0 \quad \quad w_i^{(j)}(380)>0 \quad \quad(j=2,4)
\end{align}

Hence each $w_i^{(j)}(t)$ does not change sign on the interval $[60,380]$ and \eqref{w13} and \eqref{w24} are proved.

\newpage

\subsection{Table of coefficients $a^{(1,2)}_{n,i}$ and $b^{(1,2)}_{n,i}$}\label{table} $ $\\
\it{Note: In the table below, a pair of the form ``$A,B$'' is used to denote the complex number $A+B\,i$}\\

{\renewcommand{\arraystretch}{1.3} \renewcommand\tabcolsep{1pt}
	\hspace{0cm}	
\begin{tabular}{|l||l|l|l|l|l|l||}
		\hline		i		  & 0 & 1 & 2 & 3 & 4 & 5 \\ \hline
		$a^{(1)}_{1,i}$ & ${\frac {49}{12}},{\frac {43}{15}}$ & ${\frac {17}{10}},{\frac {31}{139}}$ & $
		{\frac {-9}{26}},{\frac {-7}{80}}$ & ${\frac {10}{111}},{\frac {2}{61}}$ & ${
			\frac {-3}{118}},{\frac {-1}{86}}$ & ${\frac {1}{135}},{\frac {3}{784}}
		$ \\ \hline
		$a^{(2)}_{1,i}$ & ${\frac {45}{11}},{\frac {17}{6}}$ & ${\frac {22}{13}},\frac{1}{4}$ & ${\frac {-13}{38
			}},{\frac {-9}{92}}$ & ${\frac {3}{34}},{\frac {3}{82}}$ & ${\frac {-2}{81}},
			{\frac {-4}{307}}$ & ${\frac {2}{279}},{\frac {2}{461}}
			$ \\ \hline
			$a^{(1)}_{2,i}$ & ${\frac {63}{20}},{\frac {-40}{11}}$ & ${\frac {19}{12}},{\frac {13}{36}}$ & $
			{\frac {-26}{89}},{\frac {-8}{65}}$ & ${\frac {14}{201}},{\frac {9}{230}}$ & $
			{\frac {-6}{313}},{\frac {-4}{343}}$ & ${\frac {2}{327}},{\frac {2}{605}}
			$ \\ \hline
			$a^{(2)}_{2,i}$ & ${\frac {60}{19}},\frac{-11}{3}$ & ${\frac {96}{61}},{\frac {18}{47}}$ & ${\frac {-25}
				{87}},{\frac {-8}{61}}$ & ${\frac {4}{59}},{\frac {3}{71}}$ & ${\frac {-6}{325
				}},{\frac {-5}{387}}$ & ${\frac {1}{171}},{\frac {2}{517}}
				$ \\ \hline
				$a^{(1)}_{3,i}$ & ${\frac {96}{37}},{\frac {42}{17}}$ & ${\frac {34}{23}},{\frac {13}{29}}$ & $
				{\frac {-14}{55}},{\frac {-9}{64}}$ & ${\frac {8}{135}},{\frac {6}{139}}$ & ${
					\frac {-6}{359}},{\frac {-6}{433}}$ & ${\frac {2}{371}},{\frac {1}{205}}
				$ \\ \hline
				$a^{(2)}_{3,i}$ & ${\frac {99}{38}},{\frac {105}{43}}$ & ${\frac {22}{15}},{\frac {22}{47}}$ & $
				{\frac {-58}{233}},{\frac {-14}{95}}$ & ${\frac {5}{88}},{\frac {6}{131}}$ & $
				{\frac {-1}{64}},{\frac {-4}{269}}$ & ${\frac {2}{409}},{\frac {2}{379}}
				$ \\ \hline
				$a^{(1)}_{4,i}$ & ${\frac {42}{19}},{\frac {44}{19}}$ & ${\frac {18}{13}},{\frac {22}{43}}$ & $
				{\frac {-19}{85}},{\frac {-11}{72}}$ & ${\frac {6}{121}},{\frac {9}{190}}$ & $
				{\frac {-8}{615}},{\frac {-3}{187}}$ & ${\frac {1}{276}},{\frac {1}{170}}
				$\\
\hline
							\end{tabular}\begin{tabular}{||l || l | l | l | l | l | l |}
							\hline		i		  & 0 & 1 & 2 & 3 & 4 & 5 \\ \hline
							$a^{(1)}_{17,i}$ & ${\frac {11}{16}},\frac{4}{3}$ & ${\frac {67}{101}},{\frac {13}{20}}$ & ${\frac {-4}{
									197}},{\frac {-3}{25}}$ & ${\frac {-2}{233}},{\frac {9}{451}}$ & ${\frac {2}{
									609}},{\frac {-1}{420}}$ & $0,0
							$\\ \hline
							$a^{(2)}_{17,i}$ & ${\frac {51}{73}},{\frac {41}{31}}$ & ${\frac {19}{29}},{\frac {28}{43}}$ & $
							{\frac {-5}{273}},{\frac {-5}{42}}$ & ${\frac {-3}{337}},{\frac {4}{205}}$ & ${
								\frac {2}{605}},{\frac {-1}{440}}$ & $0,0
							$\\ \hline
							$a^{(1)}_{18,i}$ & ${\frac {31}{48}},{\frac {22}{17}}$ & ${\frac {43}{68}},{\frac {38}{59}}$ & $
							{\frac {-3}{199}},{\frac {-4}{35}}$ & ${\frac {-3}{343}},{\frac {5}{279}}$ & ${
								\frac {1}{336}},{\frac {-1}{511}}$ & $0,0
							$\\ \hline
							$a^{(2)}_{18,i}$ & ${\frac {23}{35}},{\frac {59}{46}}$ & $\frac{5}{8},{\frac {20}{31}}$ & ${\frac {-4}{
									301}},{\frac {-11}{97}}$ & ${\frac {-1}{111}},{\frac {9}{514}}$ & ${\frac {1}{
									334}},{\frac {-1}{536}}$ & $0,0
							$\\ \hline
							$a^{(1)}_{19,i}$ & ${\frac {14}{23}},{\frac {64}{51}}$ & ${\frac {35}{58}},{\frac {58}{91}}$ & $
							{\frac {-3}{283}},{\frac {-16}{147}}$ & ${\frac {-3}{343}},{\frac {15}{931}
							}$ & ${\frac {1}{373}},{\frac {-1}{616}}$ & $0,0
							$\\ \hline
							$a^{(2)}_{19,i}$ & ${\frac {31}{50}},{\frac {61}{49}}$ & ${\frac {28}{47}},{\frac {44}{69}}$ & $
							{\frac {-5}{556}},{\frac {-11}{102}}$ & ${\frac {-2}{223}},{\frac {2}{127}}
							$ & ${\frac {1}{370}},{\frac {-1}{646}}$ & $0,0
							$\\ \hline
							$a^{(1)}_{20,i}$ & ${\frac {27}{47}},{\frac {39}{32}}$ & ${\frac {34}{59}},{\frac {17}{27}}$ & $
							{\frac {-3}{443}},{\frac {-26}{251}}$ & ${\frac {-6}{695}},{\frac {8}{551}}
							$ & ${\frac {1}{415}},{\frac {-1}{732}}$ & $0,0
							$ \\
\hline
						\end{tabular}

						\hspace{0cm} \begin{tabular}{|l || l | l | l | l | l | l ||}
							\hline		i		  & 0 & 1 & 2 & 3 & 4 & 5 \\ \hline
							$a^{(1)}_{33,i}$ & ${\frac {29}{94}},{\frac {47}{53}}$ & ${\frac {21}{61}},{\frac {19}{36}}$ & ${
								\frac {3}{242}},{\frac {-7}{123}}$ & ${\frac {-1}{197}},{\frac {3}{601}}$ & $0
							,0$ & $0,0
							$ \\ \hline
							$a^{(2)}_{33,i}$ & ${\frac {20}{63}},{\frac {15}{17}}$ & ${\frac {16}{47}},{\frac {10}{19}}$ & ${
								\frac {3}{232}},{\frac {-4}{71}}$ & ${\frac {-3}{581}},{\frac {2}{411}}$ & $0,0
							$ & $0,0
							$ \\ \hline
							$a^{(1)}_{34,i}$ & ${\frac {8}{27}},{\frac {33}{38}}$ & ${\frac {144}{433}},{\frac {13}{25}}$ & $
							{\frac {2}{157}},{\frac {-7}{128}}$ & ${\frac {-2}{411}},{\frac {2}{425}}$ & $0
							,0$ & $0,0
							$ \\ \hline
							$a^{(2)}_{34,i}$ & ${\frac {18}{59}},{\frac {19}{22}}$ & ${\frac {24}{73}},{\frac {14}{27}}$ & ${
								\frac {2}{151}},{\frac {-15}{277}}$ & ${\frac {-4}{809}},{\frac {1}{218}}$ & $0
							,0$ & $0,0
							$ \\ \hline
							$a^{(1)}_{35,i}$ & ${\frac {39}{137}},{\frac {23}{27}}$ & ${\frac {9}{28}},{\frac {20}{39}}$ & ${
								\frac {9}{692}},{\frac {-27}{514}}$ & ${\frac {-2}{429}},{\frac {1}{225}}$ & $0
							,0$ & $0,0
							$ \\ \hline
							$a^{(2)}_{35,i}$ & ${\frac {22}{75}},{\frac {11}{13}}$ & ${\frac {7}{22}},{\frac {22}{43}}$ & ${
								\frac {4}{297}},{\frac {-9}{173}}$ & ${\frac {-1}{211}},{\frac {1}{231}}$ & $0
							,0$ & $0,0
							$ \\ \hline
							$a^{(1)}_{36,i}$ & ${\frac {20}{73}},{\frac {151}{181}}$ & ${\frac {14}{45}},{\frac {46}{91}}
							$ & ${\frac {3}{227}},{\frac {-5}{99}}$ & ${\frac {-2}{447}},{\frac {3}{712}}$ & $0
							,0$ & $0,0
							$ \\
\hline
						\end{tabular}\begin{tabular}{||l || l | l | l | l | l | l |}
						\hline		i		  & 0 & 1 & 2 & 3 & 4 & 5 \\ \hline
						$b^{(1)}_{1,i}$ & ${\frac {83}{11}},{\frac {34}{11}}$ & ${\frac {59}{37}},{\frac {3}{77}}$ & ${\frac {-13}{41}},{\frac {-2}{129}}$ & ${\frac {9}{107}},{\frac {2}{325}}$ & ${
							\frac {-22}{879}},{\frac {-1}{409}}$ & ${\frac {1}{126}},0
						$ \\ \hline
						$b^{(2)}_{1,i}$ & ${\frac {83}{11}},{\frac {71}{23}}$ & ${\frac {51}{32}},{\frac {5}{114}}$ & $
						{\frac {-13}{41}},{\frac {-3}{172}}$ & ${\frac {10}{119}},{\frac {2}{289}}$ & $
						\frac{-1}{40},{\frac {-1}{364}}$ & ${\frac {3}{379}},{\frac {1}{918}}
						$ \\ \hline
						$b^{(1)}_{2,i}$ & ${\frac {189}{29}},{\frac {-68}{21}}$ & ${\frac {27}{17}},{\frac {28}{363}
						}$ & ${\frac {-11}{35}},{\frac {-6}{197}}$ & ${\frac {7}{85}},{\frac {2}{167}}
						$ & ${\frac {-4}{165}},{\frac {-3}{638}}$ & ${\frac {3}{397}},{\frac {1}{544}}
						$ \\ \hline
						$b^{(2)}_{2,i}$ & ${\frac {176}{27}},{\frac {-159}{49}}$ & ${\frac {27}{17}},{\frac {5}{61}}
						$ & ${\frac {-16}{51}},{\frac {-14}{433}}$ & ${\frac {6}{73}},{\frac {3}{236}}
						$ & ${\frac {-5}{207}},{\frac {-2}{401}}$ & ${\frac {1}{133}},{\frac {1}{513}}
						$ \\ \hline
						$b^{(1)}_{3,i}$ & ${\frac {88}{15}},{\frac {920}{307}}$ & ${\frac {30}{19}},{\frac {5}{44}}$ & $
						{\frac {-17}{55}},{\frac {-9}{203}}$ & ${\frac {23}{288}},{\frac {4}{233}}
						$ & ${\frac {-7}{304}},{\frac {-2}{303}}$ & ${\frac {2}{285}},{\frac {1}{396}}
						$ \\ \hline
						$b^{(2)}_{3,i}$ & ${\frac {88}{15}},{\frac {323}{108}}$ & ${\frac {30}{19}},{\frac {9}{76}}$ & $
						{\frac {-25}{81}},{\frac {-3}{65}}$ & ${\frac {9}{113}},{\frac {1}{56}}$ & ${
							\frac {-3}{131}},{\frac {-2}{291}}$ & ${\frac {2}{287}},{\frac {1}{381}}
						$ \\
\hline
				\end{tabular}
				
				\hspace{0cm}\begin{tabular}{|l || l | l | l | l | l | l ||}
					\hline		i		  & 0 & 1 & 2 & 3 & 4 & 5 \\ \hline
					$b^{(1)}_{19,i}$ &  ${\frac {146}{53}},{\frac {46}{19}}$ &  ${\frac {38}{29}},{\frac {25}{57}}$ &  $
					{\frac {-24}{121}},{\frac {-17}{132}}$ &  ${\frac {2}{51}},{\frac {5}{131}}
					$ &  ${\frac {-3}{344}},{\frac {-5}{423}}$ &  ${\frac {1}{498}},{\frac {1}{260}}
					$ \\ \hline
					$b^{(2)}_{19,i}$ & ${\frac {102}{37}},{\frac {29}{12}}$ &  ${\frac {72}{55}},{\frac {15}{34}}$ &  $
					{\frac {-15}{76}},{\frac {-7}{54}}$ &  ${\frac {4}{103}},\frac{1}{26}$ &  ${\frac {-3}{
							350}},{\frac {-7}{589}}$ &  ${\frac {1}{513}},{\frac {1}{259}}
					$ \\ \hline
					$b^{(1)}_{20,i}$ & ${\frac {123}{46}},{\frac {67}{28}}$ & ${\frac {22}{17}},{\frac {22}{49}}$ & $
					{\frac {-17}{88}},{\frac {-3}{23}}$ & ${\frac {7}{186}},{\frac {33}{859}}$ & $
					{\frac {-1}{122}},{\frac {-3}{253}}$ & ${\frac {1}{552}},{\frac {2}{519}}
					$ \\ \hline
					$b^{(2)}_{20,i}$ & ${\frac {91}{34}},{\frac {43}{18}}$ & ${\frac {22}{17}},{\frac {14}{31}}$ & $
					{\frac {-5}{26}},{\frac {-8}{61}}$ & ${\frac {6}{161}},{\frac {8}{207}}$ & ${
						\frac {-5}{621}},{\frac {-6}{503}}$ & ${\frac {1}{570}},{\frac {1}{259}}
					$ \\ \hline
					$b^{(1)}_{21,i}$ & ${\frac {13}{5}},{\frac {-47}{12}}$ & ${\frac {32}{25}},{\frac {11}{24}}$ & $
					{\frac {-13}{69}},{\frac {-5}{38}}$ & ${\frac {3}{83}},{\frac {10}{259}}$ & ${
						\frac {-1}{130}},{\frac {-1}{84}}$ & ${\frac {1}{618}},{\frac {1}{260}}
					$ \\ \hline
					$b^{(2)}_{21,i}$ & ${\frac {13}{5}},{\frac {-47}{12}}$ & ${\frac {23}{18}},{\frac {6}{13}}$ & $\frac{-3}
					{16},{\frac {-9}{68}}$ & ${\frac {30}{839}},{\frac {4}{103}}$ & ${\frac {-2}{
							265}},{\frac {-3}{251}}$ & ${\frac {1}{640}},{\frac {2}{519}}
					$ \\
\hline
					\end{tabular}\begin{tabular}{||l || l | l | l | l | l | l |}
					\hline		i		  & 0 & 1 & 2 & 3 & 4 & 5 \\ \hline
					$b^{(1)}_{35,i}$ & ${\frac {63}{34}},{\frac {35}{17}}$ & ${\frac {37}{34}},{\frac {17}{31}}$ & $
					{\frac {-11}{84}},{\frac {-9}{64}}$ & ${\frac {4}{219}},{\frac {7}{183}}$ & ${
						\frac {-1}{624}},{\frac {-6}{557}}$ & $0,{\frac {2}{673}}
					$ \\ \hline
					$b^{(2)}_{35,i}$ & ${\frac {13}{7}},{\frac {37}{18}}$ & ${\frac {25}{23}},{\frac {11}{20}}$ & ${
						\frac {-13}{100}},{\frac {-10}{71}}$ & ${\frac {3}{167}},{\frac {8}{209}}$ & $
					{\frac {-1}{667}},{\frac {-1}{93}}$ & $0,{\frac {1}{339}}
					$ \\ \hline
					$b^{(1)}_{36,i}$ & ${\frac {49}{27}},{\frac {51}{25}}$ & ${\frac {14}{13}},{\frac {21}{38}}$ & $
					{\frac {-13}{102}},{\frac {-9}{64}}$ & ${\frac {8}{465}},{\frac {7}{184}}$ & $
					{\frac {-1}{789}},{\frac {-2}{189}}$ & $0,{\frac {1}{348}}
					$ \\ \hline
					$b^{(2)}_{36,i}$ & ${\frac {20}{11}},{\frac {57}{28}}$ & ${\frac {29}{27}},{\frac {46}{83}}$ & $
					{\frac {-10}{79}},{\frac {-11}{78}}$ & ${\frac {6}{355}},{\frac {7}{184}}$ & $
					{\frac {-1}{856}},{\frac {-3}{284}}$ & $0,{\frac {2}{701}}$
					\\ \hline
					$b^{(1)}_{37,i}$ & ${\frac {16}{9}},{\frac {97}{48}}$ & ${\frac {33}{31}},{\frac {59}{106}}$ & $
					{\frac {-16}{129}},{\frac {-9}{64}}$ & ${\frac {5}{309}},{\frac {17}{450}}$ & $0
					,{\frac {-1}{96}}$ & $0,{\frac {2}{721}}
					$ \\ \hline
					$b^{(2)}_{37,i}$ & ${\frac {41}{23}},{\frac {119}{59}}$ & ${\frac {17}{16}},{\frac {24}{43}}$ & $
					{\frac {-9}{73}},{\frac {-11}{78}}$ & ${\frac {1}{63}},{\frac {17}{450}}$ & $0
					,{\frac {-3}{289}}$ & $0,{\frac {1}{363}}
					$ \\ \hline
					$b^{(1)}_{38,i}$ & ${\frac {68}{39}},{\frac {793}{396}}$ & ${\frac {20}{19}},{\frac {14}{25}}
					$ & ${\frac {-7}{58}},{\frac {-9}{64}}$ & ${\frac {13}{857}},{\frac {3}{80}}$ & $0
					,{\frac {-10}{979}}$ & $0,{\frac {1}{374}}
					$ \\
\hline
				\end{tabular}

\section{Acknowledgments}
{\rm The work of O.C. was partly supported by the NSF grant DMS 1108794.}


\begin{thebibliography}{99}

\bibitem{BIZON} Piotr  Bizo\'n. 
\newblock An unusual eigenvalue problem. 
\newblock {\em Acta Phys. Polon. B}, 36(1):5-15, 2005.

\bibitem{CHS} Ovidiu Costin, Min Huang and Wilhelm Schlag.
\newblock On the spectral properties of $L_{\pm}$ in three dimensions.
\newblock Nonlinearity 25, pp. 125--164, 2012.


\bibitem{Elaydi} Saber Elaydi 
\newblock {\em An Introduction to Difference Equations.}
\newblock Springer, New York, 2005.

\bibitem{Lor92}
Lisa Lorentzen and Haakon Wadeland.
\newblock {\em Continued fractions with applications}.
\newblock North-Holland Publishing Co., Amsterdam, 1992.

\bibitem{Don11}
Roland Donninger.
\newblock On stable self-similar blowup for equivariant wave maps.
\newblock {\em Comm. Pure Appl. Math.}, 64(8):1095--1147, 2011.

\bibitem{DonSchAic11}
Roland Donninger, Birgit Sch\"orkhuber, and Peter Aichelburg.
\newblock On stable self-similar blow up for equivariant wave maps: The
  linearized problem.
\newblock {\em Ann.~Henri Poincar\'e}, 13:103--144, 2012.

\bibitem{BizChmTab00}
Piotr Bizo{\'n}, Tadeusz Chmaj, and Zbis{\l}aw Tabor.
\newblock Dispersion and collapse of wave maps.
\newblock {\em Nonlinearity}, 13(4):1411--1423, 2000.

\bibitem{CazShaTah98}
Thierry Cazenave, Jalal Shatah, and A.~Shadi Tahvildar-Zadeh.
\newblock Harmonic maps of the hyperbolic space and development of
  singularities in wave maps and {Y}ang-{M}ills fields.
\newblock {\em Ann. Inst. H. Poincar\'e Phys. Th\'eor.}, 68(3):315--349, 1998.

\bibitem{DonAic08}
Roland Donninger and Peter~C. Aichelburg.
\newblock On the mode stability of a self-similar wave map.
\newblock {\em J. Math. Phys.}, 49(4):043515, 9, 2008.

\bibitem{DonAic09}
Roland Donninger and Peter~C. Aichelburg.
\newblock Spectral properties and linear stability of self-similar wave maps.
\newblock {\em J. Hyperbolic Differ. Equ.}, 6(2):359--370, 2009.

\bibitem{DonAic10}
Roland Donninger and Peter~C. Aichelburg.
\newblock A note on the eigenvalues for equivariant maps of the {SU(2)}
  sigma-model.
\newblock {\em Appl.~Math.~Comp.~Sciences}, 1(1):73--82, 2010.

\bibitem{GelLev60}
M.~Gell-Mann and M.~L{\'e}vy.
\newblock The axial vector current in beta decay.
\newblock {\em Nuovo Cimento (10)}, 16:705--726, 1960.

\bibitem{Kri08}
J.~Krieger.
\newblock Global regularity and singularity development for wave maps.
\newblock In {\em Surveys in differential geometry. {V}ol. {XII}. {G}eometric
  flows}, volume~12 of {\em Surv. Differ. Geom.}, pages 167--201. Int. Press,
  Somerville, MA, 2008.

\bibitem{Mis78}
Charles~W. Misner.
\newblock Harmonic maps as models for physical theories.
\newblock {\em Phys. Rev. D (3)}, 18(12):4510--4524, 1978.

\bibitem{Sha88}
Jalal Shatah.
\newblock Weak solutions and development of singularities of the {${\rm
  SU}(2)$} {$\sigma$}-model.
\newblock {\em Comm. Pure Appl. Math.}, 41(4):459--469, 1988.

\bibitem{TurSpe90}
Neil Turok and David Spergel.
\newblock Global texture and the microwave background.
\newblock {\em Phys. Rev. Lett.}, 64(23):2736--2739, Jun 1990.

\end{thebibliography}
\end{document}